\documentclass[twoside,10pt]{book}
\usepackage{amsmath,amsthm}
\usepackage{graphicx,amscd,amssymb}
\graphicspath{{figures/}{fr/}}
\usepackage[all]{xy}

\newcommand{\ass}{\mathcal{A}ss}
\newcommand{\Aut}{\operatorname{Aut}}
\newcommand{\Beta}{\operatorname{B}}
\newcommand{\CC}{\mathcal{C}}
\newcommand{\Char}{\operatorname{char}}
\newcommand{\comm}{\mathcal{C}omm}
\newcommand{\Corr}{\mathcal{C}orr}
\newcommand{\DD}{\mathcal{D}}

\newcommand{\End}[1]{\mathcal{E}nd\,_{#1}}
\newcommand{\ev}{\operatorname{ev}}
\newcommand{\Hom}{\operatorname{Hom}}
\newcommand{\id}{\operatorname{id}}
\newcommand{\II}{\mathcal{I}}
\newcommand{\In}{\operatorname{In}}

\newcommand{\lie}{\mathcal{L}ie}
\newcommand{\Map}{\operatorname{Map}}
\newcommand{\MM}{\mathcal{M}}
\newcommand{\Mor}{\operatorname{Mor}}
\newcommand{\OO}{\mathcal{O}}
\newcommand{\out}{\operatorname{out}}
\newcommand{\PB}{\operatorname{PB}}
\newcommand{\PP}{\mathcal{P}}
\newcommand{\PRB}{\operatorname{PRB}}

\newcommand{\RB}{\operatorname{RB}}
\newcommand{\Si}{\Sigma}
\newcommand{\SO}{\operatorname{SO}}
\newcommand{\tr}{\operatorname{tr}}
\newcommand{\UU}{\mathcal{U}}
\newcommand{\varin}{\operatorname{in}}
\newcommand{\virt}{{\operatorname{virt}}}

\newcommand{\nc}{{\mathbb{C}}}
\newcommand{\nq}{{\mathbb{Q}}}
\newcommand{\nr}{{\mathbb{R}}}
\newcommand{\nz}{{\mathbb{Z}}}

\newcommand{\bc}{\mathbb{C}}
\newcommand{\bp}{\mathbb{ P}}
\newcommand{\bz}{\mathbb{Z}}
\newcommand{\br}{\mathbb{R}}

\newcommand{\bh}{\mathbb{H}}

\newcommand{\p}{\partial}
\newcommand{\cc}{\mathcal{C}}

\newcommand{\cs}{\mathcal{S}}

\newcommand{\cp}{\mathcal{P}}

\newcommand{\cm}{\mathcal{M}}

\newcommand{\sr}{\mathcal{R}}

\newcommand{\cf}{\mathcal{F}}

\newcommand{\ca}{\mathcal{A}}
\newcommand{\hk}{\hookrightarrow}
\newcommand{\bg}{\bigskip}
\newcommand{\med}{\medskip}
\newcommand{\la}{\longrightarrow}
\newcommand{\bfl}{\begin{flushleft}}
\newcommand{\efl}{\end{flushleft}} 
\newcommand{\eps}{\epsilon}
\newcommand{\mtm}{M^{-TM}}
\newcommand{\ltm}{LM^{-TM}}

\newcommand{\om}{\Omega}

\newcommand{\xr}{\xrightarrow}

\newcommand{\met}{Map(8,M)}
\newcommand{\heq}{H^{S^1}_*(LM)}
\newcommand{\G}{\Gamma}
\newcommand{\scrb}{\mathcal{B}}
\newtheorem{thm}{Theorem}[section]
\newtheorem{lm}[thm]{Lemma}
\newtheorem{prop}[thm]{Proposition}
\newtheorem{crl}[thm]{Corollary}
 
\newtheorem{conj}[thm]{Conjecture}

\theoremstyle{definition}

\newtheorem{definition}{Definition}[section]
\newtheorem{ex}{Example}[section]
\newtheorem{xca}{Exercise}

\theoremstyle{remark}

\newtheorem{rem}{Remark}%[section]
\newtheorem{ack}{Acknowledgments}

\newcommand{\abs}[1]{\lvert#1\rvert}
\begin{document}

\frontmatter

\title{Notes on String Topology}

\author{Ralph L. Cohen
\thanks{Partially supported by a grant from the NSF.}
\\
%\address{
Department of Mathematics\\
Stanford University\\
Building 380\\
Stanford, CA 94305-2125\\
%\email
\texttt{ralph@math.stanford.edu}
%\author{
\and Alexander A. Voronov $^1$
\\
%\address{
School of Mathematics\\
University of Minnesota\\
206 Church St.\ S.E.\\
Minneapolis, MN 55455\\
%\email
\texttt{voronov@umn.edu}}

\date{March 23, 2005}

%\begin{abstract}
%\end{abstract}

\setcounter{page}{5}
\maketitle
\tableofcontents
\thispagestyle{empty}

\mainmatter
%version 3/21/05

 \chapter*{Introduction}
 
 String topology is the study of algebraic and differential topological properties of spaces of paths and loops in manifolds.
 It was initiated by the beautiful paper of Chas and Sullivan \cite{chas-sullivan} in which algebraic structures in  both the nonequivariant and equivariant    homology (and indeed chains)  of the (free) loop space, $LM$,  of a  closed, oriented manifolds were uncovered.  
 This has lead to considerable work by many authors over the past    five years.  The   goals  of this paper are twofold.  First, this paper is meant to be an introduction to this new and exciting field.  Second, we will attempt to give a ``status report".  That is, we will describe what has been learned over the last few years, and also give our views about future directions of research.
 This paper is a joint account of each of the author's lecture series given at the 2003 Summer School on String Topology and Hochschild Homology, in Almer\'{\i}a, Spain. 
 
   \med
   In our view there are two basic reasons for the excitement about the development of string topology.  First, it uses most of the modern techniques of algebraic topology, and relates them to several other areas of mathematics.  For example, the description of the structure involved in the string topology operations uses such concepts as operads, PROPs, field theories, and Gerstenhaber and Batalin-Vilkovisky algebras.  The fundamental  role played by   moduli spaces of Riemann surfaces in string topology, relates it  to basic objects of study in algebraic and symplectic geometry.   Techniques in low dimensional topology  such as the use of graphs to  study these  moduli spaces are also used in an essential way.   
  Moreover there are both formal and computational relationships between string topology and Gromov-Witten theory that are only beginning to be uncovered.  Gromov-Witten theory is a basic tool in string theory, algebraic geometry, and symplectic geometry, and  understanding its relationship to string topology is an exciting area of current and probably future research.
  
  The second reason for the attention the development of string topology has been receiving has to do with the historical significance,  in both mathematics and physics, played by spaces of paths and loops in manifolds.  The systematic study of the differential topology of path and loop spaces began in the 1930's with Morse, who used his newly developed theory of  ``calculus of variations in the large" to prove among other things that for any Riemannian metric on the $n$-sphere, there are an infinite number of geodesics connecting any two points. 
  In the 1950's, R. Bott studied Morse theory on the loop spaces of Lie groups and symmetric spaces to prove his celebrated periodicity theorem.  In  the 1970's and 1980's, the   $K$-theoretic tools developed by Waldhausen to study diffeomorphisms of high dimensional manifolds were found to be closely related to the equivariant stable homotopy type of the free loop space.  Finally, within the development of string theory in physics,  the basic configuration spaces are spaces of paths and loops in a manifold. Some of the topological issues
  this theory has raised are the following.

  \begin{enumerate}
  \item What mathematical structure should the appropriate notions of field and field strength have in this theory?  This has been addressed by the notion of a ``B-field", or a ``gerbe with connection".  These are structures on principal bundles over the loop space.
  
  \item How does one view elliptic operators, such as the Dirac operator, on the loop space of a manifold?  The corresponding index theory has been developed in the context of elliptic cohomology theory.
  
  \item  How does one understand geometrically and topologically, intersection theory in the infinite dimensional loop space and path space of a manifold?
  \end{enumerate}
  
  It is this last  question  that is the subject of string topology.
 The goal of these notes is to give an introduction to the exciting developments in this new theory.  They are organized as follows. 
 In Chapter 1 we review basic intersection theory, including the Thom-Pontrjagin construction, for compact manifolds.  We then
 develop and review the results and constructions of Chas and Sullivan's original paper.  In Chapter 2 we review the concepts of operads and PROPS, discuss many examples,  and study in detail the important example of the ``cacti  operad", which plays a central role in string topology.  In Chapter 3 we discuss field theories in general, and the field theoretic properties of string topology.  Included  are discussions of ``fat graphs", and how they give a model for the moduli space of Riemann surfaces, and 
 of ``open-closed" string topology, which involves spaces of paths in a manifold with prescribed boundary conditions. In Chapter 4 we discuss a Morse theoretic interpretation of string topology, incorporating the classical energy functional on the loop space, originally studied by Morse himself.  In this chapter we also discuss how this perspective suggests a potentially deep relationship  with the Gromov-Witten theory of the cotangent bundle.  Finally in Chapter 5 we study similar structures on spaces of maps of higher dimensional spheres to manifolds.

\begin{ack}
We are very grateful to David Chataur, Jos\'e Luis Rodr\'{\i}guez, and
J\'er\^ome Scherer for organizing and inviting us to participate in
such an active and inspiring summer school.
\end{ack}

%%% Version of 03/17/05

  \chapter{Intersection theory in loop spaces}
 String topology is ultimately about the differential and algebraic topology of spaces of paths and loops in compact,   oriented manifolds.  The basic spaces of paths that we consider are $C^\infty (\br , M)$, $C^\infty ([0,1], M)$, which we denote by $\cp (M)$, $C^\infty (S^1, M)$, which we denote by $LM $, and 
  $\om (M, x_0) = \{\alpha \in LM : \alpha (0) = x_0\}$. 
By the $C^\infty$ notation we actually mean spaces of \sl piecewise \rm smooth maps. For example a map   $f : [x_0, x_k] \to M$ is \sl piecewise smooth \rm if $f$ is continuous and  if  there exists $x_0 < x_1 < \cdots < x_{k-1} <x_k$ with $f_{|_{(x_i, x_{i+1})}}$ infinitely differentiable for all $i$. 
 These spaces of paths are infinite dimensional smooth manifolds.  See, for example \cite{kling}.
 
The most basic algebraic topological property of closed, oriented manifolds is Poincare duality.  This manifests itself in a homological intersection theory. In their seminal paper, \cite{chas-sullivan}, Chas and Sullivan showed that  certain intersection constructions also exist in the chains and homology of  loop spaces of closed, oriented manifolds.  This endows the homology of the loop space with a rich structure that goes under the heading of ``string topology". 

In this chapter we review Chas and Sullivan's constructions, as well as certain
homotopy theoretic interpretations and generalizations found in \cite{cohen-jones}, \cite{cohen-godin}.  In particular we  recall from   \cite{cohen1}  the ring spectrum structure in the Atiyah dual of a closed manifold, which realizes the intersection pairing in homology, and recall from \cite{cohen-jones} the existence of a related ring spectrum realizing the Chas-Sullivan intersection product (``loop product") in the homology of a loop space.  We also discuss the relationship with Hochschild cohomology
proved in \cite{cohen-jones}, and studied further by \cite{merkulov}, \cite{french}, as well as  the homotopy invariance properties proved in \cite{cohenkleinsull}.
We begin by recalling some basic facts about intersection theory in finite dimensions.

\section{Intersections in compact manifolds}

Let $e : P^p \subset M^d$ be an embedding of closed, oriented manifolds of dimensions
 $p$ and $n$ respectively.  Let  $k$ be  the   codimension,   $k = d-p$.  
 
 Suppose $\theta \in H_q(M^d)$ is represented by an oriented manifold,
 $f : Q^q \to M^d$.  That is, $\theta = f_*([Q])$, where $[Q] \in H_q(Q)$ is the fundamental class.  We may assume that the map $f$ is transverse to the submanifold
 $P \subset M$,  otherwise we perturb $f$ within its homotopy class to achieve transversality.  We then consider the  ``pull-back" manifold 
 $$
 Q\cap P  = \{x \in Q: \, f(x) \in P \subset M \}.
 $$
 This is a dimension $q-k$  manifold, and the map $f$ restricts to give a map
 $f : Q\cap P \to P$.  One therefore has the induced homology class,
 $$
 e_!  (\theta) = f_*([Q\cap P]) \in H_{q-k}(P).
 $$
 
 More generally, on the chain level, the idea is to take 
  a $q$-cycle  in $M$, which is transverse to $P$ in an appropriate sense, and take the intersection to produce a  $q-k$-cycle in $P$.  Homologically, one can make this rigorous by using Poincare duality, to define the intersection   or ``umkehr" map,
$$
e_! : H_q(M) \to H_{q-k}(P)
$$
by the  composition
$$ 
e_! : H_q(M) \cong H^{d-q}(M) \xr{e^*} H^{d-q}(P) \cong H_{q-k}(P)
$$
where the first and last isomorphisms are given by Poincare duality.  

\med
Perhaps the most important example is the diagonal embedding,
$$
\Delta : M \to M \times M.
$$
If we take field coefficients, the induced umkehr map  is the \sl intersection pairing  \rm
$$
\mu = \Delta_! : H_p(M) \otimes H_q(M) \to H_{p+q-d}(M).
$$
Since the diagonal map induces cup product in cohomology, the following diagram commutes:

$$
\begin{CD}
 H_p(M) \otimes H_q(M)   @>\mu >> H_{p+q-d}(M) \\
 @VP.D VV   @VV P.D V \\
 H^{d-p} \times H^{d-q} @>cup >> H^{2d-p-q}
\end{CD}
$$

 In order to deal with the shift in grading, we let $\bh_*(M) = H_{*+d}(M)$.  So $\bh_*(M)$ is nonpositively graded.  
   
\begin{prop} Let $k$ be a field, and $M^d$ an closed, oriented, connected manifold.  Then  $\bh_*(M^d; k)$ is an associative, commutative graded algebra over $k$, together with a    map $\eps: \bh_*(M;k) \to k$ such that the composition
$$
\begin{CD}
\bh_*(M) \times \bh_*(M)  @>\mu >> \bh_*(M) @>\eps >> k
\end{CD}
$$
is a nonsingular bilinear form.  If $k = \bz/2$ the orientation assumption can be dropped.
\end{prop}
In this proposition the map $\eps : \bh_q(M) \to k$ is zero unless $q=-d$, in which case it is the isomorphism
$$
\bh_{-d}(M) = H_0(M) \cong k.
$$

Such an algebraic  structure, namely a commutative algebra $A$ together with a    map $\eps : A \to k$ making the pairing $\langle a, b \rangle = \eps (a\cdot b)$ a nondegenerate bilinear form,  is called a    \sl \bf Frobenius algebra. \rm  

 We leave to the reader the exercise of proving the following, (see \cite{abrams}).
 
  \begin{prop}  A $k$-vector space $A$ is a Frobenious algebra if and only if it is an commutative algebra with unit, it is a co-algebra
$$
\Delta : A \to A \otimes A
$$
with co-unit $\eps : A \to k$, so that $\Delta$ is a map of $A$- bimodules. 
\end{prop}

\med

Intersection theory can also be realized by the ``Thom collapse" map. 
Consider again the embedding of compact manifolds, $e : P \hk M$, and extend $e$ to a 
  tubular neighborhood, $P \subset \eta_e \subset M$.   Consider  the projection map,
\begin{equation}
\tau_e : M \to M/(M-\eta_e).
\end{equation}
Notice that  $M/(M-\eta_e)$ is the one point compactification of the tubular neighborhood,
$M/(M-\eta_e) \cong \eta_e \cup \infty$.  Furthermore,
by the tubular neighborhood theorem, this   space is homeomorphic to the Thom space $P^{\nu_e}$  of the normal bundle, $\nu_e \to P$,
$$M/(M-\eta_e)\cong  \eta_e \cup \infty \cong P^{\nu_e}.$$   So the Thom collapse map can be viewed as a map,
\begin{equation}\label{thomcollapse}
\tau_e : M \to P^{\nu_e}.
\end{equation}

Then the homology intersection map $e_!$ is equal to the composition,
\begin{equation}\label{umkehr}
e_! : H_q(M) \xr{(\tau_e)_*} H_q(P^{\nu_e}) \cong H_{q-k}(P)
\end{equation}
where the last isomorphism is given by the Thom isomorphism theorem.  In fact this description of the  umkehr  map $e_!$ shows that it can be defined in \sl any \rm generalized homology theory, for which there exists a Thom isomorphism for the normal bundle.  This is an orientation condition. 
In these notes we will usually restrict our attention to ordinary homology, but intersection theories in such (co)homology theories as $K$-theory and cobordism theory are very important as well. 

 \section{The Chas-Sullivan loop product}
 
The Chas-Sullivan ``loop product" in the homology of the free loop space of a closed oriented $d$-dimensional manifold,
\begin{equation}\label{loop}
\mu : H_p(LM) \otimes H_q(LM) \to H_{p+q-d}(LM)
\end{equation}
is defined as follows.  

Let $Map (8, M)$ be the mapping space from the figure 8  (i.e the wedge of two circles) to the manifold $M$.  As mentioned above,  the maps are required to be piecewise smooth (see \cite{cohen-jones}).  
Notice that  $Map (8, M)$  can be viewed as  the  subspace of $LM \times LM$ consisting of those pairs of loops that agree at the  basepoint $1 \in S^1$.  In other words, there is a pullback square 
\begin{equation}\label{pullback}
\begin{CD}
 Map (8, M)   @>e >>  LM \times LM \\
 @V ev VV  @VV ev \times ev V \\
 M @>>\Delta > M \times M 
 \end{CD}
\end{equation}
 where $ev : LM \to M$ is the fibration given by evaluating a loop at $1 \in S^1$.  In fact, it can be shown that $ev$ is a locally trivial fiber bundle \cite{kling}.  The map $ev : Map (8, M)\to M$ evaluates the map at the crossing point of the figure 8.   Since $ev \times ev$ is a fibre bundle,
$e : Map(8,M) \hk LM \times LM$ can be viewed as a codimension $d$ embedding, with normal bundle $ev^*(\nu_\Delta) \cong ev^*(TM)$. 

The basic Chas-Sullivan idea, is to take a chain $c \in C_p(LM \times LM)$ that is transverse to the submanifold
$Map (8,M)$ in an appropriate sense, and take the intersection to define a chain $e_!(c)\in C_{p-d}(\met)$.  This will allow
the definition of a map in homology, $e_! : H_*(LM \times LM) \to H_{*-d}(\met)$.  
The striking thing about the Chas-Sullivan construction is that this umkehr map exists in the absence of Poincare duality in this infinite dimensional context.   

\med

   As was done in \cite{cohen-jones},  one can also use the Thom collapse approach to define the umkehr map in this setting.
  They observed that the existence of this pullback diagram of fiber bundles, 
  means that there is a natural tubular neighborhood of the embedding, $e : \met \to LM \times LM$, namely the inverse image of a tubular neighborhood of the diagonal embedding, $\Delta : M \to M \times M$.  That is, 
  $\eta_e = (ev \times ev)^{-1}(\eta_\Delta).$    Because $ev$ is a locally trivial fibration, the tubular neighborhood $\eta_e$  is homeomorphic to the total  space of the normal bundle $ev^*(TM)$.  This induces   a homeomorphism of the  quotient space to the Thom space,
\begin{equation}
  (LM \times LM)/ ((LM \times LM) - \eta_e)   \cong (\met)^{ev^*(TM)}.
  \end{equation}

Combining this homeomorphism with the projection onto this quotient space, defines a Thom-collapse map
\begin{equation}\label{loopthom}
\tau_e : LM \times LM \to Map(8,M)^{ev^*(TM)}.
\end{equation}
For ease of notation, we refer to the Thom space of the pullback bundle, $ev^*(TM) \to \met$ as $(\met)^{TM}.$

Notice that if  if $h_*$ is any generalized homology theory that supports an orientation of $M$ (i.e the tangent bundle $TM$),
then one can define an umkehr map,
\begin{equation}
e_! : h_*(LM \times LM) \xr{\tau_e} h_*((\met)^{TM}) \xr{\cap u} h_{*-d}(\met)
\end{equation} 
where $u \in h^d((\met)^{TM})$ is the Thom class given by the orientation.

\med 
 Chas and Sullivan also observed that given a map from the figure 8 to $M$ then one obtains a loop in $M$ by starting at the intersection point, traversing the top loop of the 8, and then traversing the  bottom loop.  This defines a map
 $$\gamma :  Map(8,M) \to LM.  $$ 
 
 Thought of in a slightly different way, the pullback diagram \ref{pullback} says that we can view $\met$ as the fiber product $\met \cong LM \times_M LM$, as was done in \cite{cohen-jones}. $\gamma$ defines a multiplication map, which by abuse of notation we also call $\gamma : LM \times_M LM \to LM$.  This map   extends the usual multiplication in the based loop space, $\gamma: \Omega M \times \om M \to \om M$. In fact if one may view $ev : LM \to M$ as a fiberwise $H$-space (actually an $H$-group), which is to say  an $H$-group in the category of spaces over $M$.   It is actually an $A_\infty$ space in this category, coming from the $A_\infty$ structure of the multiplication in $\om M$, which is the fiber of $ev : LM \to M$. This aspect of the theory is studied further in \cite{gruher}.  
 
 Chas and Sullivan also observed that the multiplication $\gamma : \met  \to LM$ is homotopy commutative, and indeed there is a canonical,  explicit homotopy.  In the formulas that follow, we identify $S^1 = \br/\bz$.  Now as above, consider   $\met$ as a subspace of $LM \times LM$,
 and suppose $(\alpha, \beta) \in \met$.  We  consider, for each $t \in [0,1]$ a loop
 $\gamma_t (\alpha, \beta)$, which starts at $\beta (-t)$, traverses the arc between
 $\beta (-t)$ and $\beta (0) = \alpha(0)$, then traverses the loop defined by $\alpha$, and then finally traverses the arc between $\beta (0)$ and $\beta (-t)$.  A formula for $\gamma_t ( \alpha, \beta)$ is given by
 
\begin{equation}\label{homotopy}
 \gamma_t(\alpha, \beta)(s) = \begin{cases} 
 \beta (2s-t), \quad \text{for} \quad 0\leq s \leq \frac{t}{2} \\
 \alpha (2s-t ), \quad \text{for} \quad \frac{t}{2} \leq s \leq \frac{t+1}{2} \\
 \beta (2s-t), \quad \text{for} \quad  \frac{t+1}{2} \leq s \leq 1.
 \end{cases}
\end{equation}
 
 One sees that $\gamma_0 (\alpha, \beta) = \gamma (\alpha, \beta)$, and $\gamma_1 (\alpha, \beta) = \gamma (\beta, \alpha)$.  
 
 \med
 The Chas-Sullivan product in homology is defined by composing the umkehr map $e_!$ with the multiplication map $\gamma$.
 
 \begin{definition}\label{loopprod}
Define the loop product in the homology of a loop space to be the composition
$$
\mu : H_*(LM) \otimes H_*(LM) \to H_*(LM \times LM) \xr{e_!} H_{*-d}(\met) \xr{\gamma_*} H_{*-d}(LM).
$$
\end{definition}

Recall that the umkehr map $e_!$ is defined for any generalize homology theory $h_*$ supporting an orientation.   Now suppose that in addition, $h_*$ is   a multiplicative theory.   That is, the corresponding cohomology theory $h^*$ has a cup product, or more precisely, $h_*$ is represented by a \sl ring \rm spectrum.   Then there is a loop product in $h_*(LM)$ as well, 
$$
\mu : h_*(LM) \otimes h_*(LM) \to h_{*-d}(LM).
$$

In order to accomodate the change in grading, one defines
$$
\bh_*(LM) = H_{*+d}(LM).
$$
 Using the naturality of the umkehr map (i.e the naturality of the Thom collapse map) as well as the homotopy commutativity of the multiplication map $\gamma$, the following is proved in \cite{chas-sullivan}.
 
 \begin{thm}\label{loopalgebra} Let $M$ be a compact, closed, oriented manifold.
 Then the loop product defines a map
 $$
 \mu_*: \bh_*(LM) \otimes \bh_*(LM) \to \bh_*(LM)
 $$
 making $\bh_*(LM)$ an associative, commutative algebra.  Furthermore,
 the evaluation map $ev : LM \to M$ defines an algebra homomorphism from the loop algebra to the intersection ring,
 $$
 ev_* : \bh_*(LM) \to \bh_*(M).
 $$
 \end{thm}
 
 As was shown in \cite{cohen-jones}, this structure also applies to $h_*(LM)$, where $h_*$  is any multiplicative generalized homology theory which supports an orientation of $M$. 
 
 \section[The BV structure and the string bracket]{The
 Batalin-Vilkovisky structure and the string bracket}
 
 One aspect of the loop space $LM$ that hasn't yet been exploited is the fact there is an obvious circle action
\begin{equation}\label{action}
 \rho: S^1 \times LM  \la LM
 \end{equation}
 defined by $\rho(t, \alpha) (s) = \alpha (t+s)$. The purpose of this section is to describe those constructions of Chas and Sullivan \cite{chas-sullivan} that exploit this action.
 
 The existence of the $S^1$-action  defines an operator
\begin{align}
 \Delta : H_q(LM) &\to H_{q+1}(LM) \notag \\
 \theta &\to \rho_*(e_1 \otimes \theta) \notag
 \end{align}
 where $e_1 \in H_1(S^1) = \bz$ is the generator.  Notice that if we apply this operator twice, then $\Delta^2 (\theta) = \rho_*(e_1^2 \otimes \theta)$,
 where the product structure in $H_*(S^1)$ is the Pontrjagin algebra structure,  induced by the group structure of $S^1$.  But obviously $e_1^2 =0$ in this algebra,
 and so the operator $\Delta$ has the property that
 $$
 \Delta^2 = 0.
 $$ 
 By regrading, this operator may be viewed as a degree one operator
 on the loop homology algebra,
 $$
 \Delta : \bh_*(LM) \to \bh_{*+1}(LM)
 $$
 and   the following was proved in \cite{chas-sullivan}.
 
 \begin{thm}\label{BV} The pair $(\bh_*(LM), \Delta)$ is a Batalin-Vilkovisky  (BV) algebra.
 That is,
 \begin{enumerate}
 \item $\bh_*(LM)$ is a graded, commutative algebra,
 \item $\Delta \circ \Delta = 0$, \quad \text{and} 
 \item The binary operator defined by the deviation from $\Delta$ being a derivation,
 $$
 \{\phi, \theta \} = (-1)^{|\phi|}\Delta(\phi \cdot \theta) - (-1)^{|\phi|}\Delta (\phi)\cdot \theta - \phi \cdot \Delta (\theta)
 $$
 is a derivation in each variable.
 \end{enumerate}
 
\end{thm}

Now a formal argument given in \cite{chas-sullivan} shows that the degree one binary operator $\{ ,  \}: \bh_*(LM) \times \bh_*(LM) \to \bh_*(LM)$
described above satisfies the (graded) Jacobi identities, and so makes $\bh_*(LM)$ into a graded Lie algebra.  Such a combination
of being a graded, commutative algebra, as well as a Lie algebra which is a derivation in each variable, is called a Gerstenhaber algebra.  So the loop homology algebra has this structure as well. 

Chas and Sullivan also gave another description of   the bracket $\{\phi, \theta \}$.  We give a  variation of their description, which is more homotopy theoretic.  

Let $\cp \subset S^1 \times LM \times LM$
be the space
$$
\cp = \{(t, \alpha, \beta): \alpha (0) = \beta (t) \}.
$$

Notice that there diffeomorphism
$$
h : S^1 \times \met \to \cp
$$
defined by
\begin{equation}\label{diffeo}
h(t, (\alpha, \beta)) = (t, \alpha, \beta_t)
\end{equation}
where $\beta_t(s) = \beta (s-t)$. 

  Notice furthermore there is a pullback square of fibrations,
\begin{equation}\label{cp}
 \begin{CD}
 \cp @>u >\subset>  S^1 \times LM \times LM \\
 @V \eps VV  @VV\eps V \\
 M @>>\Delta >  M \times M
 \end{CD}
 \end{equation}
 where $\eps :   S^1 \times LM \times LM  \to M \times M$ is given by $(t, \alpha, \beta) \to (\alpha (0), \beta (t))$. 
 
 Notice that this pullback diagram has a $\bz/2$-action induced by the diagonal action on $S^1 \times LM \times LM$ given by the antipodal map on $S^1$,  and the permutation action on $LM \times LM$.  The action on $M \times M$ is also the permutation action under which  $\Delta (M)$ is the fixed point set. The diffeomorphism
 $h : S^1 \times \met \to \cp$ is equivariant, where $\bz/2$ acts antipodally on $S^1$ and permutes the two components of the figure 8.

 Like in the previous section, the Thom collapse map may be defined in this situation
 and we get a map 
 $$
 \tau_u : S^1 \times LM \times LM \la \cp^{TM}.
 $$

 Using the Thom isomorphism this defines
 an umkehr map,
 $$
u_!:  H_*(S^1 \times LM \times LM) \to H_{*-d}(\cp)\cong H_{*-d}(S^1 \times \met).
$$
Because of the equivariance, this descends to gives a map on the homology of the orbits,
$$u_! : H_*(S^1 \times_{\bz/2} LM \times LM) \to H_{*-d}(S^1 \times_{\bz/2} \met).$$

This in turn defines a homomorphism
$$
\nu_* :H_*(LM \times LM) \to H_{*+1-d}((S^1 \times_{\bz/2} \met) )
$$
by 
\begin{equation}\label{nu}
\nu_*(\phi \otimes \theta) = u_!(e_1 \otimes (\phi \otimes \theta - (-1)^{(|\phi|+1)(|\theta|+1)}\theta \otimes \phi ).
\end{equation}

We observe that  since the maps $u_!$ and $\nu_*$ are defined in terms of the Thom collapse map $\tau_u$, they can be defined in any generalized homology theory that supports an orientation of $M$. 

\med
Now let $G : [0,1] \times \met \to LM$ be the homotopy given by (\ref{homotopy}).
By identifying $[0,1]$ with the upper semicircle of $S^1$, $G$ defines a map
$$
G : S^1 \times_{\bz/2} \met \to LM.
$$ 
By composing we then have an operation, 
$$
G_* \circ \nu_* : H_*(LM \times LM) \to H_{*+1-d}((S^1 \times_{\bz/2} \met) \to H_{*+1-d}(LM).
$$

This operation is easily seen to be  a global description of the bracket operation defined in definition 4.1 of \cite{chas-sullivan}.
Corollary 5.3  of \cite{chas-sullivan} then gives the following.

\begin{thm}
($G_* \circ \nu_*) (\phi \otimes \theta) = \{\phi, \theta \}.$
\end{thm}

Notice that from this global description, which defines the Batalin-Vilkovisky structure ultimately in terms of the Thom collapse map, we have the following generalization of theorem \ref{BV}, originally  proved in a somewhat different way  by Cohen and Jones in \cite{cohen-jones}.

\med
\begin{thm} Let $h_*$ be any multiplicative generalized homology theory that supports an orientation of $M$. Then $h_*(LM)$ is a Batalin-Vilkovisky algebra.
\end{thm}

\bg
We now turn to the effect of the Batalin-Vilkovisky structure on the equivariant homology,
$H^{S^1}_*(LM)$.     Chas and Sullivan refer to this as the ``string homology" of $M$.
It is the homology of the homotopy orbit space,
$$
\heq = H_*(ES^1 \times_{S^1} LM).
$$
We view this homotopy orbit space as the ``space of closed strings" in $M$ for the following reason.  

Consider the space of embeddings of $S^1$ into an infinite dimensional Euclidean space,
$Emb(S^1, \br^\infty)$.   This is a contractible  space with  an obvious free action of $S^1$, given by reparameterization of the embeddings.   Indeed this action extends to a free action of the homotopy equivalent, but much larger group of orientation preserving diffeomorphisms, $Diff^+(S^1)$.  So the homotopy orbit space $ES^1 \times_{S^1} LM$ is homotopy equivalent to the orbit space,
$Emb(S^1, \br^\infty) \times_{Diff^+(S^1)} LM$ which can be described on the point set level as follows:

$$
Emb(S^1, \br^\infty) \times_{Diff^+(S^1)} LM =  \{(S,f)\}
$$
 where $S\subset \br^\infty$ is a closed,  
 oriented, connected  
  one dimensional submanifold of $\br^\infty$, and $f : S \to M$ is a continuous map.  In other words,
  this orbit space is the space of oriented closed curves (closed strings) in $M$.  Thus the equivariant
  homology, $\heq$ is the homology of the space of closed strings in $M$. 
  
  \med
  Now consider the principal $S^1$ bundle,
  
  \begin{equation}\label{bundle}
  S^1 \to ES^1 \times LM \to ES^1 \times_{S^1} LM.
  \end{equation}
  This gives rise to the Gysin sequence in homology,
  \begin{equation}\label{gysin}
  \to \cdots H_q(LM) \xr{\iota_*} H_q^{S^1}(LM) \xr{j_*} H_{q-2}^{S^1}(LM) \xr{\tau_*} H_{q-1}^{S^1}(LM) \xr{\iota_*} \cdots
  \end{equation}
  We recall that the connecting homomorphism $\tau : H_k^{S^1}(LM) \to H_{k+1}(LM)$ is induced by the $S^1$ - transfer map for this principal bundle , which is defined via a Thom  collapse map.  This is a stable map (map of suspension spectra),
  $$
  \tau : \Sigma^\infty (\Sigma (ES^1 \times_{S^1} X)_+) \to \Sigma^\infty (X_+)
  $$
  which exists for any space $X$ with an $S^1$-action \cite{beckgott}.  Here $\Sigma^\infty$ refers to the suspension spectrum, and $Y_+$ is  $Y$ with a disjoint basepoint.  
  The map $\tau_*$ in the Gysin sequence is homomorphism induced by $\tau$ in homology.  
  
  We note that Chas and Sullivan denoted the homomorphism $\iota_*$ by $E$,  standing for ``erase", and the homomorphism $\tau_*$ by $M$, standing for ``mark".  The motivation for this terminology is that 
  one might consider the loop space $LM$ as the space of ``marked" closed strings, since  the space of ``markings" (i.e choices of marked point) on   a closed, oriented, connected one dimensional manifold $S$, is homeomorphic to $S^1$, which is homotopy equivalent to the space of parameterizations of $S$  by a diffeomorphism from the circle, $S^1 \xr{\cong} S$.   Thus the homomorphism $\iota_* = E$ can be viewed as erasing the marking, and the transfer homomorphism, $\tau_* = M$ can be viewed as taking the map in homology induced by taking all possible markings on $S\subset \br^\infty$.  Since the space of markings is one dimensional, this accounts for the dimension shift in the homomorphism  $\tau_*$.  
  
  Using these maps, Chas and Sullivan define an operator
\begin{align} \label{bracket}
  [\, , \, ] : H_q^{S^1}(LM) \times H_r^{S^1}(LM) &\xr{\tau_* \times \tau_*} H_{q+1}(LM) \times H_{r+1}(LM)  \\  
 & \xr{\circ} H_{q+r+2-d}(LM) \xr{\iota_*} H_{q+r+2-d}^{S^1}(LM). \notag
  \end{align} 
  where  ``$\circ $" is the loop product described above. Using just the formal properties of the Batalin-Vilkovisky structure on $H_*(LM)$, Chas and Sullivan proved the following in \cite{chas-sullivan}.
  
  \begin{thm} The operator $$  [\, , \, ] : H_q^{S^1}(LM) \times H_r^{S^1}(LM) \la   H_{q+r+2-d}^{S^1}(LM)$$ gives the string homology $H_*^{S^1}(LM)$ the structure of a graded Lie algebra of degree $(2-d)$.
  \end{thm}
  
  We note that since the transfer map $\tau$ comes from a map of spectra, there is a corresponding
  string bracket on $h_*^{S^1}(LM) = h_*(ES^1 \times_{S^1}LM)$ for $h_*$ any multiplicative generalized homology theory  that supports an orientation of $M$.

  The string bracket gives the equivariant homology a very rich structure.  For example, if $M$ is a surface, the bracket restricts to define  a Lie algebra structure on the vector space  generated by the path components of closed curves in a surface,
$$   [\, , \, ] : H_0^{S^1}(LM^2) \times H_0^{S^1}(LM^2) \la   H_{0}^{S^1}(LM^2).$$
This bracket operation was originally discovered by Wolpert \cite{wolpert} and Goldman \cite{goldman}
and is highly nontrivial.   It has been studied in more depth by Chas in \cite{chas}.   

\section{A stable homotopy   point  of view}
In this section we describe a homotopy theoretic realization of the loop product due to Cohen and Jones \cite{cohen-jones}.  This takes the form that the Thom spectrum of a certain virtual bundle over the loop space is a ring spectrum, which in homology realizes the loop product.

Recall that the definition of the loop product (\ref{loopprod}), an essential ingredient was the Thom collapse map (\ref{loopthom})
$$\tau_e : LM \times LM \to Map(8,M)^{ev^*(TM)}.$$

The general Thom collapse map (\ref{thomcollapse}) can be modified
by 
  twisting with bundles in the following manner. 

\med
Given $e : N \hk M$, and $\zeta \to M$ vector bundle, consider embedding of total spaces,
$$
\begin{CD}
e^*\zeta @>e_\zeta >\hk > \zeta \\
@VVV  @VVV \\
N @>e>\hk >  M
\end{CD}
$$

The tubular neighborhood of $e_\zeta$ is homeomorphic to total space of $e^*\zeta \oplus \eta_e$.  So we get Thom collapse map
$$
\tau : M^\zeta \to N^{e^*\zeta \oplus \eta_e}
$$

\med
This construction can be carried out even if 
  $\zeta$ is a \sl virtual \rm bundle,  
$\zeta = \gamma_1 - \gamma_2$, where $\gamma_i$'s are vector bundles. We view a virtual bundle as an element of $K$-theory, and there is a dimension homomorphism,
  $\dim : K(M) \to \bz$, where $\dim (\zeta) = \dim (\gamma_1) -  \dim (\gamma_2).$  So the dimension is possibly negative.   
  
 In this setting, the
Thom isomorphism still holds,  $H^q(M) \cong H^{q+\dim \zeta}M^\zeta,$    where again, $\dim \, \zeta$ might be negative.    

\med
\bfl \bf Example.
\rm
\efl
Let $\zeta \to M$ be a vector bundle.  Consider the product bundle $\zeta \times \zeta \to M \times M$, and embedding $\Delta : M \hk M\times M$.  The Thom collapse is a map of spectra,
$$
\tau : M^\zeta \wedge M^\zeta \la  M^{2\zeta \oplus TM}
$$

Now
take $\zeta$ to be the virtual bundle,  $\zeta = -TM$.  The Thom collapse is then  a map  

\begin{equation}\label{mult}
\tau : \mtm \wedge \mtm \la \mtm
\end{equation}

This defines a ring structure on $\mtm$ which was studied in detail in \cite{cohen1}.
We describe this structure more explicitly as follows.
 
\ 
\rm Let $e : M \to \br^L$ be an embedding, with  $\eta_e$ its normal bundle. Then there is an equivalence between the $L$-fold suspension of the Thom spectrum $\mtm$ and the Thom space of the normal bundle, $M^{\eta_e}$,
$$
\Sigma^L\mtm \simeq  M^{\eta_e}.
$$

For $\eps > 0$, 
let $\nu_e^\eps $ be the $\eps$-tubular neighborhood:
$=\{y \in \br^L \, : \, d(y, e(M)) < \eps \}$.  Then for $\eps$ sufficiently small,
$$
M^{\eta_e} \cong \br^L/(\br^L - \nu_e^\eps).
$$

 Let $B_\eps (0)$ be the ball of radius $\eps$ around the origin in $\br^L$.  
In the 1930's, Alexander considered the map
\begin{align}
(\br^L - \nu_e^\eps) \times M &\la \br^L - B_\eps (0) \simeq S^{L-1} \notag \\
v, x &\la v - e(x)  \notag
\end{align}

This defines a map in homology,
$$H_{L-q-1}(\br^L - \nu_e^\eps) \otimes H_q(M) \la H_{L-1}S^{L-1} \cong \bz$$ or by taking the adjoint,
$$
H_q(M) \la H^{L-q-1}(\br^L -\nu_e^\eps)
$$
This is the famous
``Alexander duality"   isomorphism. 
\med

In early 1960's, Atiyah \cite{atiyah} considered the same map defined on the quotient spaces,
\begin{align}
\br^L \times M / ((\br^L - \nu_e^\eps) \times M) = M^{\eta_e}\wedge M_+  &\to \br^L /(\br^L - B_\eps (0)) \notag \\
& \cong S^L \notag \\
v, x &\la v - e(x)  \notag
\end{align}

This yields a map
$$
\alpha : M^{\eta_e} \to Map (M_+, S^L).
$$

On the spectrum level, $\alpha$ induces a map
$$
\alpha :  \mtm \to F(M_+, S),
$$
   where  $F(M_+, S)$ is the space of  stable maps from $M_+$ to the sphere spectrum $S$.  As a spectrum, its $k^{th}$ space is given by $Map (M_+, S^k)$. 
   
   The following theorem, proved by Atiyah in \cite{atiyah}   says that this map is a homotopy equivalence.

\med

\begin{thm} For a closed manifold $M$, the map $\alpha : \mtm \to  F(M_+, S), $ is homotopy equivalence from the Thom spectrum of $-TM$  to the Spanier-Whitehead dual,
of $M_+$.  
\end{thm}

\med
The Spanier-Whitehead dual, $ F(M_+, S)$ is clearly a ring spectrum.  Its multiplication is dual to the diagonal
map $\Delta : M \to M \times M$.  Since this diagonal map is cocommutative, this says that the ring structure
on $ F(M_+, S)$ also has commutativity properties.  In recent years, symmetric monoidal categories of spectra have been developed (\cite{maygroup}, \cite{hss}) in which the appropriate coherence issues regarding homotopy commutativity can be addressed.  In \cite{cohen1}, Cohen used the notion of symmetric spectrum developed in \cite{hss} and proved the following. 

\med

\begin{thm}\label{atring}   $\alpha$ is an equivalence   of commutative, symmetric ring spectra
and of bimodules over $F(M_+, S)$.  
\end{thm}

\med

That is, the ring structure on $\mtm$ induced by the Thom collapse map (\ref{mult}), can be rigidified
to give a commutative ring structure which coincides, via the classical map of Alexander, to the commutative ring structure on $F(M_+,S)$ given by the dual of the diagonal map.

\med
Now as seen above, the  ring map (\ref{mult}) is determined by the Thom collapse map
for the pullback diagram of virtual bundles,
$$
\begin{CD}
-2 TM  @>\Delta >> -TM \times -TM \\
@VVV   @VVV \\
M @>\hk>\Delta > M \times M.
\end{CD}
$$

We now pull this diagram of virtual bundles back, using the evaluation map, $ev : LM \to M$,
to the total spaces of pullback square 
 (\ref{pullback}).  
 
  $$
 \begin{CD}
 ev^*(-2TM) @>e>> ev^*(-TM )\times ev^*(-TM) \\
 @VVV  @VVV \\
 \met @>\hk >e >  LM \times LM
 \end{CD}
 $$

 This defines a Thom collapse map,
 
$$
\tau_e : \ltm \wedge \ltm \to \met^{-TM}
$$
which yields a product
$$
\begin{CD}
\mu : \ltm \wedge \ltm @>\tau_e >>  \met^{-TM} @>\gamma>> \ltm.
\end{CD}
$$

In \cite{cohen-jones} Cohen and Jones proved the following.

\begin{thm} For any closed manifold $M$,  $\ltm$ is a ring spectrum.  If $M$ is oriented, then the Thom isomorphism in homology,    
$$
H_*(\ltm) \cong H_{*+n}LM \cong \bh_*(LM)
$$
is an isomorphism of commutative graded algebras. The analogous statement holds for any multiplicative generalized homology theory $h_*$ that supports an orientation of $M$.   Furthermore the evaluation map 
$$
ev  : \ltm \to \mtm
$$
is a map of ring spectra.  Also, 
  if one considers pullback diagram
$$
\begin{CD}
\Omega M  @>\hk >> LM \\
@VVV  @VVev_0 V \\
point  @>>\hk > M
\end{CD}
$$
then  the induced Thom collapse map
$$
\ltm \to \Sigma^\infty (\Omega M_+)
$$
is a map of ring spectra.
\end{thm}

\med
The existence of ring structures on Thom spectra of more general fiber products was examined by Klein in \cite{klein}.

\med
This theorem says that even on the stable homotopy level, the  string topology structure on the loop space is compatible with both the intersection product on $M$,  and Pontrjagin product on $\Omega M$.  

 This compatibility was strengthened by Cohen, Jones, and Yan in \cite{cohenjonesyan} for calculational purposes.  Namely, 
 consider Serre  spectral sequence for the  fibration
$$\begin{CD}
\om M \to LM @>ev >> M
\end{CD}
$$   Assume $M$ is simply connected. Then the $E_2$ term is given by  $E_2^{*,*} = H_*(M, H_*(\om M)) $, and the spectral sequence converges to  $H_*(LM)$.  By applying  Poincare duality, one gets  a second quadrant spectral sequence
$$
E^{-s,t}_2 = H^s(M; H_t(\om M)) \rightrightarrows H_{n-s+t}(LM) = \bh_{t-s}(LM).
$$
Notice that this $E_2$-term is an algebra, using cup product in cohomology and the  Pontrjagin product on the coefficients.  It was shown in \cite{cohenjonesyan}  that this spectral sequence is multiplicative.   That is, each $E_r^{*,*}$ is a graded ring, and the differentials are derivations.  It converges multiplicatively to the loop product on $\bh_*(LM)$.  Cohen, Jones, and Yan then used this spectral sequence to calculate the ring structure in $\bh_*(L\bc \bp^n)$ and $\bh_*(L S^m)$ for all $n$ and $m$, and demonstrated what a rich structure this algebra detects.  

\med
 
\section{ Relation to Hochschild cohomology}

In this section we describe an algebraic point of view of the string product.   Namely we recall how Hochschild cohomology of the cochains of a manifold, $H^*(C^*(M), C^*(M))$ was shown in \cite{cohen-jones} to be isomorphic as rings to $\bh_*(LM)$.  
This relationship is an outgrowth of the cosimplicial model of the loop space and its relationship to Hochschild homology established by Jones in \cite{jones}. 

We begin by recalling the basic construction.  Consider the standard simplicial decomposition of the circle,
$S^1$  which has one   $0$-simplex and one nondegenerate
$1$-simplex.

In a simplicial set $\cs_*$, there are face maps $\p_i : \cs_k \to \cs_{k-1},$ $i = 0, \cdots k$, and degeneracy maps $\sigma_j : \cs_k \to \cs_{k+1}$, $j=0, \cdots , k$.
So $S^1_*$ is generated  with respect to the degeneracies  by the one zero simplex and one nondegenerate one-simplex.  Using the relations satisfied by the face and degeneracy operators, it turns out that the set of $k$- simplices, $S^1_k$ is a set with $k+1$ elements.  
$$
S^1_k =  \{{\bf k+1}\}
$$
All the simplices are degenerate if $k>1$.  

Now
recall the coface and codegeneracy maps    on the standard simplices:
\begin{align}
d^i : \Delta^{k-1} &\to \Delta^k  \quad  i=0, \cdots , k \notag \\
s^j : \Delta^{k+1} &\to \Delta^k  \quad  j=0, \cdots , k. \notag
\end{align}

We then have the resulting homeomorphism of the geometric realization
$$
S^1 \cong \bigcup_k \Delta^k \times S^1_k / \sim \quad  = \bigcup_k \Delta^k \times \{{\bf k+1}\} / \sim
$$
where $(s^j(t), x) \sim (t, \sigma_j(x))$, $(d^i(t), y) \sim (t, d_i(y))$. 

\med
This gives an embedding of the loop space,

\begin{align}\label{fk}
f = \prod_k f_k : LX = Map (S^1, X) &= Map( \bigcup_k \Delta^k \times \{\bf k+1\}\rm / \sim, \, X)   \\
&\subset Map (\coprod_k \Delta^k \times  \{\bf k+1\}\rm ; \, X) \notag \\
&= \prod_k Map (\Delta^k; \, X^{k+1}) \notag
\end{align}

The image of this embedding are   those sequences of maps that commute with the  coface and codegeneracy operators.

 The component maps $f_k : LX \to Map(\Delta^k, X)$ can be described explicitly as follows. 
\begin{align}
f_k : \Delta^k \times LX &\la X^{k+1} \notag \\
(0 \leq t_1 \leq \cdots \leq t_k \leq 1 ; \, \gamma)  &\to ( \gamma (t_1), \cdots \gamma (t_k), \gamma (1) ). \notag
\end{align}

If we consider the homomorphism on singular cochains induced by $f_k$, and then take the cap product
with the canonical $k$-simplex in the chains of $\Delta^k$, we get a diagram, which in \cite{jones} was shown to commute:   $$
\begin{CD}
C^{*-k}(LX)   @<f_k^*<< C^*(X)^{\otimes k+1} \\
@V\delta VV    @VV total \, diff V \\
C^{*-k+1}(LX)  @<<f_{k+1}^* <  C^*(X)^{\otimes k}
\end{CD}
$$
Here  the right hand vertical map is the total differential in the Hochschild chain complex.  Recall that the Hochschild chain complex    of a differential graded algebra, with coefficients in  a bimodule $C$ is the complex $CH_*(A, C)$:
$$
\la \cdots \xr{\p} C \otimes A^{\otimes n}   \xr{\p} C \otimes A^{\otimes n-1}  \xr{\p} \cdots
$$
where $\p$ is the total differential, given by the sum of the internal differential on $A^{\otimes n} \otimes C$
and the Hochschild boundary operator
\begin{align}
b(c \otimes a_1 \otimes \cdots \otimes a_n )  =  c\cdot a_1 \otimes a_2 \otimes \cdots \otimes a_n \quad + &\sum_{i=1}^{n-1}(-1)^{i} c \otimes a_1 \otimes \cdots a_i\cdot a_{i+1}\otimes \cdots \otimes a_n  \notag \\
  + &(-1)^n a_n \cdot c \otimes  a_1 \otimes \cdots \otimes a_{n-1}. \notag
\end{align}

In the case above, we are considering the Hochschild complex $CH_*(C^*(X), C^*(X))$.  The following was proved by Jones in \cite{jones}.

\begin{thm}\label{hochjones}  For simply connected $X$,
$$
f^* : CH_*(C^*(X), C^*(X)) \la C^*(LX)
$$
is a chain homotopy equivalence.   It therefore induces an isomorphism
$$
f^*:  H_*(C^*(X), C^*(X))    \xr{\cong}  H^*(LX).
$$
\end{thm}

Now the loop product lives in the homology $H_*(LM)$, so to get a model for this, we dualize the
Hochschild complex, and we get a complex

$$
\la \cdots \la Hom (C^*(X)^{\otimes q }, k)  \xr{\delta} Hom (C^*(X)^{\otimes k+1}, k)\xr{\delta} \cdots
$$
which computes $H_*(LX; k)$.  But with respect to the obvious identification,
$Hom (C^*(X)^{\otimes q+1}, k) \cong Hom (C^*(X)^{\otimes q}; C_*(X))$, this complex is the Hochschild cochain complex of $C^*(X)$ with coefficients in the bimodule $C_*(X)$.  We therefore have the following corollary.

\begin{crl}  For simply connected $X$, there is an isomorphism
$$f_* : H_*(LX) \xr{\cong} H^*(C^*(X); C_*(X)).$$
\end{crl}

\med
Now let $X = M^n$ be a simply connected,  oriented, closed manifold.  Notice that the following diagram commutes:
$$
\begin{CD}
\Delta^q \times LM @>f_q >>  M^{q+1} \\
@Vev VV  @VV p_{q+1} V \\
M  @>> = > M
\end{CD}
$$
where $p_{q+1}$ is the projection onto the last factor.  This implies we have a map of Thom spectra,
$$
f_q : \Delta^q_+ \wedge LM^{-TM} \to M^q_+ \wedge \mtm 
$$
for each $q$.  It was shown in \cite{cohen-jones} that as a consequence of the embedding (\ref{fk}) and theorem \ref{hochjones}, we have the following:

 \begin{thm} For $M$ simply connected, the map $f$ induces an isomorphism of rings,
$$
\begin{CD} 
f_* : \bh_*(LM) \cong H_*(\ltm) @>\cong >> H^*(C^*(M); C_*(\mtm))
\end{CD}
$$    The ring structure of the Hochschild cohomology is given by cup product, where one 
is using the ring spectrum structure of the Atiyah dual, $\mtm$ to give a ring structure to the coefficients,
$C_*(\mtm)$.
\end{thm}

 As a consequence of the fact that the  Atiyah duality map $\alpha$ is an equivalence of ring spectra, (theorem \ref{atring}), it was shown in \cite{cohen1} that $\alpha$ induces an isomorphism,
 $$
 \alpha_* :  H^*(C^*(M); C_*(\mtm)) \xr{\cong} H^*(C^*(M), C^*(M)).
 $$
 This then implies the following theorem \cite{cohen-jones}: 

\begin{thm}  The composition
 $$ \psi:   \bh_*(LM) \xr{f_*} H^*(C^*(M); C_*(\mtm)) \xr{\alpha_*} H^*(C^*(M), C^*(M))$$ is an isomorphism of graded algebras.
\end{thm}

\med

We end this section with a few comments and observations about recent work in this direction.

\bfl
\bf Comments. \rm

\efl
\begin{enumerate}
\item Clearly the Hochschild cohomology, $H^*(C^*(M), C*(M))$ only depends on the homotopy type of $M$.  However the definition of the loop product and Batalin-Vilkovisky
structure on $\bh_*(LM)$ involves intersection theory, and as we've described it above,
the Thom collapse maps.  These constructions involve the smooth structure of $M$.
Indeed, even the definition of the isomorphism, $\psi : \bh_*(LM) \to H^*(C^*(M), C*(M))$
given in \cite{cohen-jones} involves the Thom collapse map, and therefore the smooth
structure of $M$.  Nonetheless, using the Poincare embedding theory of Klein \cite{klein2},  Cohen, Klein, and Sullivan have recently proved that if $h_*$ is a multiplicative generalized homology theory supporting an orientation of $M$, and
$f : M_1 \to M_2$ is an $h_*$-orientation preserving  homotopy equivalence of simply connected, closed manifolds, then the induced homotopy equivalence of loop spaces,
$Lf : LM_1 \to LM_2$ induces an isomorphism of BV-algebras, $(Lf)_* :  h_*(LM_1) \to h_*(LM_2) $ \cite{cohenkleinsull}.

\item Felix, Menichi, and Thomas \cite{french} proved that the Hochschild cohomology
$H^*(C^*(M), C^*(M))$ is a BV algebra.  It is expected that the ring isomorphism
$$
\psi : \bh_*(LM) \to H^*(C^*(M), C^*(M))
$$
preserves the BV structure.  However this has not yet been proved.  They also show
that a kind of Koszul duality implies that
there is an isomorphism of Hochschild cohomologies,
\begin{equation}\label{koszul}
H^*(C^*(M), C^*(M)) \cong H^*(C_*(\Omega M), C_*(\om M))
\end{equation}
where the chains on the based loop space $C_*(\om M)$ has the Pontrjagin
algebra structure.   This is significant because
of the alternative description of the homology  $H_*(LX)$ in terms of Hochschild
homology \cite{goodwillie}
$$
H_*(LX) \cong  H_*(C_*(\Omega X), C_*(\om X))
$$
for any $X$.  ($X$ need not be simply connected in this case.)

\item The topological Hochschild cohomology of a ring spectrum $R$,   $THH^*(R)$ (see \cite{bokhsiangmad} ) is a cosimplicial spectrum which has the the structure of an algebra over the little disk operad by work of McClure and Smith.  The cosimplicial model for the suspension spectrum $\ltm$ constructed by Cohen and Jones \cite{cohen-jones}
is equivalent to $THH^*(\mtm)$, which by Atiyah duality (theorem (\ref{atring})) is equivalent to $THH^*(F(M_+, S))$.  It was observed originally by Dwyer and Miller, as well as Klein \cite{klein} that the topological Hochschild cohomology $THH^*(\Sigma^\infty (\om (M_+))$ is also homotopy equivalent to $\ltm$.  The spectrum homology of these spectra are given by
\begin{align}
H_*(THH^*(F(M_+, S))) &\cong H^*(C^*(M), C^*(M)) \notag \\
H_*(THH^*(\Sigma^\infty (\om (M_+)))) &\cong  H^*(C_*(\Omega X), C_*(\om X)) \notag
\end{align}
and so the equivalence,
$$
 H_*(THH^*(F(M_+, S)))  \simeq  \Sigma^\infty (LM_+) \simeq H_*(THH^*(\Sigma^\infty (\om M_+ ))
 $$
 realizes, on the spectrum level, the Koszul duality isomorphism (\ref{koszul}) above. 

\end{enumerate}

% Version of 03/17/05

\chapter{The cacti operad}
\label{cactioperad}

\section{PROPs and operads}
\label{operads}

Operads in general are spaces of operations with certain rules on how
to compose the operations. In this sense operads are directly related
to Lawvere's algebraic theories and represent true objects of
universal algebra. However, operads as such appeared in topology in
the works of J.~P. May \cite{may}, J.~M. Boardman and R.~M. Vogt
\cite{boardman-vogt} as a recognition tool for based multiple loop
spaces. Stasheff \cite{stasheff:HI} earlier described the first example
of an operad, the associahedra, which recognized based loop
spaces. About the same time, Gerstenhaber
\cite{gerstenhaber:composition}, studying the algebra of the
Hochschild complex, introduced the notion of a composition algebra,
which is equivalent to the notion of an operad of graded vector
spaces.

\subsection{PROP's}

We will start with defining the notion of a PROP (=PROducts and
Permutations) and think of an operad as certain part of a PROP.
However, later we will give an independent definition of an operad.

\begin{definition}
  A \emph{PROP} is a symmetric monoidal (sometimes called tensor)
  category whose set of objects is identified with the set $\nz_+$ of
  nonnegative integers. The monoidal law on $\nz_+$ is given by
  addition and the associativity transformation $\alpha$ is equal to
  identity. See the founding fathers' sources, such as, J.~F. Adams'
  book \cite{Ad} or S.~Mac Lane's paper \cite{maclane} for more
  detail.
\end{definition}

Usually, PROP's are enriched over another symmetric monoidal category,
that is, the morphisms in the PROP are taken as objects of the other
symmetric monoidal category. This gives the notions of a PROP of sets,
vector spaces, complexes, topological spaces, manifolds, etc.
Examples of PROP's include the following. We will only specify the
morphisms, because the objects are already given by the definition.

\begin{ex}
  The \emph{endomorphism PROP} of a vector space $V$ has the space of
  morphisms $\Mor (m,n) = \Hom(V^{\otimes m}, V^{\otimes n})$. This is
  a PROP of vector spaces. The composition and tensor product of
  morphisms are defined as the corresponding operations on linear
  maps.
\end{ex}

\begin{ex}[Segal \cite{segal:elliptic}]
\label{Segal:PROP}
  The \emph{Segal PROP} is a PROP of infinite dimensional complex
  manifolds. A morphism is defined as a point in the moduli space
  $\PP_{m,n}$ of isomorphism classes of complex Riemann surfaces
  bounding $m+n$ labeled nonoverlapping holomorphic holes. The
  surfaces should be understood as compact smooth complex curves, not
  necessarily connected, along with $m+n$ biholomorphic maps of the
  closed unit disk to the surface, thought of as holes. The
  biholomorphic maps are part of the data, which in particular means
  that choosing a different biholomorphic map for the same hole is
  likely to change the point in the moduli space. The more precise
  nonoverlapping condition is that the closed disks in the inputs do
  not intersect pairwise and the closed disks in the outputs do not
  intersect pairwise, however, an input and an output disk may have
  common boundary, but are still not allowed to intersect at an
  interior point. This technicality brings in the identity morphisms
  to the PROP, but does not create singular Riemann surfaces by
  composition. The composition of morphisms in this PROP is given by
  sewing the Riemann surfaces along the boundaries, using the equation
  $zw=1$ in the holomorphic parameters coming from the standard one on
  the unit disk. The tensor product of morphisms is the disjoint
  union. This PROP plays a crucial role in Conformal Field Theory, as
  we will see now.
\end{ex}

\subsection{Algebras over a PROP}

We need to define another important notion before we proceed.
\begin{definition}
  We say that a vector space $V$ is an \emph{algebra over a PROP} $P$,
  if a morphism of PROP's from $P$ to the endomorphism PROP of $V$ is
  given. A morphism of PROP's is a functor respecting the symmetric
  monoidal structures and also equal to the identity map on the
  objects.
\end{definition}

An algebra over a PROP could have been called a \emph{representation},
but since algebras over operads, which are similar objects, are
nothing but familiar types of algebras, it is more common to use the
term ``algebra.''

\begin{ex}
  An example of an algebra over a PROP is a \emph{Conformal Field
    Theory $($CFT$)$}, which may be defined (in the case of a
  vanishing central charge) as an algebra over the Segal PROP. The
  fact that the functor respects compositions of morphisms translates
  into the sewing axiom of CFT in the sense of G.~Segal. Usually, one
  also asks for the functor to depend smoothly on the point in the
  moduli space $\PP_{m,n}$. One needs to extend the Segal PROP by a
  line bundle to cover the case of an arbitrary charge, see
  \cite{huang:book}.
\end{ex}

\begin{ex}[Sullivan]
  Another example of an algebra over a PROP is a Lie bialgebra.  There
  is a nice graph description of the corresponding PROP, about which
  we learned from Sullivan, see \cite{markl-voronov}.
\end{ex}

\subsection{Operads}

Now we are ready to deal with operads, which formalize the notion of a
space of operations, as we mentioned in the introduction to
Section~\ref{operads}. Informally, an operad is the part $\Mor(n,1)$,
$n \ge 0$, of a PROP. Of course, given only the collection of
morphisms $\Mor(n,1)$, it is not clear how to compose them. The idea
is to take the union of $m$ elements from $\Mor(n,1)$ and compose them
with an element of $\Mor(m,1)$. This leads to cumbersome notation and
ugly axioms, compared to those of a PROP. However operads are in a
sense more basic than the corresponding PROP's; the difference is
similar to the difference between Lie algebras and the universal
enveloping algebras.

\begin{definition}[May \cite{may}]
An \emph{operad} $\OO$ is a collection of sets (vector spaces,
complexes, topological spaces, manifolds, \dots, objects of a
symmetric monoidal category) $\OO (n)$, $n \ge 0$, with
\begin{enumerate}
\item A composition law:
\end{enumerate}
\noindent
\begin{equation*}
\gamma: \OO(m) \otimes \OO(n_1) \otimes \dotsb \otimes \OO(n_m) \to
\OO(n_1 + \dotsb + n_m).
\end{equation*}
\begin{enumerate}
\setcounter{enumi}{1}
\item A right action of the symmetric group $\Si_n$ on $\OO (n)$.
\item A unit $e \in \OO (1)$.
\end{enumerate}
such that the following properties are satisfied:
\begin{enumerate}
\item The composition is associative, \emph{i.e}., the following
diagram is commutative:
\end{enumerate}
\noindent
\[
\begin{CD}
\left\{
\begin{aligned}
\OO(l) & \otimes \OO(m_1) \otimes \dotsb \otimes \OO(m_l) \\
& \otimes \OO(n_{11}) \otimes \dotsb \otimes \OO(n_{l,n_l})
\end{aligned}
\right\}
@>{\id \otimes \gamma^l}>>
\OO(l) \otimes \OO(n_1) \otimes \dotsb \otimes \OO(n_l)\\
@V{\gamma \otimes \id}VV        @VV{\gamma}V , \\
\OO(m) \otimes \OO(n_{11}) \otimes \dotsb \otimes \OO(n_{m,n_m})
@>\gamma>> \OO(n)
\end{CD}
\]
\noindent
\begin{quote}
where $m = \sum_i m_i$, $n_i = \sum_j n_{ij}$, and $n = \sum_i n_i$.
\end{quote}
\noindent
\begin{enumerate}
\setcounter{enumi}{1}
\item The composition is equivariant with respect to the symmetric
  group actions: the groups $\Si_m$, $\Si_{n_1}$, \dots, $\Si_{n_m}$ act on
  the left-hand side and map naturally to $\Si_{n_1 + \dotsb + n_m}$,
  acting on the right-hand side.
\item The unit $e$ satisfies natural properties with respect to the
composition: $\gamma(e; f) \linebreak[0] = f$ and $\gamma(f; e, \dots,
e) = f$ for each $f \in \OO (k)$.
\end{enumerate}

The notion of a \emph{morphism of operads} is introduced naturally.
\end{definition}

\begin{rem}
One can consider \emph{non-$\Sigma$ operads}, not assuming the action
of the symmetric groups. Not requiring the existence of a unit $e$, we
arrive at \emph{nonunital operads}. Do not mix this up with operads
with no $\OO(0)$, algebras over which (see next section) have no
unit. There are also good examples of operads having only $n \ge 2$
components $\OO(n)$.

An equivalent definition of an operad may be given in terms of
operations $f \circ_i g \linebreak[1] = \linebreak[0] \gamma (f; \id,
\dots, \linebreak[1] \id, g, \id , \dots, \id)$, $i = 1, \dots, m$,
for $f \in \OO(m), g \in \OO(n)$. Then the associativity condition
translates as $f \circ_i (g \circ_j h) = (f \circ_i g) \circ_{i+j-1}
h$ plus a natural symmetry condition for $(f \circ_i g) \circ_j h$,
when $g$ and $h$ ``fall into separate slots'' in $f$, see \emph{e.g}.,
\cite{ksv2}.
\end{rem}

%\begin{sloppypar}
\begin{ex}[The Riemann surface and the endomorphism operads]
$\PP(n)$ is the space of Riemann spheres with $n+1$ boundary
components, \emph{i.e}., $n$ inputs and 1 output.  Another example is
the \emph{endomorphism operad of a vector space $V$}: $\End{V}(n) =
\Hom(V^{\otimes n}, V)$, the space of $n$-linear mappings from $V$ to
$V$.
\end{ex}
%\end{sloppypar}

\subsection{Algebras over an operad}

\begin{definition}
An \emph{algebra over an operad $\OO$} (in other terminology, a
\emph{representation of an operad}) is a morphism of operads $\OO \to
\End{V}$, that is, a collection of maps
\[
\OO(n) \to \End{V}(n) \qquad \text{for $ n \ge 0$}
\]
compatible with the symmetric group action, the unit elements, and the
compositions. If the operad $\OO$ is an operad of vector spaces, then
we would usually require the morphism $\OO \to \End{V}$ to be a
morphism of operads of vector spaces. Otherwise, we would think of
this morphism as a morphism of operads of sets. Sometimes, we may also
need a morphism to be continuous or respect differentials, or have
other compatibility conditions. We will also consider a nonlinear
version of the notion of an operad algebra, which may be defined in
any symmetric monoidal category. For example, an $\OO$-algebra $X$ in
the category of topological spaces would be an operad morphism $\OO(n)
\to \Map (X^n, X)$ for $n \ge 0$, where $\Map$ is the space of
continuous maps.
\end{definition}

\subsubsection{The commutative operad}
The \emph{commutative operad} is the operad of $k$-vector spaces with
the $n$th component $\comm (n) = k$ for all $n \ge 0$. We assume that
the symmetric group acts trivially on $k$ and the compositions are
just the multiplication of elements in the ground field $k$.  The term
``commutative operad'' may seem confusing to some people, but it has
been in use for a while. An algebra over the commutative operad is
nothing but a commutative associative algebra with a unit, as we see
from Exercise 2 below.

Another version of the commutative operad is $\comm(n) =
\{\text{point}\}$ for all $n \ge 0$. This is an operad of sets. It is
equivalent to the previous version in the sense that an algebra over
it is the same as a commutative associative unital algebra.

\begin{xca}
Show that the operad $\mathcal{T}op (n) \linebreak[1] = \linebreak[0]$
\{the set of diffeomorphism classes of Riemann spheres with $n$ input
holes and 1 output hole\} is isomorphic to the commutative operad of
sets.
\end{xca}

\begin{xca}
Prove that the structure of an algebra over the commutative operad
$\comm$ on a vector space is equivalent to the structure of a
commutative associative algebra with a unit.
\end{xca}

\subsubsection{The associative operad}
\label{ass:operad}

\begin{figure}
\centerline{\includegraphics[height=2.5cm]{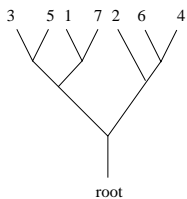}}
\caption{A planar binary tree}
\label{tree}
\end{figure}

The \emph{associative operad} $\ass$ can be considered as a planar
one-dimensional analogue of the commutative operad $\mathcal{T}op$.
$\ass(n)$ is the set of equivalence classes of connected planar binary
(each vertex being of valence 3) trees that have a root edge and $n$
leaves labeled by integers 1 through $n$, see Figure~\ref{tree}.
\noindent
If $n=1$, there is only one tree --- it has no vertices and only one
edge connecting a leaf and a root. If $n=0$, the only tree is the one
with no vertices and no leaves --- it only has a root.

\begin{figure}
\centerline{\includegraphics[height=2cm]{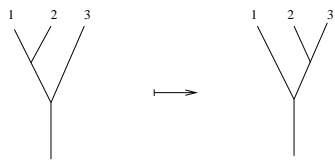}}
\caption{An equivalence move}
\label{move}
\end{figure}

Two trees are equivalent if they are related by a sequence of moves of
the kind pictured on Figure~\ref{move},
performed over pairs of two adjacent vertices of a tree. The symmetric
group acts by relabeling the leaves, as usual. The composition is
obtained by grafting the roots of $m$ trees to the leaves of an
$m$-tree, no new vertices being created at the grafting points. Note
that this is similar to sewing Riemann surfaces and erasing the seam,
just as we did to define operad composition in that case. By
definition, grafting a 0-tree to a leaf just removes the leaf and, if
this operation creates a vertex of valence 2, we should erase the
vertex.

\begin{xca}
Prove that the structure of an algebra over the associative operad
$\ass$ on a vector space is equivalent to the structure of an
associative algebra with a unit.
\end{xca}

\subsubsection{The Lie operad}
\label{lie:operad}

\begin{figure}
\centerline{\includegraphics[height=2.5cm]{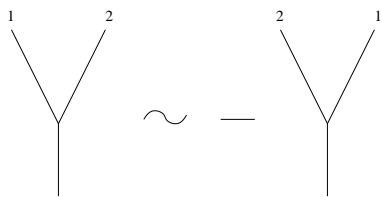}}
\caption{Skew symmetry}
\label{skew}
\end{figure}

\begin{figure}
\centerline{\includegraphics[height=2.5cm]{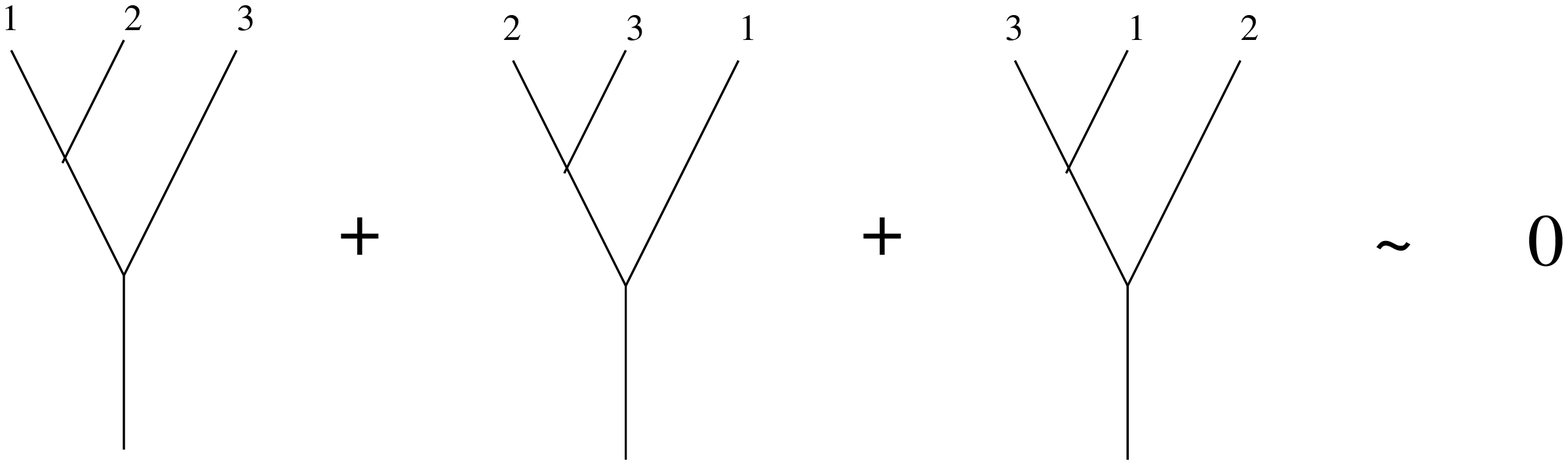}}
\caption{The Jacobi identity}
\label{jacobi}
\end{figure}

The \emph{Lie operad} $\lie$ is another variation on the theme of a
tree operad. Consider the vector space spanned by the same planar
binary trees as for the associative operad, except that we do not
include a 0-tree, \emph{i.e}., the operad has only positive components
$\lie(n)$, $n \ge 1$, and there are now two kinds of equivalence
relations, see Figures \ref{skew} and \ref{jacobi}.
Now that we have arithmetic operations in the equivalence relations,
we consider the Lie operad as an operad of vector spaces. We also
assume that the ground field is of a characteristic other than 2,
because otherwise we will arrive at the wrong definition of a Lie
algebra.

\begin{xca}
  Prove that the structure of an algebra over the Lie operad $\lie$ on
  a vector space over a field of a characteristic other than 2 is
  equivalent to the structure of a Lie algebra.
\end{xca}

\begin{xca}
Describe algebraically an algebra over the operad $\lie$, if we modify
it by including a 0-tree, whose composition with any other tree is
defined as (a) zero, (b) the one for the associative operad.
\end{xca}

\subsubsection{The Poisson operad}
\label{poisson:operad}

Recall that a \emph{Poisson algebra} is a vector space $V$ (over a
field of $\Char \ne 2$) with a unit element $e$, a dot product $ab$,
and a bracket $[a,b]$ defined, so that the dot product defines the
structure of a commutative associative unital algebra, the bracket
defines the structure of a Lie algebra, and the bracket is a
derivation of the dot product:
\[
[a,bc] = [a,b] c + b [a,c] \qquad \text{for all $a$, $b$, and $c \in
V$}.
\]

\begin{xca}
Define the \emph{Poisson operad}, using a tree model similar to the
previous examples. Show that an algebra over it is nothing but a
Poisson algebra. [\emph{Hint}: Use two kinds of vertices, one for the
dot product and the other one for the bracket.]
\end{xca}

\subsubsection{The Riemann surface operad and vertex operator algebras}

\begin{figure}
\centerline{\includegraphics[height=1.8cm]{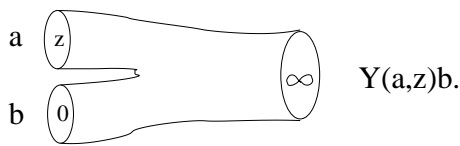}}
\caption{The vertex operator}
\label{voapants}
\end{figure}

Just for a change, let us return to the operad $\PP$ of Riemann
surfaces, more exactly, isomorphism classes of Riemann spheres with
holomorphic holes. What is an algebra over it? Since there are
infinitely many nonisomorphic pairs of pants, there are infinitely
many (at least) binary operations. In fact, we have an infinite
dimensional family of binary operations parameterized by classes of
pairs of pants. However modulo the unary operations, those which
correspond to cylinders, we have only one fundamental binary operation
corresponding to a fixed pair of pants. An algebra over this operad
$\PP$ is part a CFT data. (For those who understand, this is the tree
level, central charge $c=0$ part). If we consider a holomorphic
algebra over this operad, that is, require that the defining mappings
$\PP(n) \to \End{V}(n)$, where $V$ is a complex vector space, be
holomorphic, then we get part of a chiral CFT, or an object which
might have been called a \emph{vertex operator algebra $($VOA$)$} in
an ideal world. This kind of object is not equivalent to what people
use to call a VOA; according to Y.-Z. Huang's Theorem, a VOA is a
holomorphic algebra over a ``partial pseudo-operad of Riemann spheres
with rescaling,'' which is a version of $\PP$, where the disks are
allowed to overlap. The fundamental operation $Y(a,z)b$ for $a, b \in
V$, $z \in \nc$ of a VOA is commonly chosen to be the one
corresponding to a pair of pants which is the Riemann sphere with a
standard holomorphic coordinate and three unit disks around the points
0, $z$, and $\infty$ (No doubt, these disks overlap badly, but we
shrink them on the figure to look better), see Figure~\ref{voapants}.

\begin{figure}
\centerline{\includegraphics[height=2.5cm]{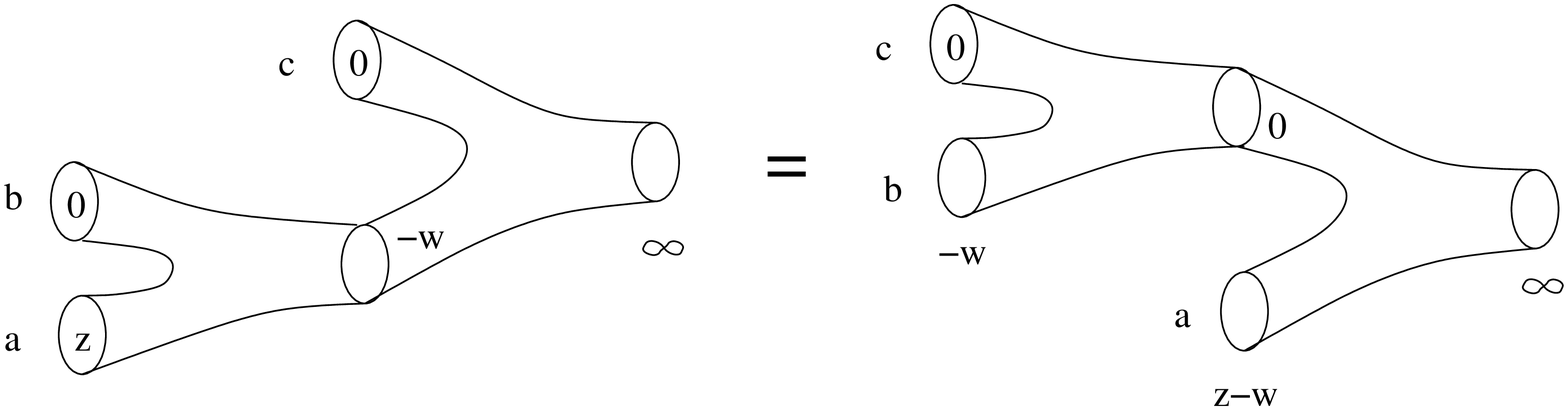}}
\caption{VOA associativity}
\label{voass}
\end{figure}

The famous associativity identity 
\[
Y(a,z-w)Y(b,-w) c = Y(Y(a,z)b,-w) c
\]
for vertex operator algebras comes from the natural
isomorphism of the Riemann surfaces sketched on Figure~\ref{voass}.

Another remarkable feature of the Riemann surface operad $\PP$ is that
an algebra over it in the category of spaces group completes to an
infinite loop space, which is Tillmann's result \cite{tillmann}. This
implies, for example, that the classifying space $B\Gamma^+_\infty$ of
the stable mapping class group (which is morally the moduli space of
Riemann surfaces of infinite genus) is an infinite loop space.

\subsubsection{The little disks operads and GBV-algebras}
\label{GBV}

\begin{figure}
\centerline{\includegraphics[width=1.2in]{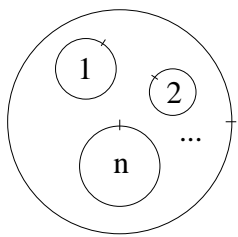}}
\caption{Framed little disks}
\label{fdisks}
\end{figure}

Here is another construction, related to the operad $\PP$ of Riemann
surfaces and more relevant to string topology. We will be talking
about a finite-dimensional retract of $\PP$, the \emph{framed little
disks operad} $f\DD$, which may be defined as follows, see Getzler
\cite{getzler:BV} and Markl-Shnider-Stasheff \cite{mss}. It is based on the
collection $f\DD = \{ f\DD(n) \; | \; n \ge 1\}$ of configuration
spaces of $n$ labeled nonintersecting (closed) disks in the standard
(closed) unit disk $D^2$ ain the plane $\nr^2$ with the choice of
marked point, thought of as framing, on the boundary of each
``little'' disk, see Figure~\ref{fdisks}.
It is also convenient to choose a marked point, e.g., $(1,0)$, on the
boundary of the ``big,'' unit disk, which should not be thought of as
extra data, because we are talking about the standard plane $\nr^2$,
with fixed $x$ and $y$ coordinates. An identity element is the framed
little disk coinciding with the big disk, together with framing. The
symmetric group acts by relabeling the framed little disks, as usual.
The operad composition $\circ_i: f\DD(m) \times f\DD(n) \to
f\DD(m+n-1)$ takes a given configuration of $n$ little disks, shrinks
it to match the size of the $i$th little disk in a given configuration
of $m$ little disks, rotates the shrunk configuration of $n$ little
disks (by a unique element of $\SO(2)$) to match the point on the
boundary of the big disk with the point on the boundary of the $i$th
little disk, and glues the configuration of $n$ disks in place of the
$i$th disk, erasing the seam afterwards. A
\emph{nonframed version}, which may be defined as the suboperad of
$f\DD$ with all the points on the boundaries of the little disks point
in the direction of the positive $x$ axis, is called the \emph{little
disks operad} $\DD$.

The little disks operad (along with its higher dimensional version)
was invented in topology and proved itself as a powerful tool for
studying iterated (based) loop spaces, \cite{may,boardman-vogt}. For
example, the little disks operad acts on every based double loop space
$\Omega^2 X = \Map_* (S^2, X)$ in the following way. Given a
configuration of $n$ little disks and $n$ pointed maps $S^2 \to X$,
which we can think of maps $D^2 \to X$ sending the boundary of the
standard unit disk $D^2$ to the basepoint of $X$, we can define a new
map $D^2 \to X$ by using the given maps on the little disks (after
appropriate translation and dilation) and extending them to a constant
map from the complement of the little disks to $X$. Of course, the
first thing one looks at in topology is homology, and the following
description of $H_* (\DD)$-algebras by F.~Cohen is very interesting.

\begin{thm}[F.~Cohen \cite{fcohen:disks}]
An algebra over the homology little disks operad $H_* (\DD; k)$ (over
a field $k$ of $\Char k \ne 2$) is equivalent to a \emph{Gerstenhaber
(}or simply \emph{G-) algebra}, i.e., a graded vector space $V$ with a
unit element $e$, a dot product $ab$, and a bracket $[a,b]$ defined,
so that the dot product defines the structure of a graded commutative
associative unital algebra, the bracket defines the structure of a
graded Lie algebra on the suspension $V[-1]$, which is the same as $V$
but with a grading shifted by $-1$, and the bracket is a degree-one
derivation of the dot product:
\[
[a,bc] = [a,b] c + (-1)^{\abs{a+1} \abs{b}} b [a,c] \qquad \text{for
all $a$, $b$, and $c \in V$}.
\]
\end{thm}

Note that if you have an algebra $X$ over an operad $\OO$ in the
category of topological spaces, you may always pass to homology and
obtain the structure of $H_*(\OO)$-algebra on the graded vector space
$H_* (X)$. Thus, the homology of every double loop space becomes a
G-algebra.

Similarly, the framed little disks operad acts on every based double
loop space, but the underlying structure on homology had not been much
of interest to topologists till it was discovered in various physical
contexts much later under the name of a BV-algebra.  A
\emph{BV-algebra} is a graded vector space $V$ (over a field of $\Char
\ne 2$) with the structure of a graded commutative algebra and a
second-order derivation $\Delta$, called a \emph{BV operator}, of
degree one and square zero. The second-order derivation property may
be defined using an old idea of Grothendieck:
\[
[[[\Delta, L_a],L_b],L_c] = 0 \quad \text{ for all } a, b, c \in V,
\]
where $L_a$ is the operator of left multiplication by $a$ and the
commutators of operators are understood in the graded
sense. Alternatively, see \cite{schwarz,getzler:BV}, one can define
the structure of a BV-algebra on a graded vector space $V$ in the
following way:
\begin{itemize}

\item
A ``dot product'' $V \otimes V \to V$ and a bracket $V \otimes V \to
V[1]$ of degree one making on $V$ making it a G-algebra;

\item
An operator $\Delta: V \to V[1]$ of degree one, which is a
differential, i.e., $\Delta^2 = 0$, a degree-one derivation of the
bracket, and satisfies the property
\[
\Delta (ab) - (\Delta a) b - (-1)^{\abs{a}} a \Delta b =
(-1)^{\abs{a}-1} [a,b].
\]
\end{itemize}

The BV-algebra structure is also known as the algebraic structure
induced on homology from the structure of an algebra over the framed
little disks operad $f\DD^2$ on a topological space:
\begin{thm}[Getzler \cite{getzler:BV}]
\label{getzler}
The category of $H_*(f\DD; k)$-algebras over a field $k$ of $\Char \ne
2$ is naturally isomorphic to the category of BV-algebras.
\end{thm}
In particular, the homology of a based double loop space $\Omega^2 X$
is naturally a BV-algebra.

\begin{figure}
\centerline{\includegraphics[width=1.2in]{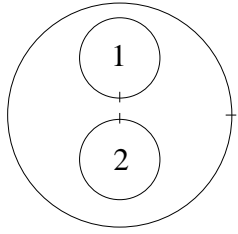}}
\caption{The dot product}
\label{dot}
\end{figure}

The proof of this theorem involves identifying the dot product and the
BV operator for a given $H_*(f\DD; k)$-algebra $V$. This is done as
follows. Look at the $n=2$ part
\[
H_*(f\DD(2);k) \to \Hom(V \otimes V, V)
\]
of the operad action. Note that $f\DD(2)$ is path connected and take
the class of a point in $H_0 (f\DD(2);k)$, for example, as on
Figure~\ref{dot}.
Define the dot product on $V$ as the resulting bilinear map $V \otimes
V \to V$.

To define the BV operator, consider the $n=1$ part
\[
H_*(f\DD(1);k) \to \Hom(V, V)
\]
of the operad action. Note that $f\DD(1)$ is homotopy equivalent to
$S^1$, choose an orientation on $S^1$, for example, the
counterclockwise rotation of the marked point on the little disk, and
take the fundamental class in $H_1 (S^1; k) \cong H_1 (f\DD(1);k)$.
The resulting linear map $V \to V$ is the BV operator, by definition.

To complete this proof, one needs to check that the topology of the
framed little disks operad forces the right identities between these
basic operations, as well as to show that there are no other
identities, except those coming from combining the identities in the
definition of a BV-algebra. We know how to do that, but we choose not
to do it in public.

\subsubsection{The little $n$-disks operad and $n$-fold loop spaces}
\label{n-algebras}

Generalizing the little disks operad of the previous section to an
arbitrary dimension $n \ge 1$, we get the little $n$-disks operad (an
equivalent version of which is known as the little $n$-cubes operad)
$\DD^n$ whose $k$th component $\DD^n (k)$ is the configuration space
of $k$ labeled ``little'' $n$-dimensional disks inside the standard
unit $n$-disk in $\nr^n$. This operad was introduced by Boardman and
Vogt \cite{boardman-vogt} and May \cite{may} to recognize $n$-fold
loop spaces $\Omega^n X = \Map_* (S^n, X)$ among other spaces: their
recognition principle states that a path-connected topological space
is weakly equivalent to an $n$-fold loop space, if and only if it is
weakly equivalent to an algebra over the little $n$-disks operad.

Passing to homology, F.~Cohen also proved that the notion of an
algebra over the homology little $n$-disks operad (over a field of
characteristic other than two) is equivalent to the notion of an
$n$-\emph{algebra}, which for $n \ge 2$ is the same as that of a
G-algebra, except that the bracket must now have degree $n-1$. A
1-algebra is the same as an associative algebra.

There is also a framed version $f\DD^n$ of the little $n$-disks
operad. Here a frame in a little $n$-disk is a positively oriented
orthonormal frame attached to the center of that disk. In other words,
each little disk comes with an element of $\SO (n)$. The operad of the
rational homology of $f\DD^n$ was characterized by the following
theorem of Salvatore and Wahl, which will be used in the last chapter
on brane topology. For better compatibility with that chapter, we are
going to quote this result for $n+1$ rather than $n$.

\begin{thm}[Salvatore-Wahl \cite{salvatore-wahl}: Theorem 6.5]
\label{sw:BV}
  The category of algebras over the operad $H_\bullet (f\DD^{n+1}; \nq)$
  for $n \ge 1$ is isomorphic to the category of graded vector spaces
  $V$ over $\nq$ with the following operations and identities. Below
  $a,b,c \in V$ are homogeneous elements and $\abs{a}$ denotes the
  degree of $a$ in $V$.
 \begin{enumerate}
   
 \item A dot product $a \cdot b$, or simply $ab$, defining the
   structure of a \textup{(}graded\textup{)} commutative associative
   algebra on $V$.
   
 \item A bracket $[a,b]$ of degree $n$, defining the structure of a
 \textup{(}graded\textup{)} Lie algebra on the shifted space $V[n]$.
   
 \item The bracket with an element $a$ must be a
 \textup{(}graded\textup{)} derivation of the dot product,
 \emph{i.e}., $[a,bc] = [a,b]c + (-1)^{(\abs{a}+n)\abs{b}} b[a,c]$.
   
 \item For $n$ odd, a collection of unary operators $\Beta_i$, $i =1,
 \dots, (n-1)/2$, of degree $4i-1$ and $\Delta$ of degree $n$, called
 a \emph{BV operator.}
   
 \item For $n$ even, a collection of unary operators $\Beta_i$, $i =1,
 \dots, n/2$, of degree $4i-1$.
   
 \item The unary operators $\Beta_i$ must square to zero: $\Beta_i^2 =
 0$ for all $i$. The operators $\Beta_i$ must be
 \textup{(}graded\textup{)} derivations of the commutative algebra
 structure on $V$ and the Lie algebra structure on $V[n]$.
   
 \item For $n$ odd, \emph{i.e}., when $\Delta$ is defined, $\Delta^2
   = 0$, $\Delta (ab) - \Delta (a) b - (-1)^{\abs{a}} a \Delta b =
   (-1)^{\abs{a}} [a,b]$, and $\Delta [a,b] = [\Delta a, b] -
   (-1)^{\abs{a}} [a, \Delta b]$.

 \end{enumerate}

\end{thm}

\begin{rem}
  The last identity means that $\Delta$ is a graded derivation of the
  bracket, while the last two equations may be interpreted as $\Delta$
  being a graded second-order derivation of the dot product.
\end{rem}

\begin{proof}[Idea of proof]
  This theorem is based on an observation that $f\DD^{n+1}$ is a
  semidirect product $\DD^{n+1} \rtimes \SO(n+1)$, see \cite{mss}. The
  elements $B_i$ (and $\Delta$, when it is defined) are the standard
  generators of $H_\bullet (\SO(n+1); \nq)$, $\Delta$ corresponding to
  the Euler class via the transfer map.
\end{proof}

\subsection{Operads via generators and relations}

The tree operads that we looked at above, such as the associative and
the Lie operads, are actually operads defined by generators and
relations. Here is a way to define such operads in general. To fix
notation, assume throughout this section that we work with operads
$\OO(n)$, $ n \ge 1$, of vector spaces.

\begin{definition}
An \emph{ideal} in an operad $\OO$ is a collection $\II$ of
$\Si_n$-invariant subspaces $\II(n) \subset \OO(n)$, for each $n \ge
1$, such that whenever $i \in \II$, its operad composition with
anything else is also in $\II$.
\end{definition}

The intersection of an arbitrary number of ideals in an operad is also
an ideal, and one can define the ideal generated by a subset in $\OO$
as the minimal ideal containing the subset.

\begin{definition}
For an operad ideal $\II \subset \OO$, the \emph{quotient operad}
$\OO/\II$ is the collection $\OO(n)/\II(n)$, $n \ge 1$, with the
structure of operad induced by that on $\OO$.
\end{definition}

\begin{figure}
\centerline{\includegraphics[height=0.5in]{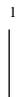}}
\caption{The operad identity}
\label{id}
\end{figure}

The \emph{free operad $F(S)$ generated by a collection $S = \{ S(n) \;
  | \; n \ge 1 \}$ of sets}, is defined as follows.
\[
F(S)(n) = \bigoplus_{\text{$n$-trees $T$}} k \cdot S(T),
\]
where the summation runs over all planar rooted trees $T$ with $n$
labeled leaves and
\[
S(T) = \Map(v(T), S),
\]
the set of maps from the set $v(T)$ of vertices of the tree $T$ to the
collection $S$ assigning to a vertex $v$ with $\In(v)$ incoming edges
an element of $S(\In(v))$ (the edges are directed toward the root). In
other words, an element of $F(S)(n)$ is a linear combination of planar
$n$-trees whose vertices are decorated with elements of $S$. There is
a special tree with no vertices, see Figure~\ref{id}.
The component $F(S)(1)$ contains, apart from $S(1)$, the
one-dimensional subspace spanned by this tree.

The following data defines an operad structure on $F(S)$.

\begin{enumerate}

\item The identity element is the special tree in $F(S)(1)$ with no 
vertices.

\item The symmetric group $\Si_n$ acts on $F(S)(n)$ by relabeling the inputs.

\item The operad composition is given by grafting the roots of trees
to the leaves of another tree. No new vertices are created.

\end{enumerate}

\begin{definition}
  Now let $R$ be a subset of $F(S)$, \emph{i.e}., a collection of
  subsets $R(n) \subset F(S)(n)$. Let $(R)$ be the ideal in $F(S)$
  generated by $R$. The quotient operad $F(S)/(R)$ is called the
  \emph{operad with generators $S$ and defining relations $R$}.
\end{definition}

\begin{ex}
  The associative operad $\ass$ is the operad generated by a point $S
  = S(2) = \{\bullet\}$ with a defining relation given by the
  associativity condition, see Section~\ref{ass:operad}, expressed in
  terms of trees. Note that equation $S=S(2)$ implies that $S(n) =
  \varnothing$ for $n \ne 2$.
\end{ex}

\begin{ex}
  The Lie operad $\lie$ is the operad also generated by a point $S =
  S(2) = \{\bullet\}$ with defining relations given by the skew symmetry
  and the Jacobi identity, see Section~\ref{lie:operad}.
\end{ex}

\begin{ex}
  The Poisson operad is the operad also generated by a two-point set
  $S = S(2) = \{\bullet, \circ\}$ with defining relations given by the
  commutativity and the associativity for simple trees decorated only
  with $\bullet$'s, the skew symmetry and the Jacobi identity for
  simple trees decorated with $\circ$'s, and the Leibniz identity for
  binary 3-trees with mixed decorations, see
  Section~\ref{poisson:operad}.
\end{ex}

\section{The cacti operad}

The construction and results in this section have been announced in
\cite{me:ua}. The BV structure arising in string topology at the level
of homology comes from an action of a \emph{cacti operad} $\CC$ at the
motivic level, quite close to the category of topological spaces.

The $k$th component $\CC(k)$ of the cacti operad $\CC$ for $k \ge 1$
may be described as follows.  $\CC(k)$ is the set of labeled by
numbers 1 through $k$ tree-like configurations of parameterized
circles, called the \emph{lobes}, of varying (positive) radii, along
with the following data: (1) the choice of a cyclic order of
components at each intersection point and (2) the choice of a marked
point on the whole configuration along with the choice of one of the
circles on which this point lies. The last choice is essential only
when the marked point happens to be an intersection point. Here
``tree-like'' means that the dual graph of this configuration, whose
vertices correspond to the lobes and the intersection points thereof
and whose edges reflect the obvious incidence relation, is a tree. We
will refer to an element of $\CC(k)$ as a \emph{cactus}. A cactus
defines a \emph{pinching map} from the standard unit circle $S^1$ to
the cactus. The pinching map starts from the marked point in the
direction of the increasing parameter on the circle on which the
marked point lies and traces the whole cactus along the parameters on
the circles, jumping from one lobe to the next in the cyclic order at
the intersection points. This produces a map from the circle of
circumference $c$ equal to the total circumference of the cactus. To
get a map from the standard unit circle, first expand (or contract) it
to a circle of circumference $c$.

\begin{figure}
\centerline{\includegraphics[width=1.5in]{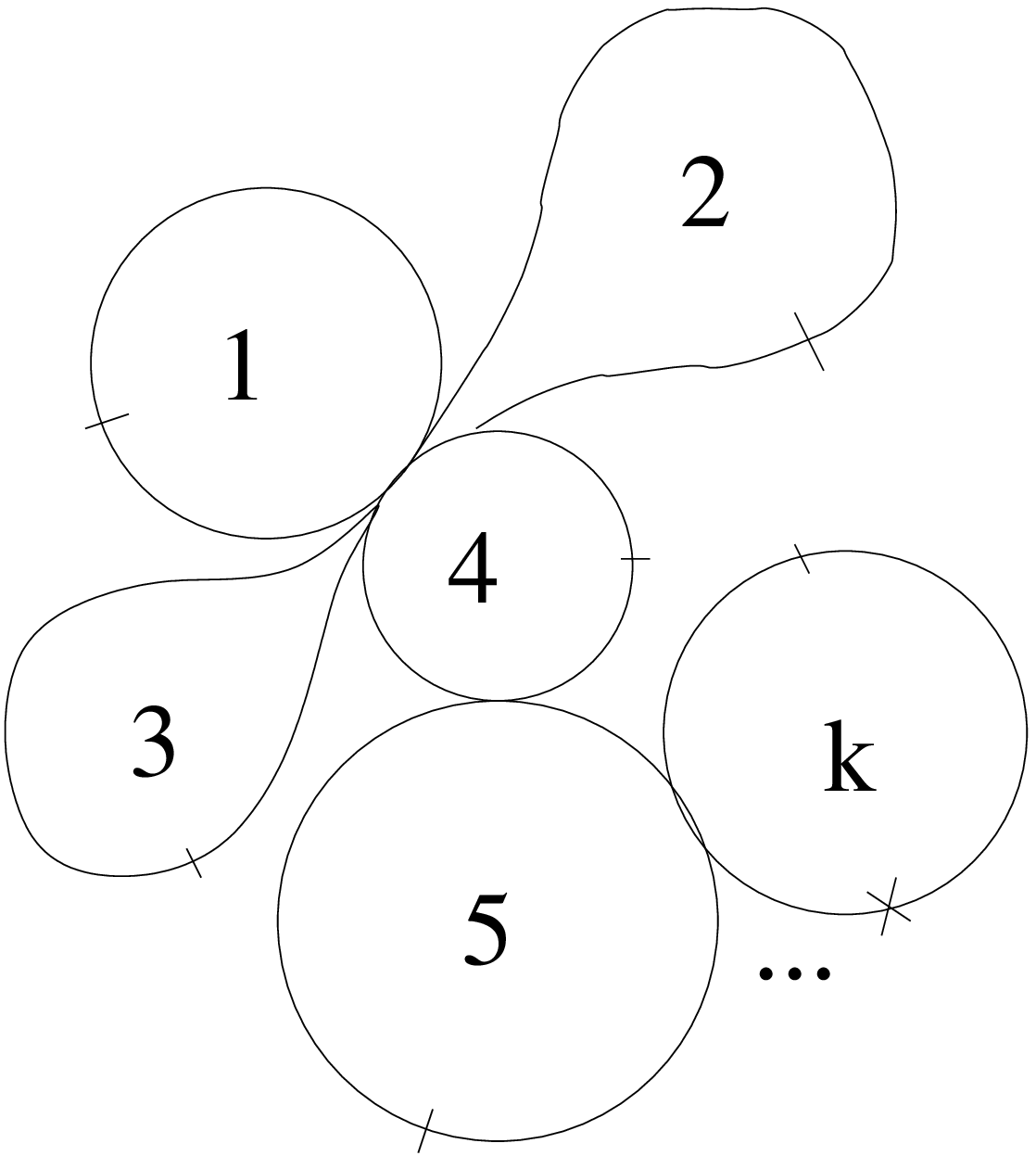}}
\caption{A cactus}
\label{cactus}
\end{figure}

The topology on the set $\CC(k)$ of cacti may be introduced in the
following way. There is a unique up to isotopy way to place a cactus
on the plane, so that the parameters of the lobes go counterclockwise
and the cyclic order of the lobes at each intersection point is also
counterclockwise, see Figure~\ref{cactus}.
Thereby a cactus defines a metric planar graph whose vertices are the
intersection points of the circles and edges are the arcs between two
adjacent vertices. The word ``metric'' refers to the fact that the
edges are provided with positive real numbers, the arclengths in the
parameters of the lobes. A cactus is determined by its metric planar
graph and the choice of a marked point on each of the $k$ circles
forming the cactus (called ``interior boundary components'' of the
graph) and a global marked point on the ``exterior boundary
component'' defined as the standard circle $S^1$ together with the
pinching map.

First of all, we will define a topology on the space of metric planar
graphs (and thereby on its subspace of metric planar graphs arising
from cacti) and then describe the space $\CC(k)$ of cacti as an
$(S^1)^{k+1}$ fiber bundle over the subspace of metric planar graphs.

Each planar graph $\Gamma$ defines an open cell $\nr_+^{e(\Gamma)}$,
where $\nr_+$ is the interval $(0,\infty)$ of the real line and
$e(\Gamma)$ is the set of edges of $\Gamma$. This open cell is
attached to the union of lower dimensional cells as follows. Consider
part of the boundary of $\nr_+^{e(\Gamma)}$ in $\nr^{e(\Gamma)}$ given
by setting some of the edge lengths to zero, except for the edges
forming simple loops, i.e., the edges which span the whole lobe less a
point. An attaching map identifies the face $l(e_0) = 0$ of the cell
$\nr_+^{e(\Gamma)} = \{ l(e) > 0 \; | \; e \in e(\Gamma)\}$ with a
cell $\nr_+^{e(\Gamma/e_0)}$, where $\Gamma/e_0$ is the planar graph
obtained from $\Gamma$ by contracting the edge $e_0$ to a point.  This
way we build the space of metric planar graphs as a generalization of
a cell complex: start with a collection of open cubes of dimension $k$
(corresponding to the planar wedges of $k$ circles whose starting
points and the global marked point are at the basepoint of the wedge),
attach a collection of open cubes of dimension $k+1$ to them along
some of their faces, then attach open cubes of dimension $k+2$ to the
result, and so on. We take the topology of the union on the resulting
space. This topologizes the space of metric planar graphs. The
``universal graph bundle'' over this space is determined by specifying
the fiber over a point to be the graph (thought of as a
one-dimensional CW complex) represented by this point in the space of
graphs. Marking a point on a specific (interior or exterior) boundary
component defines an $S^1$-bundle over the space of graphs. Thus,
$\CC(k)$ becomes a product of $k+1$ such bundles over the
corresponding subspace of metric planar graphs.

\begin{rem}
This construction identifies the space $\CC(k)$ of cacti as a certain
subspace of the space of metric ribbon graphs of genus zero with $k+1$
labeled boundary components and one point marked on each boundary
component. See Chapter~\ref{chapterthree} for a discussion of general
ribbon (fat) graphs.
\end{rem}

The \emph{operad structure} on the cacti comes from the pinching
map. Given two cacti and the $i$th lobe in the first one, the operad
composition $\circ_i$ will be given by further gluing the $i$th lobe
of the first cactus (identified with the standard circle $S^1$ via a
suitable dilation) to the second cactus along the pinching map from
$S^1$ to the second cactus.

\begin{thm}[\cite{me:ua}]
\label{cacti}
The cacti operad $\CC$ is homotopy equivalent to the framed little
disks operad $f\DD^2$.
\end{thm}

\begin{rem}
\label{Salvatore}
In principle, one can prove this theorem explicitly, as indicated in
\cite{cohen-jones}. For example, P.~Salvatore (private communication)
suggests constructing a map from the configuration space of $k$
labeled points in the plane by placing at each point a particle, which
creates a radial repulsive field of magnitude $1/r$, where $r$ is the
distance from the particle, and looking at the degenerate trajectories
of the superposition field. This is, in fact, dual to the construction
of a ribbon graph on the Riemann sphere with $k+1$ punctures using a
Strebel differential and its degenerate horizontal trajectories:
Salvatore's approach uses the vertical ones.
\end{rem}

We will take a less constructive approach and deduce the statement
from Salvatore-Wahl's recognition principle \cite{salvatore-wahl} for
the framed little disks operad $f\DD^2$, which generalizes
Fiedorowicz's recognition principle \cite{fied-symm,fied-constr} for
the little disks operad $\DD^2$. We will use the recognition principle
in the following form. Let $B_k$ be the braid group on $k$ strands, $k
\ge 1$, then the ribbon braid group $\RB_k$ on $k$ ribbons is a
semidirect product $\nz^k \rtimes B_k$ determined by the permutation
action of $B_k$ on $\nz^k$. It is worth pointing out that the ribbon
braid group $\RB_k$ is a subgroup of $B_{2k}$ by the induced braid on
the edges of the ribbons. Let $\PRB_k$ denote the kernel of the
natural epimorphism $\RB_k \to \Si_k$. Then $\PRB_k \cong \nz^k \times
\PB_k$, where $\PB_k$ is the pure braid group, the kernel of the
natural epimorphism $B_k \to \Si_k$.

\begin{thm}[Salvatore-Wahl \cite{salvatore-wahl}]
\label{sw}
Let $\OO = \{ \OO(k) \; | \; k \ge 1\}$ be a topological operad,
satisfying the following conditions:
\begin{enumerate}

% \item
% The symmetric group $\Si_k$ acts freely on $\OO(k)$ for each $k \ge
% 1$.

\item
Each quotient $\OO(k)/\Si_k$ is $K(\RB_k,1)$, the normal subgroup
$\PRB_k \subset \RB_k$ corresponding to the covering $\OO(k) \to
\OO(k) /
\Si_k$.

\item
There exists an operad morphism $\DD^1 \to \OO$ from the little
intervals operad $\DD^1$.

\end{enumerate}
Then the operad $\OO$ is homotopy equivalent to the framed little
disks operad $f\DD^2$.
\end{thm}

\begin{rem}
If we require in Condition 2 that a cofibrant model of the little
intervals operad $\DD^1$ admits a morphism to $\OO$, the two
conditions become not only sufficient, but also necessary. Examples of
cofibrant models of $\DD^1$ include the operad of metrized planar
trees with labeled leaves (a cellular operad, whose chain operad is
$A_\infty$) and Boardman-Vogt's $W$-resolution $W\DD^1$ of $\DD^1$.
\end{rem}

\begin{rem}
This recognition principle can be generalized to any non-$\Si$ operad
of groups $H(k)$ with $H(1)$ abelian and the operad composition and
the operad unit for $H(1)$ given by the group law and the group unit
therein, respectively. This operad of groups should also be provided
with an operad epimorphism $H(k) \to \Si_k$. $\{\Si_\bullet\}$,
$\{B_\bullet\}$, and $\{\RB_\bullet\}$ are examples of such operads of
groups. The proof of the recognition principle is the same.
\end{rem}

\begin{proof}
We will only sketch a proof, following Salvatore and Wahl: the reader
may collect missing details from
\cite{fied-symm,fied-constr,mcs,salvatore-wahl}. One first introduces
the notion of a ribbon braided operad, in which the symmetric groups
get replaced with the ribbon braid groups. Condition 2 amounts to a
consistent choice of basepoints in the universal coverings $\widetilde
\OO(k)$, which gives the structure of a ribbon braided operad on
$\widetilde \OO$. Condition 1 implies that each space $\widetilde
\OO (k)$ is contractible with a free action of $\RB_k$. Then the
following diagram of homotopy equivalences of ribbon braided operads
\[
\widetilde{f\DD^2} \leftarrow \widetilde{f\DD^2} \times \widetilde{\OO}
\to \widetilde{\OO},
\]
after taking quotients by the pure ribbon braid groups $\PRB_k$, gives
homotopy equivalences of usual operads:
\[
f\DD^2 \leftarrow (\widetilde{f\DD^2} \times \widetilde \OO)/\PRB_\bullet
\rightarrow \OO.
\]
%\noindent
\end{proof}

\begin{proof}[Proof of Theorem~\textup{\ref{cacti}}]
As promised, we will use Salvatore-Wahl's recognition principle,
Theorem~\ref{sw}.

To see that $\Si_k$ acts freely on $\CC(k)$, note that the structure
map from $S^1$ to the cactus defines an ordering on the lobes.

There is an obvious homomorphism $\RB_k \to \pi_1 (\CC(k)/\Si_k)$ in
which a ribbon braid moves in the plane the lobes of a fixed cactus,
\emph{e.g}., one like on Figure~\ref{standard},
\begin{figure}
\centerline{\includegraphics[width=1.3in]{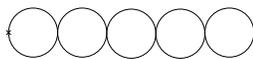}}
\caption{The basepoint in the cacti operad}
\label{standard}
\end{figure}
keeping the marked point somewhere on the left to return it to the
original position at the end of the move. Here is an example of three
intermediate snapshots of the above cactus through a move, with the
corresponding braid on the side, see Figure~\ref{move-c}.

\begin{figure}
\centerline{\includegraphics[width=2.1in]{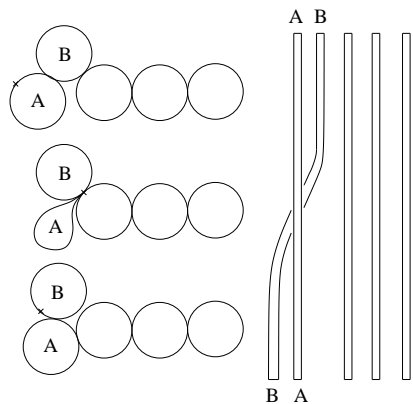}}
\caption{A path in the cacti space corresponding to a braid}
\label{move-c}
\end{figure}

We will prove that $\CC(k)/\Si_k$ is $K(\RB_k,1)$ by induction on
$k$. For $k=1$ the space $\CC(1)/\Si_1 = \CC(1)$ is obviously
homeomorphic to $S^1 \times \nr_+$, which is $K(\nz,1)$. Since $\RB_1
\cong \nz$, we get the induction base.

To make the induction step, assume that the above homomorphism $\RB_k
\to \pi_1 (\CC(k)/\Si_k)$ is bijective and consider a forgetful
fibration $\CC(k+1) \to \CC(k)$ which contracts the $k+1$st lobe to
a point, \emph{i.e}., takes a quotient of the cactus by its $k+1$st
lobe. The marked point, if it happened to be on that lobe, will be at
the resulting contraction point, assigned to the lobe that follows the
one being contracted, if one moves from the marked point along the
parameter of the wrapping map from $S^1$ to the cactus. It is easy to
observe that the fiber of this fibration is homotopy equivalent to
$\left(\bigvee_n S^1\right) \times S^1$, where the last circle is
given by the angular parameter on the contracted circle.

This fibration admits a section, for example, by placing a unit circle
at the marked point, so that the wrapping map will trace this circle
completely from its angular parameter $0$ and at $2\pi$ bumps into the
marked point.

Therefore, the long exact homotopy sequence of this fibration splits
into short exact sequences as follows:
\begin{gather*}
1 \to \pi_1 \left(\left({\textstyle\bigvee_k} S^1\right) \times S^1
\right) \to \pi_1(\CC(k+1)) \to \pi_1 (\CC(k)) \to 1,\\
0 \to 0 \to \pi_i(\CC(k+1)) \to \pi_i (\CC(k)) \to 0 \qquad
\text{for $i \ge 2$,}
\end{gather*}
with the arrow $\pi_1(\CC(k+1)) \to \pi_1 (\CC(k))$ in the first
sequence having a right inverse. By induction, $\CC(k)$ is
$K(\PRB_k,1)$. Then from the short exact sequences, we see that
$\pi_i(\CC(k+1)) = 0$ for $i \ge 2$.

To treat the fundamental-group case, note that we have a natural
morphism of group extensions:
\begin{equation*}
\begin{CD}
1 @>>> \pi_1 \left(\left(\bigvee_k S^1\right) \times S^1\right) @>>>
\PRB_{k+1} @>>> \PRB_k @>>> 1\\
@. @| @VVV @VVV \\
1 @>>> \pi_1 \left(\left(\bigvee_k S^1\right) \times S^1\right) @>>> 
\pi_1(\CC(k+1)) @>>> \pi_1 (\CC(k)) @>>> 1,
\end{CD}
\end{equation*}
where the first line is obtained by applying the fundamental group
functor to a similar forgetful fibration $f\DD(k+1) \to f\DD(k)$. The
right vertical arrow is an isomorphism by the induction assumption,
therefore, so is the middle one.

Now the following morphism of extensions:
\begin{equation*}
\begin{CD}
1 @>>> \PRB_{k+1} @>>> \RB_{k+1} @>>> \Si_{k+1} @>>> 1\\
@. @VVV @VVV @| \\
1 @>>> \pi_1(\CC(k+1)) @>>> \pi_1 (\CC(k+1)/\Si_{k+1}) @>>>
\Si_{k+1} @>>> 1
\end{CD}
\end{equation*}
shows that the middle vertical arrow is an isomorphism. This completes
the induction step and thereby verifies Condition~1.

Finally, let us check Condition~2. Given a configuration of $k$ little
intervals, sketch $k$ circles based on the little intervals as
diameters. Then contract the spaces between the circles and put the
marked point at the left corner of the first circle on the left, as
on Figure~\ref{intervals}.
\begin{figure}
\centerline{\includegraphics[width=2.1in]{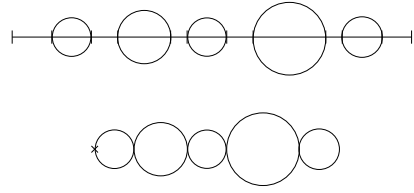}}
\caption{Little intervals to cacti}
\label{intervals}
\end{figure}
Amazingly enough, this defines a morphism of operads $\DD^1 \to \CC$.
\end{proof}

\section{The cacti action on the loop space}
\label{cacti:action}

Let $M$ be an oriented smooth manifold of dimension $d$. We would like
to study continuous $k$-ary operations on the \emph{free loop space}
$LM: = M^{S^1} := \Map (S^1, M)$ of continuous loops, functorial with
respect to $M$. By passing to singular chains or homology, these
operations will induce functorial operations on the chain and homology
level, respectively. We will generalize a homotopy-theoretic approach
developed by R.~Cohen and J.~D.~S. Jones \cite{cohen-jones} for the
loop product.

\subsection{Action via correspondences}

First of all, consider the following diagram
\begin{equation}
\label{rho}
\CC(k) \times (LM)^k \xleftarrow{\rho_{\varin}} \CC (k)M
\xrightarrow{\rho_{\out}} LM,
\end{equation}
where $\CC$ is the cacti operad, $k \ge 1$, and $\CC (k)M$ is the
space of pairs $(c,f)$ with $c \in \CC(k)$ being a cactus and $f: c
\to M$ a continuous map of the cactus $c$ to $M$. The map
$\rho_{\varin}: \CC (k)M \to \CC(k) \times (LM)^k$ takes a pair
$(c,f)$ to $c \in \CC(k)$ and the restrictions of $f$ to the $k$ lobes
of $c$. It is an embedding of codimension $d(k-1)$. The map
$\rho_{\out}: \CC (k)M \to LM$ composes the pinching map $S^1 \to c$
with $f: c \to M$. This diagram is, in principle, all you need to
define an operad ``action'', \emph{i.e}., the structure of a
``$\CC(k)$-algebra'' on the loop space $LM$. It would be an honest
operad action, should the arrow $\rho_{\varin}$ be pointing in the
opposite direction and the resulting composite map be compatible with
the operad structure on $\CC$. To check this compatibility, we indeed
need to find a way to invert $\rho_{\varin}$ and take the composite
map. The main idea is to take the motivic standpoint and treat
diagram~\eqref{rho} as a single morphism going from left to right.

We will consider the \emph{category $\Corr$ of correspondences}. The
objects of this category are just topological spaces, and a morphism
between two objects $X$ and $Y$ is a \emph{correspondence}, by which
here we mean a diagram $X \leftarrow X' \to Y$ of continuous maps for
some space $X'$, or, equivalently, a map $X' \to X \times Y$. One
composes two correspondences $X \leftarrow X' \to Y$ and $Y \leftarrow
Y' \to Z$ by taking a pullback
\begin{equation*}
\begin{CD}
X' \times_Y Y' @>>> Y'\\
@VVV @VVV\\
X' @>>> Y
\end{CD}
\end{equation*}
which defines a new correspondence $X \leftarrow X' \times_Y Y' \to Z$.

\begin{thm}
\label{homotopyaction}
\begin{enumerate}
\item
Diagram \eqref{rho}, considered as a morphism $\CC(k) \times (LM)^k
\to LM$ in $\Corr$, defines the structure of a $\CC$-algebra on the
loop space $LM$ in $\Corr$.
\item
This $\CC$-algebra structure on the loop space $LM$ in $\Corr$ induces
an $h_* (\CC)$-algebra structure on the shifted homology $h_{*+d}
(LM)$ for any multiplicative generalized homology theory $h_*$ which
supports an orientation of $M$.
\end{enumerate}
\end{thm}

\begin{proof}
\begin{enumerate}

\item
We just need to see that the following diagram commutes in $\Corr$:
\begin{equation}
\label{corr:comm}
\begin{CD}
\CC(k) \times \CC(l) \times (LM)^{k+l-1}
@>{\circ_i \times \id}>> \CC(k+l-1) \times (LM)^{k+l-1}\\
@VVV @VVV\\
\CC(k) \times (LM)^{k} @>>> LM,
\end{CD}
\end{equation}
where the unmarked arrows are correspondences defined using diagram
\eqref{rho}, the left vertical arrow is such a correspondence
preceded by an appropriate permutation, to make sure that $\CC(l)$
acts on the $l$ components $LM$ in $LM^{k+l-1}$, starting from the
$i$th one. To verify the commutativity of diagram~\eqref{corr:comm},
we compose the correspondences, using pullbacks, and see that both
compositions (the lower left one and the upper right one) are equal to
the following correspondence:
\[
\CC(k) \times \CC(l) \times (LM)^{k+l-1} \leftarrow \CC(k) \circ_i
\CC(l) M \to LM,
\]
where $\CC(k) \circ_i \CC(l) M := \{ (C_1, C_2, f) \; | \; C_1 \in
\CC(k), C_2 \in \CC(l), \text{ and } f: C_1 \circ_i C_2 \to M \text{
is continuous}\}$ and the maps are obvious.

\item
Verifying the statement at the homology level is a little subtler, as
not any correspondence induces a morphism on homology. Suppose we have
a correspondence $X \xleftarrow{e} X' \xrightarrow{\gamma} Y$ between
smooth (infinite dimensional) manifolds, so that $e$ is a regular
embedding of codimension $p$. (One can say that such a
\emph{correspondence is of degree} $-p$.) A \emph{regular embedding}
$X \hookrightarrow X'$ is, locally in $X'$ at $X$, the product of a
neighborhood in $X$ with a Euclidean space of codimension $p$. This
condition assures that the tubular neighborhood theorem applies. Then
we can apply the Thom collapse construction to $e$ and get a
composition
\[
h_* (X) \xrightarrow{e^!} h_{*-p}(X') \xrightarrow{\gamma_*} h_{*-p}(Y).
\]
Note that the inclusions $\CC(k) \circ_i \CC(l) M \hookrightarrow
\CC(k) \times \CC(l) M \times LM^{k-1} \hookrightarrow \CC(k) \times
\CC(l) \times (LM)^{k+l-1}$ participating in the lower left path on
diagram~\eqref{corr:comm} are regular embeddings of codimensions
$d(k-1)$ and $d(l-1)$, respectively, because they are pullbacks of
regular (finite-dimensional) embeddings along fiber bundles, while the
composite inclusion, which also shows up in the upper right path on
\eqref{corr:comm}, is a regular embedding of codimension
$d(k+l-2)$. Because of the functoriality of the homology with respect
to Thom collapse maps and the naturality of Thom collapse maps on
pullback diagrams with regular embeddings on two parallel sides, we
conclude that diagram~\eqref{corr:comm} induces a commutative diagram
in homology. This checks that $h_{*+d} (LM)$ is an algebra over the
operad $h_*(\CC)$.
\end{enumerate}

\end{proof}

\subsection{The BV structure}

String topology originated from Chas and Sullivan's construction of a
BV structure on the shifted homology of a loop space in a compact,
oriented manifold. A BV-algebra is nothing but another algebraic
structure, see Section~\ref{GBV}. However, its constant appearance
within different contexts of mathematical physics makes it worth a
little attention.

Combining Theorems \ref{cacti} and \ref{homotopyaction} with
Theorem~\ref{getzler} and checking what the basic operations (the dot
product and the BV operator) really are, as in the discussion of the
proof of Theorem~\ref{getzler}, we conclude with the following result.

\begin{crl}
For an oriented $d$-manifold $M$ and a field $k$ of characteristic
$\ne 2$, the space $H_{*+d}(LM; k) = H_*(LM; k)[d]$ has the natural
structure of a BV-algebra. This structure coincides with the one
constructed by Chas and Sullivan.
\end{crl}

%%%%%% Version of 3/21/05

\chapter[String topology as field theory]
{Field theoretic properties of string topology}
\label{chapterthree}

\section{Field theories}

\subsection{Topological Field Theories}

The first axiomatic definition of Topological (Quantum) Field Theories
was due to Atiyah \cite{atiyah:TQFT} and based on the work of Witten and
Segal. An $n$-\emph{dimensional Topological Field Theory} (\emph{TFT})
comprises the following data.\footnote{One usually puts in more data,
so what we are describing is a geometric background of a TFT.}
\begin{enumerate}

\item
An assignment: \{a closed oriented $(n-1)$-dimensional manifold $X$,
considered up to diffeomorphism\} $\mapsto$ \{a complex vector space
$V(X)$\}, a \emph{state space}, interpreted as the space of functions
(or other linearizations of a space, such as differential forms,
(singular) cochains, or chains) on a space of fields associated to
$X$. In the case of the so-called sigma-model, most relevant to string
topology, the space of fields is the loop space $LX$.

\item
An assignment: \{an $n$-dimensional oriented cobordism $Y$, considered
up to diffeomorphism, bounded by two $(n-1)$-manifolds $X_1$ and
$X_2$\} $\mapsto$ \{a linear operator $\Psi_Y: V(X_1) \to V(X_2)$\}.
One usually thinks of $X_1$ as the \emph{input} and $X_2$ as the
\emph{output} of cobordism $Y$.

\end{enumerate}
These data need to satisfy the following axioms:
\[
V(X_1 \coprod X_2) \xrightarrow{\sim} V(X_1) \otimes V(X_2)
\]
and $\Psi_Y$ must be compatible with the tensor products and
composition of cobordisms.

\begin{rem}
The categories of cobordisms and vector spaces are tensor (more
precisely, symmetric monoidal) categories and an $n$-dimensional TFT
is nothing but a tensor functor between them.
\end{rem}

\begin{rem}
For $n=2$, a closed oriented one-dimensional manifold $X$ is
diffeomorphic to a disjoint union $\coprod S^1$ of finitely many
circles $S^1$, interpreted as closed strings. In this case, the
category of cobordisms is a PROP, see Section~\ref{operads}, and a TFT
is a morphism from the PROP of compact oriented surfaces with boundary
to $\End{V}$, the endomorphism PROP of the vector space $V$. The PROP
$\End{V}$ may be identified with the tensor subcategory of the
category of vector spaces generated by $V$.
\end{rem}

\begin{thm}[Folklore]
A 2-dimensional TFT is equivalent to a complex Frobenius algebra $A$,
which is a finite-dimensional commutative associative algebra $A$ with
a unit and a nonsingular trace $\theta: A \to \nc$, i.e., a linear map
defining a nonsingular bilinear product $A \otimes A \to \nc$ via
$\theta (ab)$.
\end{thm}

\begin{rem}
This is a folk theorem, whose proof was passed on since the early
nineties by word of mouth by many folks, at least including Dijkgraaf,
Segal, and Witten, and was written down in the mathematical literature
by Abrams \cite{abrams} and Segal \cite{segal:stanford}.
\end{rem}

\begin{proof}[Outline of Proof] Let $V$ be
the state space of a TFT. We would like to construct the structure of
a Frobenius algebra on $V$. Consider the oriented surfaces and the
corresponding operators, which we denote $ab$, $\langle a ,b \rangle$,
\emph{etc}., as on Figure~\ref{Frob}.
\begin{figure}
\caption{Basic operations corresponding to surfaces}
\label{Frob}
\[
\begin{array}{ccll}
\parbox{2in}{\centerline{\includegraphics[height=.7in,width=1.5in]{pants.eps}}} &
 \quad \longmapsto \quad & V \otimes V \to V, & a \otimes b \mapsto ab
\smallskip\\
\parbox{2in}{\centerline{\includegraphics[width=.9in,height=.7in]{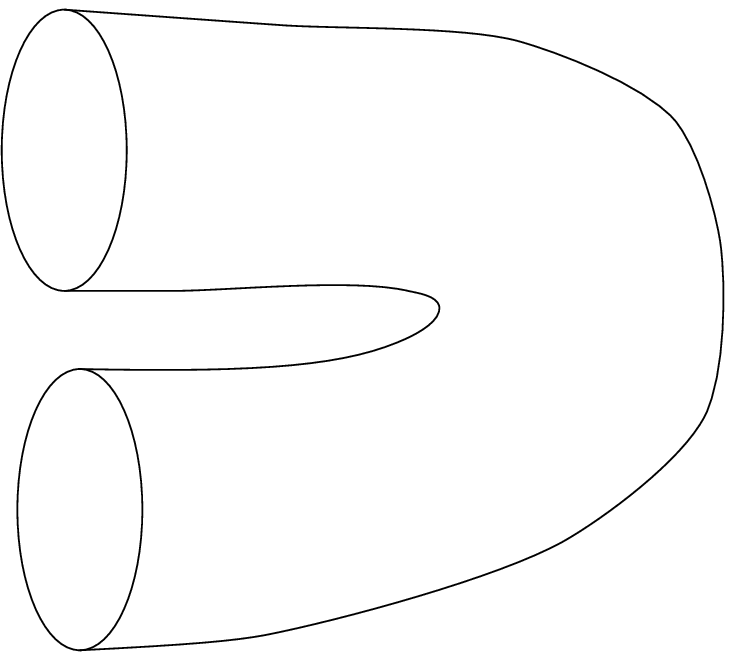}}} &
 \quad \longmapsto \quad & V \otimes V \to k, & a \otimes b \mapsto
\langle a,b \rangle
\medskip\\
\parbox{.5in}{\centerline{\includegraphics[width=-.9in,height=.7in]{bublik.eps}}} & \quad \longmapsto \quad & k
\to V \otimes V, & 1 \mapsto \psi
\medskip\\
\parbox{1.5in}{\centerline{\includegraphics[width=1.2in,height=.3in]{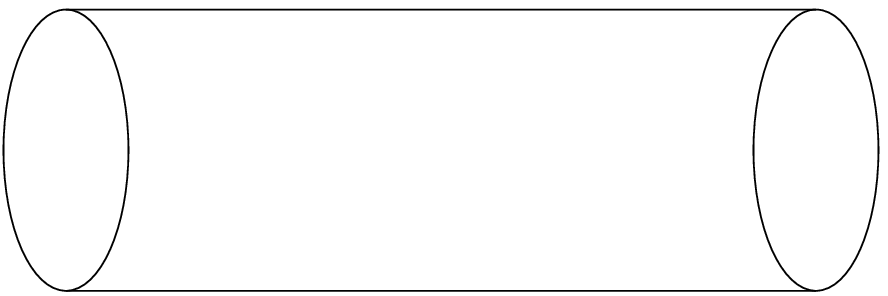}}} &
 \quad \longmapsto \quad & \id: V \to V
\medskip\\
\parbox{1.5in}{\centerline{\includegraphics[width=.8in,height=.4in]{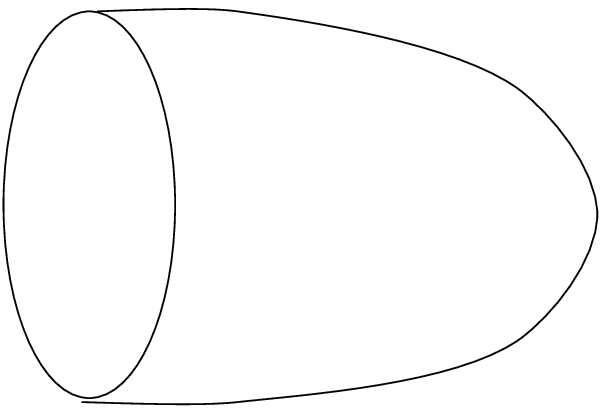}}} &
 \quad \longmapsto \quad & \tr: V \to k
\smallskip\\
\parbox{1.5in}{\centerline{\includegraphics[width=.7in,height=.4in]{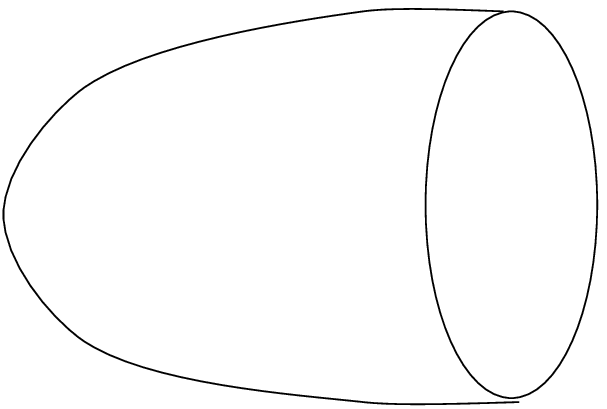}}} &
 \quad \longmapsto \quad & k \to V, & 1 \mapsto e
\end{array}
\]
\end{figure}

We claim that these operators define the structure of a Frobenius
algebra on $V$. Indeed, the multiplication is commutative, because if
we interchange labels at the legs of a pair of pants we will get a
homeomorphic oriented surface. Therefore, the corresponding operator
$a \otimes b \mapsto ba$ will be equal to $ab$. Similarly, the
associativity $(ab)c = a(bc)$ of multiplication is based on the fact
that the two surfaces on Figure~\ref{coherence} are homeomorphic:
\begin{figure}
\caption{Associativity}
\label{coherence}
\begin{equation*}
\parbox{7cm}{\centerline{\parbox{3.6cm}{\includegraphics[width=3.5cm]{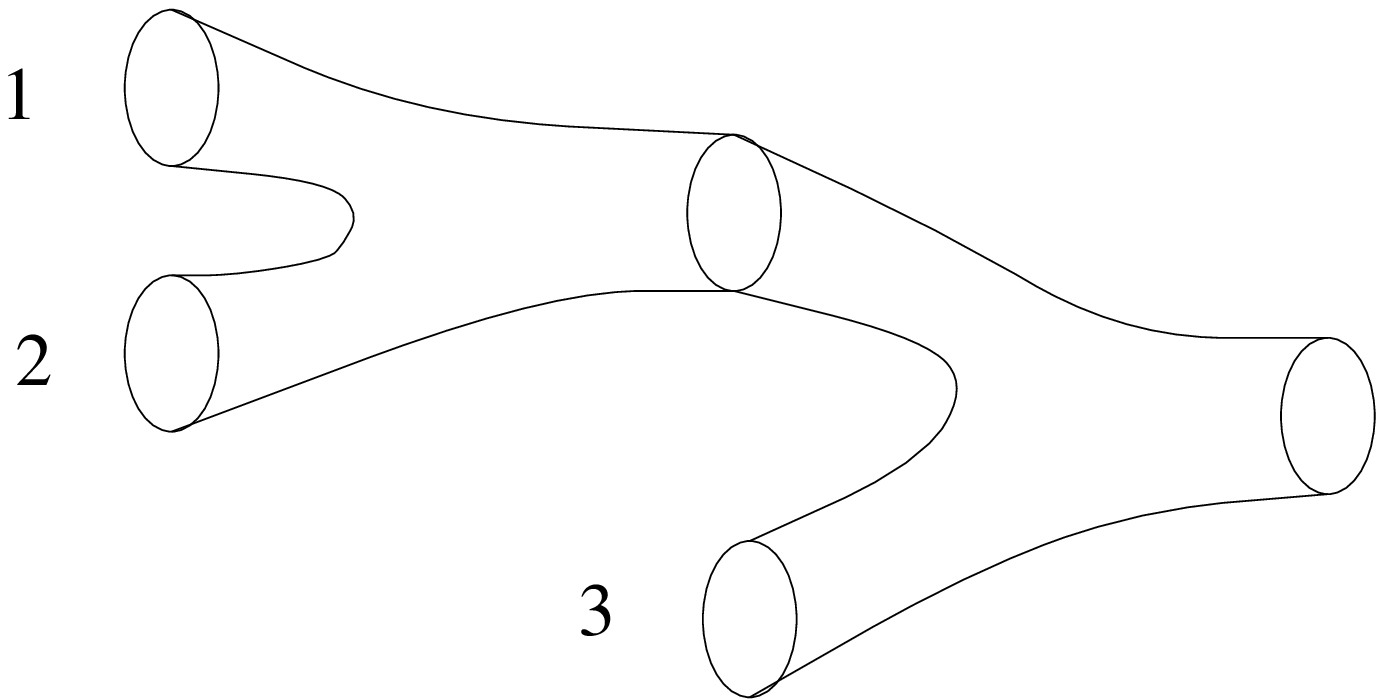}}}
$\cong$
\ \ \
\parbox{3.5cm}{\includegraphics[width=3.5cm]{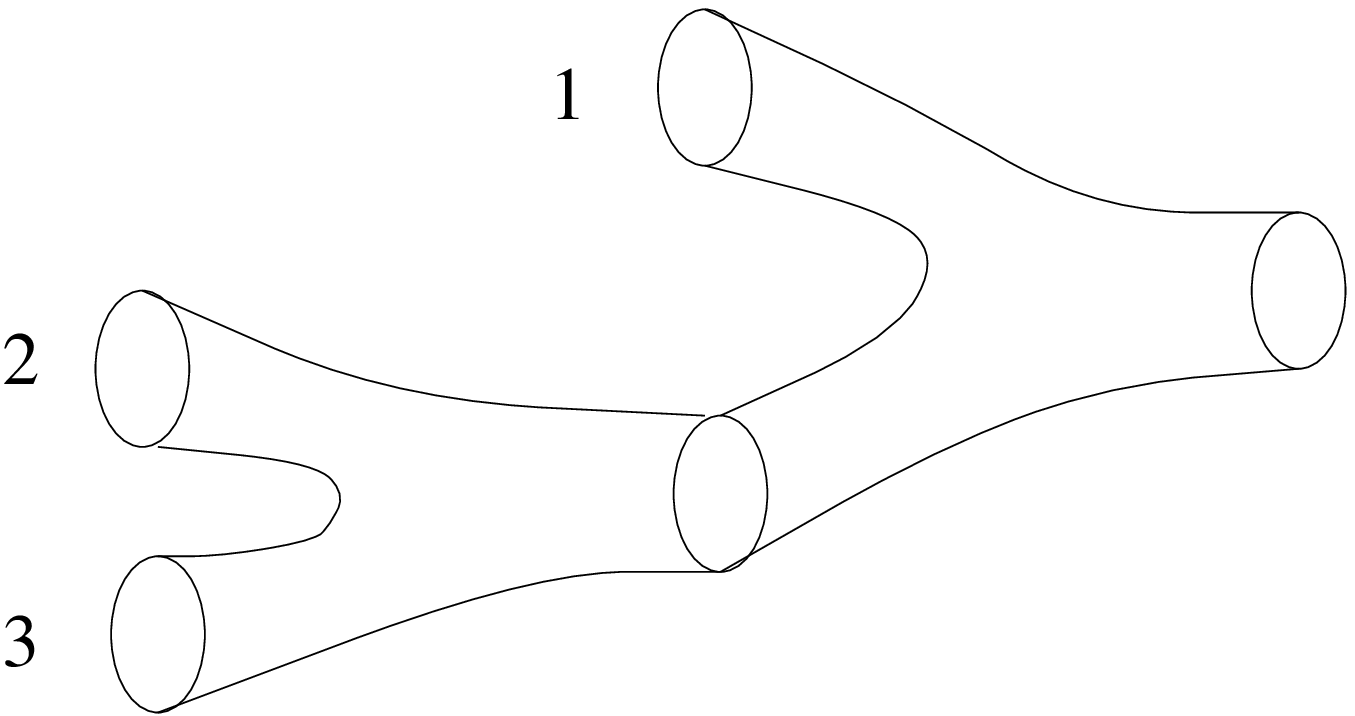}}}
\end{equation*}
\end{figure}
The property $ae = ea = a$ of the unit element comes from the
diffeomorphism on Figure~\ref{unit}.
\begin{figure}
\caption{The unit axiom}
\label{unit}
\centerline{\parbox{3cm}{\includegraphics[width=2.5cm]{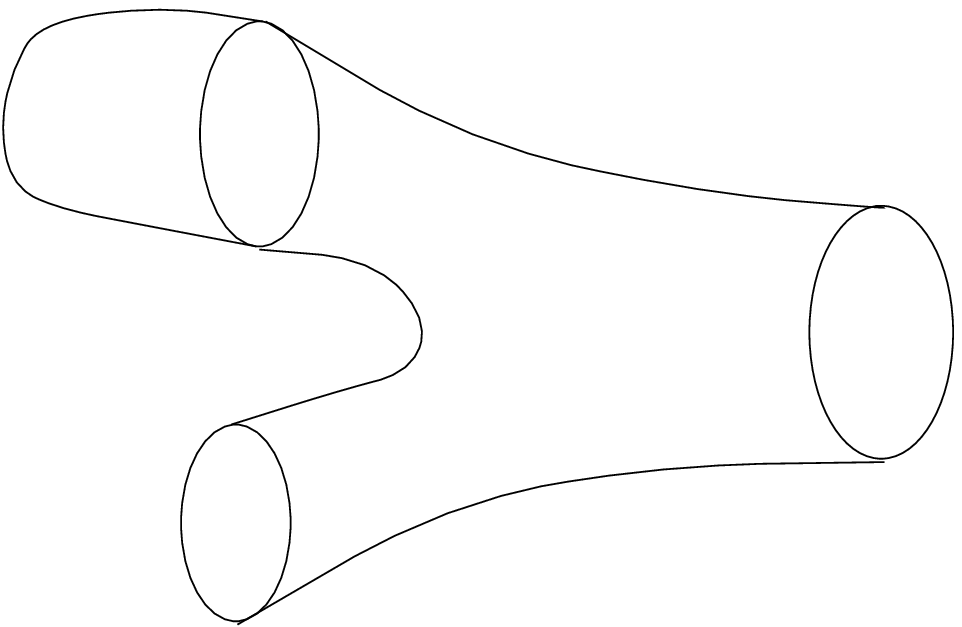}} $\cong$
\ \ \ \ \ \ 
\parbox{3.5cm}{\includegraphics[width=1.8cm]{normal.eps}}}
\end{figure}
Thus, we see that $V$ is a commutative associative unital
algebra. Now, the diffeomorphism on Figure~\ref{trace}
\begin{figure}
\caption{The trace}
\label{trace}
\centerline{\parbox{3cm}{\includegraphics[width=1.5cm]{bublik.eps}} $\cong$
\ \ \ \ \ \ 
\parbox{3.5cm}{\includegraphics[width=1.8cm]{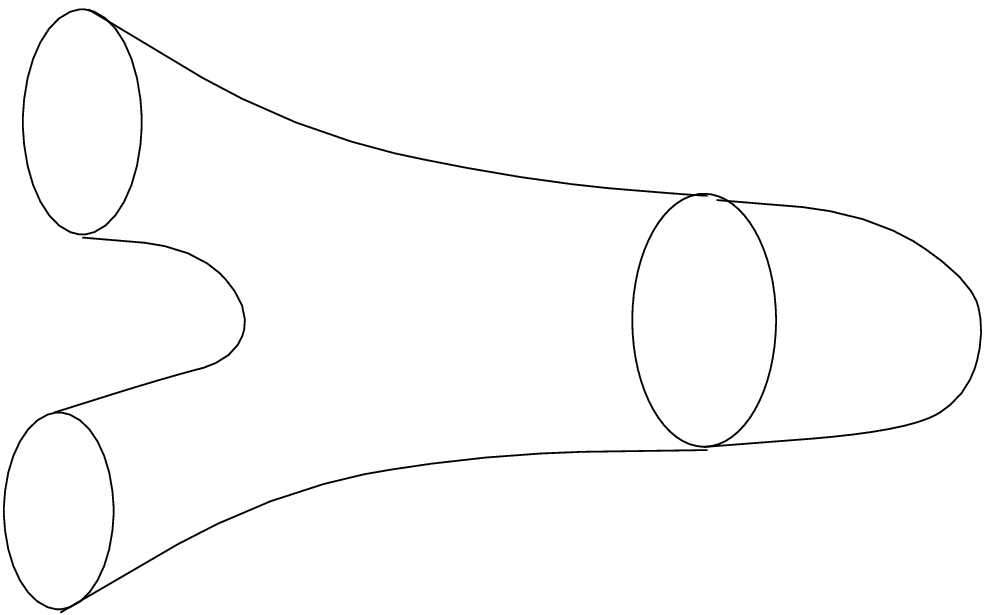}}}
\end{figure}
proves the identity $\langle a , b \rangle = \tr (ab)$, which, along
with the associativity, implies $\langle ab,c \rangle = \langle a,bc
\rangle$. The fact that the inner product $\langle \, , \, \rangle$ is
nondegenerate follows from the diffeomorphism on Figure~\ref{nondeg},
\begin{figure}
\caption{Nondegeneracy}
\label{nondeg}
\centerline{\parbox{3cm}{\includegraphics[width=2cm]{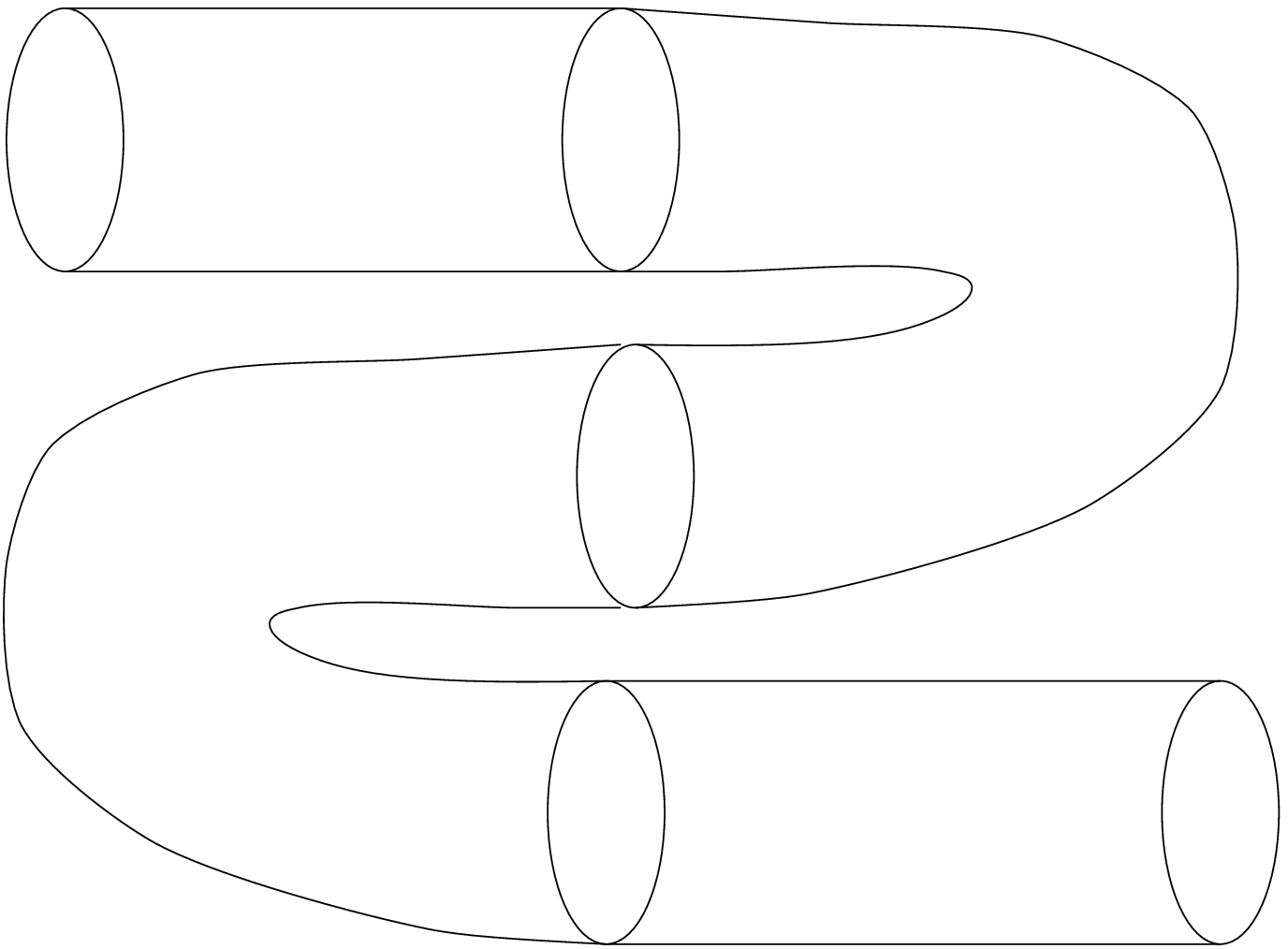}} $\cong$
\ \ \ \ \ \ 
\parbox{2.5cm}{\includegraphics[width=1.8cm]{normal.eps}}}
\end{figure}
which implies that the composite mapping
\[
\begin{array}{ccccc}
V & \xrightarrow{\id\otimes \psi} & V \otimes V \otimes V &
\xrightarrow{\langle,\rangle \otimes \id} &V,\\
v & \longmapsto & \sum_{i=1}^n v \otimes u_i \otimes v_i & \longmapsto &
\sum_{i=1}^n \langle v, u_i
\rangle v_i
\end{array}
\]
is equal to $\id : V \to V$. It takes a little linear-algebra exercise
to see that this implies that the symmetric bilinear form is
nondegenerate. This also implies that $\dim V < \infty$, because the
$v_i$'s, $i = 1, \dots, n$, must span $V$. This completes the
construction of the structure of a finite-dimensional Frobenius
algebra on the state space $V$ of a TFT.

Conversely, if we have a finite-dimensional Frobenius algebra $V$, we
can define the structure of a TFT on the vector space $V$ by (1)
cutting a surface down into pairs of pants, cylinders, and caps; (2)
defining the operators corresponding to those basic objects using the
multiplication (or its linear dual), the identity map, and the unit
element $e \in V$ (or the dual of the map $k \to V $, $ 1 \mapsto e$,
as the trace functional), respectively; and (3) using the sewing
axiom. The fact that the composite operator is independent of the way
we cut down the surface follows from Figure~\ref{coherence} and the
associativity of multiplication.
\end{proof}

\begin{ex}[The Toy Model, Dijkgraaf and Witten \cite{dijkgraaf-witten}]
Let $G$ be a finite group.  Assign to any closed one-manifold $X =
\coprod S^1$ the space $V(X)$ of $\nc$-valued functions on the
isomorphism classes of principal $G$-bundles over $X$. The set of
these isomorphism classes may be identified with $\Hom (\pi_1 (X),
G)/G$, the set of conjugacy classes of homomorphisms $\pi_1 (X) \to
G$. To an oriented surface $Y$ with no inputs, assign the operator
$\Psi_Y: \nc \to V(\p Y)$ defined by the equation
\[
\Psi_Y (P) := \sum_{Q : Q|_{\p Y} = P} \frac{1}{\abs{\Aut Q}}
\]
as a function on the set of principal $G$ bundles $P$ over $\p Y$. In
this case, $A = V(S^1)$ becomes the center of the group algebra
$\nc[G]$ with Frobenius' trace $\theta(\sum_g \lambda_g g) =
\lambda_e/ \abs{G}$, where $e \in G$ is the unit element. Note that
the same construction yields an $n$-dimensional TFT for any $n \ge 1$.
\end{ex}

\subsection{(Topological) Conformal Field Theories}

Let us concentrate on the case $n=2$ and, instead of the very small
PROP of oriented smooth surfaces, consider the Segal PROP $\PP$ of
Riemann surfaces (i.e., complex curves) with holomorphic holes from
Example~\ref{Segal:PROP}. Here we will be discussing what is known as
QFTs with zero central charge.

Recall that in Chapter~\ref{cactioperad} we introduced the following
notion.
\begin{definition}
A \emph{Conformal Field Theory} (\emph{CFT}) is an algebra over the
PROP $\PP$.
\end{definition}
We would like to consider variations on this theme.

\begin{definition}
A \emph{Cohomological Field Theory-I} (\emph{CohFT-I}) is an algebra
over the homology PROP $H_*(\PP)$. The Roman numeral one in the name
is to distinguish this theory from a standard \emph{Cohomological
Field Theory} (\emph{CohFT}), which is an algebra over the PROP
$H_*(\overline{\MM})$, where $\overline{\MM}$ is the PROP of moduli
space of stable compact algebraic curves with punctures with respect
to the operation of attaching curves at punctures. This latter theory
is also known as \emph{quantum gravity}. A \emph{Topological Conformal
Field Theory} (\emph{TCFT}) is an algebra over the chain PROP
$C_*(\PP)$ for a suitable version of chains, e.g., singular. This
theory is also referred to as a
\emph{string background}.
\end{definition}

The definitions of a CFT and a TCFT are basically rewordings of those
introduced by Segal \cite{segal:CFT}. The definition of a CohFT is
essentially a rewording of that of Kontsevich and Manin
\cite{konm,manin:frob}, which was perhaps motivated by Witten's paper
\cite{witten:cohft}.

Note that a TFT may be regarded as an algebra over the PROP $H_0
(\PP)$. Often, one gets a TFT by integrating over the moduli space, in
which case one obtains a TFT from the top homology groups of $\PP$ or
$\overline{\MM}$.

\subsection{Examples}

\subsubsection{Gromov-Witten Theory}

\emph{Gromov-Witten theory} is basically what physicists know as a
\emph{sigma-model}. Part of the structure, essentially the genus zero
part, is also known as \emph{quantum cohomology}. Gromov-Witten theory
starts with a complex projective manifold $M$ and nonnegative integers
$g$, $m$, and $n$ and uses a diagram
\[
\begin{CD}
@. [\overline{\MM}_{g,m,n} (M; \beta)]_\virt @.\\ @. @VVV @.\\
\overline{\MM}_{g,m,n} @<f<< \overline{\MM}_{g,m,n} (M; \beta) @>\ev>>
M^{m+n},
\end{CD}
\]
in which $\overline{\MM}_{g,m,n}$ is the moduli space of stable
compact complex curves of genus $g$ with $m+n$ labeled punctures, the
first $m$ of them thought of as the inputs, the last $n$ as the
outputs. Here $\overline{\MM}_{g,m,n} (M; \beta)$ is the moduli space
of stable pairs (a compact complex curve $C$ from
$\overline{\MM}_{g,m,n}$, a holomorphic map $C \to M$ taking the
fundamental class of $C$ to a given homology class $\beta \in H_2 (M;
\nz)$), and $[\overline{\MM}_{g,m,n} (M; \beta)]_\virt$ is the
so-called
\emph{virtual fundamental class}, a certain homology class of
$\overline{\MM}_{g,m,n} (M; \beta)$. This is morally the fundamental
class of the ideal moduli space associated to the corresponding
deformation problem. In the spirit of Chapter~\ref{cactioperad}, the
above diagram defines a correspondence, which allows one to define a
pairing
\[
\Omega^* (M)^{\otimes m+n} \otimes \Omega^* (\overline{\MM}_{g,m,n})
\to \nc[-d_\virt],
\]
given by
\[
\int_{[\overline{\MM}_{g,m,n} (M; \beta)]_\virt} ev^* (\omega_1 \otimes
\dots \otimes \omega_{m+n}) \wedge f^*(\phi),
\]
where $\Omega^*$ refers to the complex of differential forms, $[\;]$
denotes the degree shift, and $d_\virt$ is the dimension of the
homology class $[\overline{\MM}_{g,m,n} (M; \beta)]_\virt$.  Passing
to cohomology and using Poincar\'e duality on $M$, one gets a CohFT
with a degree shift, see
\cite{behrend,li-tian,manin:frob}.

\subsubsection{Gromov-Witten Theory: a richer version}

The following is just a variation on the theme of the previous
example, closer to string topology and what physicists originally
thought of as sigma-model. Instead of the holomorphic mapping moduli
space $\overline{\MM}_{g,m,n} (M; \beta)$, use an infinite dimensional
version $\PP_{g,m,n} (M; \beta)$ of it, the moduli space of
holomorphic maps from nonsingular complex compact curves with
holomorphic holes. Here the holes are not considered removed from the
complex curve, so that the fundamental class of the curve is defined
and mapped to a homology class $\beta$ of $M$. As in the previous
example, we get a correspondence
\[
\PP_{g,m,n} \leftarrow \PP_{g,m,n} (M; \beta) \rightarrow (LM)^{m+n},
\]
which should create something resembling a TCFT via passing to
semi-infinite differential forms on the free loop space $LM = \Map
(S^1, M)$:
\[
(\Omega^{\infty/2 + *} (LM))^{\otimes {m+n}} \otimes \Omega^*
(\PP_{g,m,n}) \to \nc[-d]
\]
for some degree $d$. Much of this has not yet been made rigorous. In
particular the theory of semi-infinite differential forms on the loop
space has never yet been made precise. However this perspective
supplies the motivation for much work in Gromov-Witten theory. For
instance, Lalonde \cite{lalonde} has described field theoretic
properties along these lines with Floer homology replacing
semi-infinite forms.

\subsubsection{String topology}

String topology generalizes to a TCFT in the following sense. Instead
of being based on a moduli space of complex curves, it is based on a
combinatorial counterpart of it, the (partial) PROP of reduced metric
chord diagrams $\sr \cc \cf_{p,q}(g)$, see next Section and also
\cite{cohen-godin,chataur}. Here \emph{partial} refers to the fact that
the PROP composition is only partially defined, namely, when the
circumferences of the outputs and the inverse circumferences of the
inputs, which one thinks of as colors, match. The correspondence
realizing the string-topology TCFT at the level of correspondences is
given by a diagram
\[
\sr \cc \cf_{m,n}(g) \times (LM)^m \hookleftarrow \sr
\cc \cf_{m,n}(g) M \to (LM)^n,
\]
where the space in the middle is the space of continuous (or
continuous, piecewise smooth, depending on the version of a free loop
space $LM$ considered) maps from reduced metric chord diagrams of
genus $g$ with $m$ inputs and $n$ outputs to the target manifold
$M$. Passing to homology and using the Thom collapse map, see Chapter
1, we get a string-topology version of a CohFT given by an algebra
over the homology PROP $H_*(\sr \cc \cf_{m,n}(g))$ defined by a
diagram
\[
H_*(\sr \cc \cf_{m,n}(g)) \otimes H_*(LM)^{\otimes m} \la
 H_{*+\chi d} (LM)^{\otimes n},
\]
where $\chi$ is the Euler characteristic of the chord diagrams of that
type and $d = \dim M$, see also Equation~\eqref{mu}.

\subsection{Motivic TCFTs}

The above examples suggest that one needs a more general definition of
a TCFT than the one given before: the main issue being the need to use
correspondences at the space level or degree shifts at the
(co)homology level. Namely a \emph{motivic TCFT} is a tensor functor
from a category of holomorphic two-dimensional cobordisms (or a
version thereof, such as the reduced metric chord diagrams) enriched
over a motivic category of manifolds with values in a motivic category
of manifolds. The tensor structure on the category of cobordisms is
given by the disjoint union and on the category of manifolds by the
Cartesian product. A motivic category of manifolds may be an
appropriate category of correspondences, as treated earlier in
Section~\ref{cacti:action}. Motivic ideas in TCFT come from Kontsevich
and Manin \cite{konm}, who used a motivic axiom to axiomatize a
Cohomological Field Theory in the context of Gromov-Witten
theory. Defining and studying motivic TCFTs is one of the more
interesting, open projects in the study of topological field theories.

\section{Generalized string topology operations}

Recall from Chapter 1 that the loop homology product is defined by
considering the diagram
$$ LM \xleftarrow{\gamma} \met \xr{e} LM
\times LM,$$
where $e : \met \hk LM \times LM$ is a codimension $d$ embedding.  We
used the Thom collapse map to define an umkehr map in a multiplicative
generalized homology theory that supports an orientation of $M$,
$$e_!: h_*(LM \times LM) \to h_{*-d}(\met).$$ The loop product was
defined to be the composition,
$$ \mu : h_*(LM) \otimes h_*(LM) \to
h_*(LM \times LM) \xr{e_!} h_{*-d}(\met) \xr{\gamma_*} h_{*-d}(LM).
$$

\begin{figure}[ht]  \centering  \includegraphics[height=3cm]{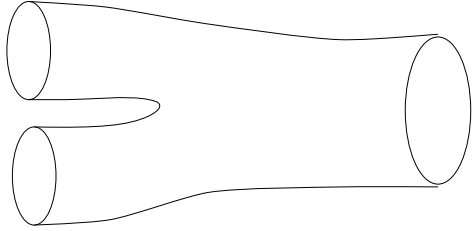}
\caption{The ``pair of pants" surface $P$}   \label{pants}\end{figure}

\med
One can think of this structure in the following way.  Consider the ``pair of pants" surface $P$, viewed as a cobordism from two circles to one circle (see Figure \ref{pants}).  We think of the two circles bounding the ``legs" of the pants as incoming, and the circle bounding the ``waist" as outgoing. 

Consider the smooth mapping space, $Map(P, M)$.  Then there are restriction
maps to the incoming and outgoing boundary circles,
$$
\rho_{in} : Map(P, M) \to LM \times LM, \quad \rho_{out} : Map(P,M) \to LM.
$$
Notice that the figure 8 is  homotopy equivalent to the surface $P$, with respect to which 
the restriction map $\rho_{in} : Map(P,M) \to LM \times LM$ is homotopic to the embedding
$e : Map (8,M) \to LM \times LM$.  Also, restriction to the outgoing boundary, $\rho_{out} : Map(P,M) \to LM$
is homotopic to $\gamma : Map (8,M) \to LM$.  So the Chas-Sullivan product can be thought of as a composition,
$$
\mu_* :  H_p(LM) \otimes H_q(LM) \xr{(\rho_{in})_!} H_{p+q-d}(Map(P,M)) \xr{(\rho_{out})_*} H_{p+q-d}(LM).
$$
The role of the figure 8 can therefore be viewed as just a technical one, that allows us to define the umkehr map
$e_! = (\rho_{in})_!$.

More generally, consider a surface $F$, viewed as a cobordism from $p$-circles to $q$-circles.  See Figure \ref{surface} below.  Again, we think of the $p$ boundary circles as incoming, and the remaining $q$ boundary circles as outgoing. 

\begin{figure}[ht]  \centering  \includegraphics[height=4cm]{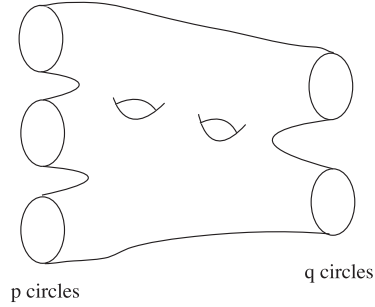}
\caption{The  surface $F$}   \label{surface}\end{figure}

We can consider the mapping space, $Map (F, M)$, and the resulting restriction maps,
\begin{equation}\label{restrict1}
(LM)^q \xleftarrow{\rho_{out}} Map (F, M) \xr{\rho_{in}} (LM)^p.
\end{equation}

In   \cite{cohen-godin} Cohen and Godin showed how to   construct an   umkehr map 
$$
(\rho_{in})_! : h_*((LM)^p) \to h_{* +\chi (F)\cdot d}(Map (F, M))
$$ where $\chi (F)$ is the Euler characteristic of the surface $F$ and as above, $d = dim (M)$.    This then allows
the definition of a string topology operation
\begin{equation}\label{operation}
\mu_F : h_*((LM)^p) \xr{(\rho_{in})_!} h_{* +\chi (F)\cdot d}(Map (F, M)) \xr{(\rho_{in})_*} h_{*+\chi (F)\cdot d}((LM)^q).
\end{equation}

\med
To construct this umkehr map, Cohen and Godin used the Chas-Sullivan
idea of representing the pair of pants surface $P$ by a figure 8, and
realized the surface $F$ by a ``fat graph" (or ribbon graph).  Fat
graphs have been used to represent surfaces for many years, and to
great success.  See for example the following important works:
\cite{harer}, \cite{strebel}, \cite{penner:CMP}, \cite{kon:witten}.

We recall the definition.
\begin{definition}   A fat graph is a finite graph with the following properties:
\begin{enumerate} 
\item Each vertex is at least trivalent 
\item Each vertex comes equipped with a cyclic order of the half edges emanating from it.
\end{enumerate}
\end{definition}

\begin{figure}[ht]  \centering  \includegraphics[height=8cm]{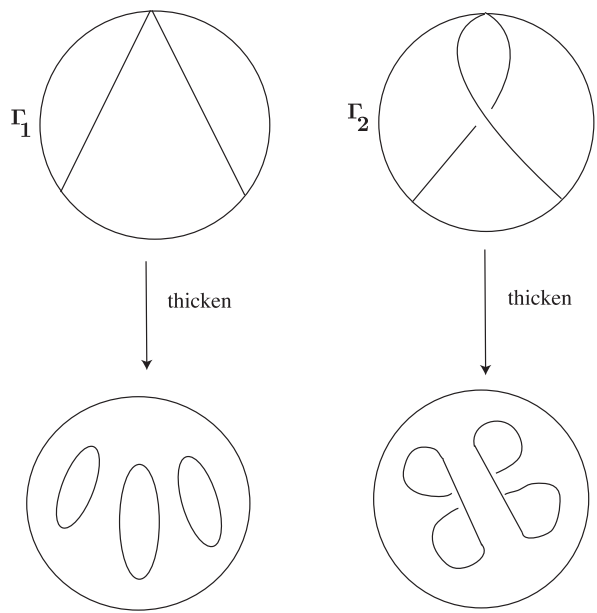}
 \caption{Thickenings of two fat graphs }   \label{thick}\end{figure}

We observe that the cyclic order of the half edges is quite important in this structure.
It allows for the graph to be ``thickened" to a surface with boundary.  This thickening   can be thought of as assigning a ``width" to the ink used in drawing a fat graph.  Thus one is actually drawing a two dimensional space, and it is not hard to see that it is homeomorphic to a smooth surface.  Consider the following two examples (Figure \ref{thick}) of fat graphs which consist of the same underlying graph, but have different cyclic orderings at the top vertex.

  The orderings of the edges are induced by the counterclockwise
orientation of the plane.  Notice that $\G_1$ thickens to a surface of
genus zero with four boundary components.  $\G_2$ thickens to a
surface of genus 1 with two boundary components.  Of course these
surfaces are homotopy equivalent, since they are each homotopy
equivalent to the same underlying graph.  But their diffeomorphism
types are different, and that is encoded by the cyclic ordering of the
vertices.

\begin{figure}[ht]  \centering  \includegraphics[height=4cm]{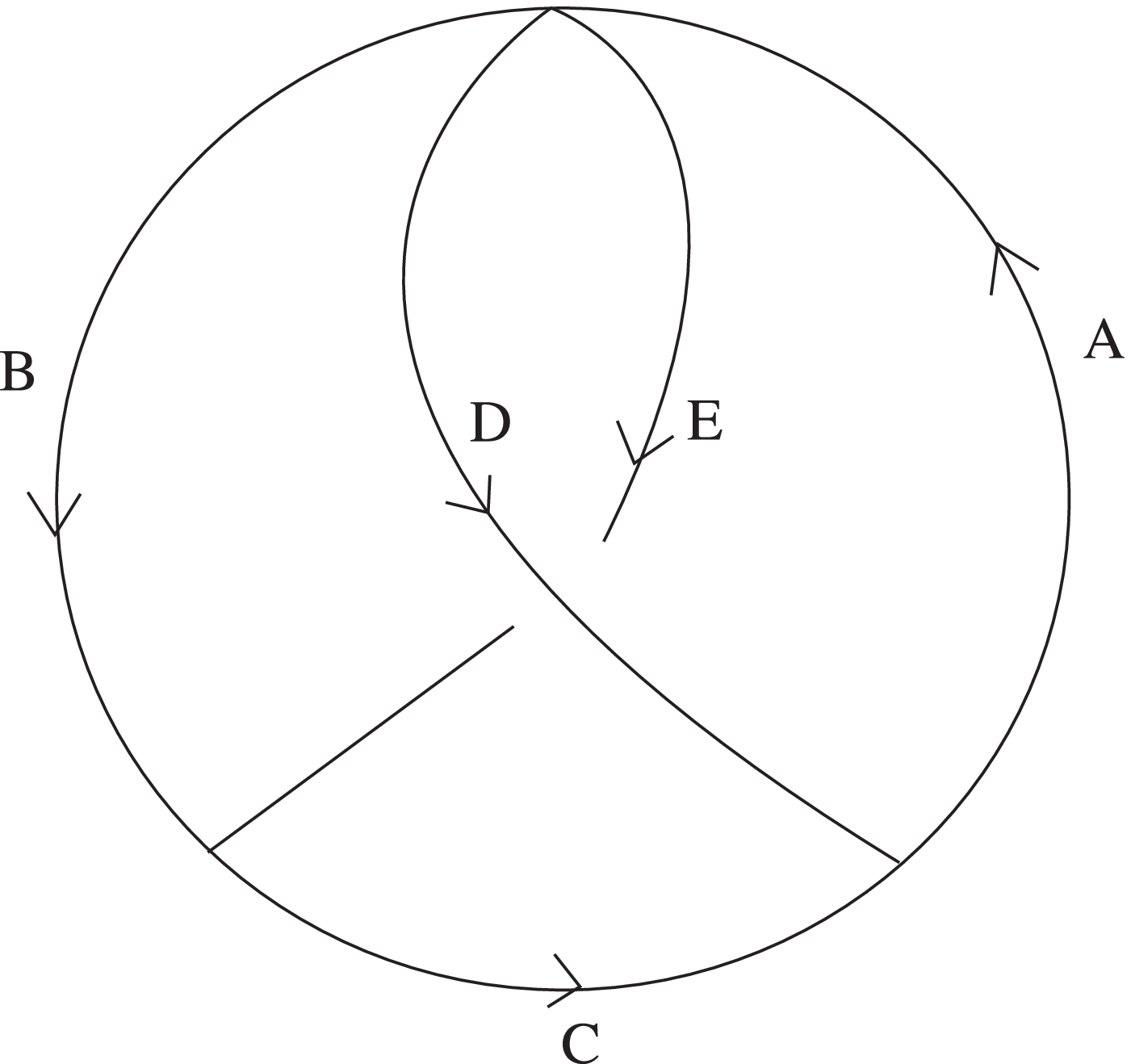}
 \caption{The  fat graph  $\G_2$}   \label{G2}\end{figure}

These examples make it clear that we need to study the combinatorics of fat graphs more carefully. 
For this purpose, for a fat graph $\G$, let $E(\G)$ be the set of edges, and let $\tilde E(\G)$ be the set of oriented edges.  $\tilde E(\G)$ is a $2$-fold cover of $E(G)$.  It has an involution $E \to \bar E$
which represents changing the orientation.  The cyclic orderings at the vertices determines
a partition of $\tilde E(\G)$ in the following way.  Consider  the  example illustrated in Figure \ref{G2}. 

 As above, the cyclic orderings at the vertices are determined by the counterclockwise orientation of the plane. 
 To obtain the partition,  notice that an oriented edge has   well defined source and target vertices. Start with an oriented edge, and follow it to its target vertex.  The next edge in the partition is the next oriented edge in the cyclic ordering at that vertex.  Continue in this way until one is back at the original oriented edge.  This will be the first cycle in the partition.  Then continue with this process
 until one exhausts all the oriented edges.  The resulting  cycles in the partition will be called  \sl ``boundary cycles" \rm as they reflect the boundary circles of the thickened surface. In the  case of $\G_2$ illustrated in Figure \ref{G2}, the partition into boundary cycles is  given by:
  
  $$
  \text{Boundary cycles of $\G_2$:} \quad (A,B,C) \, (\bar A, \bar D, E, \bar B, D, \bar C, \bar E ).
  $$
  
 So one can compute combinatorially the number of boundary components in the thickened surface of a fat graph.  Furthermore the graph and the surface have the same homotopy type, so one can compute the Euler characteristic of the surface directly from the graph.  Then using the formula
 $\chi (F) = 2 - 2g -n$, where $n$ is the number of boundary components, we can solve for the genus directly in terms of the graph.  
 
 A \sl metric fat graph \rm is a connected fat graph  $\G$ endowed with a metric so that the open edges are isometrically equivalent to open intervals in the real line.  Also, if $x, y \in \G$, then each minimal path from $x$ to $y$ in $\G$ is isometrically equivalent to a closed interval in $\br$ of length $d = d(x,y)$.  The set of metric fat graphs can be given a natural topology.  
 The main theorem about these spaces  is the following (see \cite{penner:CMP}, \cite{strebel}).
 
 \begin{thm}\label{penner:CMP}  For $g \geq 2$, the space of metric fat graphs $Fat_{g,n}$ of genus $g$ and $n$ boundary cycles is homotopy equivalent to the moduli space $\cm_{g}^n$ of closed Riemann surfaces of genus $g$ with $n$ marked points.
 \end{thm} 
 
 Notice that  on  a metric fat graph $\G$, the boundary cycles  nearly have well defined parametrizations.      For example, the boundary cycle $(A, B, C)$ of the graph $\G_2$ can be represented by a map $S^1 \to \G_2$ where the circle is of circumference equal to the sum of the lengths of sides $A$, $B$, and $C$.  The ambiguity of the parameterization is the choice of where to send the basepoint $1 \in S^1$.  In her thesis \cite{godin}, Godin described the notion of a ``marked" fat graph, and proved the following analogue of theorem \ref{penner:CMP}.
 
\begin{thm}\label{marked}   Let $Fat^*_{g,n}$ be the space of marked metric fat graphs of genus $g$ and $n$ boundary components.  Then there is a homotopy equivalence
 $$
 Fat^*_{g,n} \simeq \cm_{g,n}
 $$
 where $\cm_{g,n}$ is the moduli space of Riemann surfaces of genus $g$ having $n$ parameterized boundary components.
 \end{thm}
 
\begin{figure}[ht]  \centering  \includegraphics[height=4cm]{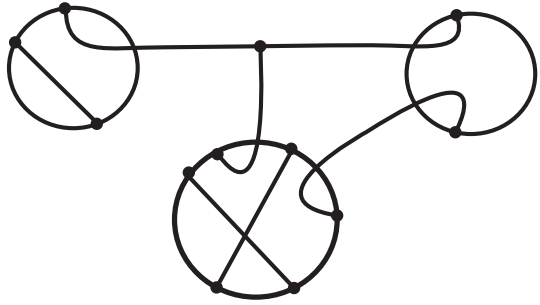}
 \caption{A Sullivan chord diagram}   \label{example}\end{figure}

 In \cite{cohen-godin} the umkehr map $\rho_{in} : h_*((LM)^p) \to
 h_{*+\chi (F)\cdot d}(Map (F, M))$ was constructed as follows.  Let
 $\G$ be a marked fat graph representing a surface $F$.  We will need
 to assume that $\G$ is a special kind of fat graph, which we refer to
 as a ``Sullivan chord diagram".  See Figure \ref{example}.

 \med
 \begin{definition}
A ``Sullivan chord diagram" of type $(g; p,q)$ is a marked fat graph
representing a surface of genus $g$ with $p+q$ boundary components,
that consists of a disjoint union of $p$ disjoint closed circles
together with the disjoint union of connected trees whose endpoints
lie on the circles.  The cyclic orderings of the edges at the vertices
must be such that each of the $p$ disjoint circles is a boundary
cycle.  These $p$ circles are referred to as the incoming boundary
cycles, and the other $q$ boundary cycles are referred to as the
outgoing boundary cycles.  An ordering of these boundary cycles is
also part of the data.
\end{definition}

 In a chord diagram, the vertices and edges that lie on one of the
 $p$-distinct circles are called ``circular vertices and edges".  The
 noncircular edges are known as ``ghost edges".  As stated in the
 definition, the union of the ghost edges is a disjoint union of
 trees.  Each of these trees is known as a ``ghost component".
 
\begin{figure}[ht]  \centering  \includegraphics[height=3cm]{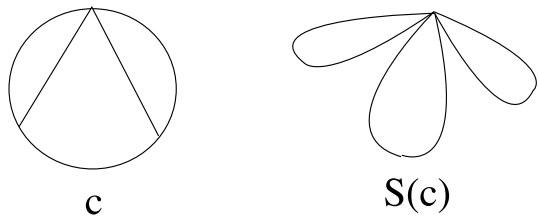}
 \caption{Reducing a chord diagram }   \label{collapse}\end{figure}

      Notice that given a chord diagram $c$,  we can construct 
  a  new fat graph $S(c)$ by collapsing all the ghost edges. That is, each ghost component is collapsed to a point. Notice that  $S(c)$ is a fat graph that represents a surface of the same topological type as that represented by $c$ (the genus and the number of boundary components are the same).  However $S(c)$ is not chord diagram.  We shall refer to $S(c)$ as a ``reduction" of the chord diagram $c$.    See Figure \ref{collapse}.  
  
It is not difficult to see that a Sullivan chord diagram $c$  has the property that if an oriented edge $E$ is contained in an incoming
boundary cycle, then $\bar E$ is contained in an outgoing boundary cycle.  In this setting the reduction $S(c)$ has the property
that $E$ is contained in an incoming boundary cycle  \sl if and only if \rm $\bar E$ is contained in an outgoing boundary cycle.
We define a reduced chord diagram to be a fat graph with this property.

 Let  $\cc \cf_{p,q}(g)$   = space of metric Sullivan chord diagrams of topological type $(g; p, q)$.  This is topologized as a subspace of $\cf at_{g, p+q} $. Let $\sr \cc \cf_{p,q}(g)$ be the corresponding space of reduced metric chord diagrams.
 The following facts are both due to Godin.  The first is contained in \cite{cohen-godin}, and the second is in \cite{godin}.
 
 \begin{prop}  \begin{enumerate}
 \item The space $\cc \cf_{p,q}(g)$ is connected.
 \item The collapse map $\pi :  \cc \cf_{p,q}(g)  \to \sr \cc \cf_{p,q}(g)$  is a homotopy equivalence. 
 \end{enumerate}
 \end{prop}
 
 \med
Let $c \in \cc \cf_{p,q}(g)$, and consider the    mapping space,  $Map (S(c), M)$.  This is the space of continuous maps that are smooth on each edge.  Equivalently, this is the space 
  of continuous maps $f : c \to M$, smooth on each edge, which is constant on each ghost edge. Notice that there is a homotopy equivalence  $  Map (S(c), M) \simeq Map (F_{g, p+q}, M)$, where $F_{g, p+q}$ is a surface of genus $g$ and $p+q$ boundary components.   
 
Markings on $S(c)$ induce   parameterizations of  the incoming and outgoing boundary  cycles   of $c$,   so restriction to these boundary cycles induces a diagram,  $$\begin{CD}(LM)^q @<\rho_{out}<<  Map(S(c) , M) @>\rho_{in} >> (LM)^p,\end{CD}$$
which is homotopic to the diagram \ref{restrict1}.  

Recall that our goal is to define the umkehr map,
  $(\rho_{in})_!$.  This was done in \cite{cohen-godin} as follows. 

  Let $v(c)$ be the collection of circular vertices of   $c$.  Let $\sigma (c)$ be the collection of vertices  of $S(c)$.  The projection   $\pi : c \to S(c)$ determines a surjective set map, $\pi_*: v(c) \to \sigma (c)$. This in turn induces a diagonal map $\Delta_c : M^{\sigma (c)} \to M^{v(c)}$.  Let $c_1, \cdots , c_p$ be the incoming circles in $c$.  The markings define parameterizations and therefore identify $(LM)^p$ with $Map (\coprod_{i=1}^p c_i, M)$. We then have a pullback diagram
  \begin{equation}  \begin{CD}  Map_*(c, M)   @>\rho_{in} >\hk > (LM)^p \\@Ve_c VV   @VVe_c V \\M^{\sigma (c)} @>\hk >\Delta_c >  M^{v (c)}\end{CD}\end{equation}
  Here $e_c$ refers to the map that evaluates at the relevant vertices.

  The  codimension of $ \Delta_c: M^{\sigma (c)} \hk  M^{v (c)}$ is  $(v(c) - \sigma (c))\cdot d$.  But    a straighforward exercise verifies that $(v(c) - \sigma (c))$  
  is equal to minus the Euler characteristic,  $v(c) - \sigma (c) = -\chi (c)  = -\chi (F_{g,p+q})$.  So the codimension of $\Delta_c$ is    $-\chi (c)\cdot d$.

\med
This pullback diagram allows us to construct the Thom collapse map,
  $$
  \tau_c : (LM)^p \la Map_*(S(c), M)^{\eta_c},
  $$
  where $\eta_c$ is the normal bundle.  So in homology we get an umkehr map,
  $$
  \begin{CD}
  (\rho_{in})_! :  H_*((LM)^p) @>(\tau_c)_*>> H_*(Map_*(S(c), M)^{\eta_c}) \\
   @>\cap u >\cong > H_{*+\chi(c)d}(Map_*(S(c), M)).
  \end{CD}
  $$
  This in turn allows us to define the string topology operation,
 \begin{align} \mu_c = (\rho_{out})_* \circ (\rho_{in})_! :  H_*((LM)^p) &\to H_{*+\chi(c)d}(Map_*(S(c), M) ) \notag \\
 & \to H_{*+\chi (c)d}((LM)^q).\end{align}
 
 \med
 
 The following was proved in \cite{cohen-godin}.

 \med
 \begin{thm}\label{stringop}  The homology operation $\mu_c$ only depends on the topology of the surface
 $F_{g, p+q}$. These operations can be defined for any generalized homology $h_*$ that supports an orientation of $M$. They  respect gluing and define a ``positive boundary" topological field theory.
 That is,  there is an operation for every surface with $p$ incoming and $q$ outgoing boundary components
 so long as $q>0$.   Equivalently, this is a Frobenius algebra whose coalgebra structure does not have a co-unit.
 \end{thm}

 \med
 The idea behind the proof of this theorem was to show that one can construct an umkehr map if one allows
 the chord diagram $c$ to vary  in  a continuous family.     The connectedness of the space of chord diagrams, $\cc \cf_{p,q}(g)$
 then estabishes that the operation $\mu_c$ is independent of $c \in \cc \cf_{p,q}(g)$.  The fact that these operations respect gluing, and thereby define a field theory, uses the naturality of the Thom collapse maps.  We refer the reader to \cite{cohen-godin}
 for details. 
 
 We remark that in \cite{cohen-godin} itone-manif was observed that by allowing $c$ to vary over $\cc \cf_{p,q}(g)$, one can actually construct operations,
\begin{equation}
\label{mu}
 \mu : H_*(\cc \cf_{p,q}(g)) \otimes H_*(LM)^{\otimes p} \la H_*(LM)^{\otimes q}.
\end{equation}
 This was verified using Jakob's bordism approach to homology theory
 by Chataur in \cite{chataur}, using the language of partial PROPs.
 Most of the basic results about string topology described in these
 notes were also verified using this theory.  We refer the reader to
 \cite{chataur} for details regarding this very appealing, geometric
 approach to string topology.
 
 \med
 \section{Open-closed string topology}
 
 In this section we describe work of Sullivan \cite{sullivan:oc} on open-closed string topology.  We also discuss generalizations
 and expansions of this theory due to Ramirez \cite{ramirez}.  Similar constructions and results were also obtained by Harrelson \cite{harrelson}. 
 %%%%%%%%%%%%%%%%%%
 In this setting our background manifold comes equipped with a collection of submanifolds,
 $$
 \scrb = \{ D_i \subset M \}.
 $$
 Such a collection is referred to as a set of ``D-branes", which in string theory supplies
 boundary conditions  for open strings.  In string topology, this is reflected by
 considering the path spaces
 $$
 \cp_M(D_i, D_j) = \{ \gamma : [0,1] \to M, \,: \,  \gamma (0) \in D_i, \quad \gamma (1) \in D_j\}.
 $$
Following Segal's viewpoint,  \cite{segal}, in a theory with $D$-branes, one associates to a connected, oriented compact one-manifold $S$  whose boundary components are labeled by $D$-branes, a vector space $V_S$.  In the case of string topology,
if $S$ is topologically a  circle, $S^1$,  the vector space $V_{S^1} = h_*(LM)$.  If $S$ is an interval with boundary points labeled
by $D_i$ and $D_j$, then $V_S = h_*(\cp_M(D_i, D_j))$.  As is usual in field theories, to a disjoint
union of such compact one-manifolds, one associates the tensor product of the vector spaces assigned to each connected component. 

Now to an appropriate cobordism, one needs to associate an operator between the vector spaces associated to the incoming and outgoing parts of the boundary.  In the presence of $D$-branes these cobordisms are cobordisms of manifolds with boundary.  More precisely, in a theory with $D$-branes, the boundary of a cobordism $F$ is partitioned into three parts:

\begin{enumerate}\label{cobordism}
\item incoming circles and intervals, written $\p_{in}(F)$,
\item outgoing circles and intervals, written $\p_{out}(F)$,
\item the ``free part" of the boundary, written $\p_{f}(F)$, each component of which is labeled by a $D$-brane.  Furthermore $\p_f(F)$ is a cobordism from the boundary of the incoming one-manifold to the boundary of  the outgoing one-manifold.   This cobordism respects the labeling.
\end{enumerate}  

\begin{figure}[ht]  \centering  \includegraphics[height=5cm]{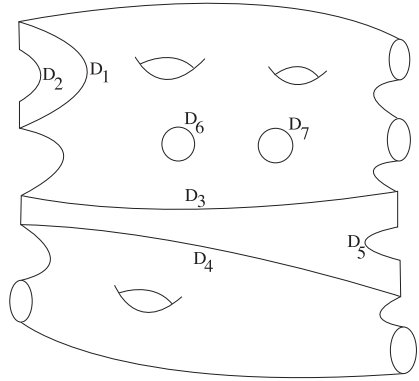}
 \caption{open-closed cobordism }   \label{openclosed}\end{figure}

We will call such a cobordism an ``open-closed cobordism" (see Figure \ref{openclosed}).  

 \med
 We remark that the topology of the category of open and closed strings has been evaluated  by Baas, Cohen, and Ramirez in \cite{baascohenramirez}.  This is actually a symmetric monoidal 2-category, where the objects are
 compact one manifolds $S$, whose boundary components are labeled by elements in a set of $D$-branes, $\scrb$.  The morphisms are open-closed cobordisms, and the 2-morphisms are diffeomorphisms of these cobordisms. Let $\cs^{oc}_{\scrb}$ denote this 2-category, and $|\cs^{oc}_{\scrb}|$ its geometric realization.  They  used  Tillmann's  work on the category of closed strings \cite{tillmann}, as well as the striking
 theorem of Madsen and Weiss \cite{madsen-weiss} proving Mumford's conjecture about the stable cohomology of mapping class groups,  to prove the following.
 
 \med
 \begin{thm}There is a homotopy equivalence of infinite loop spaces,    $$    \Omega|\cs^{oc}_{\scrb}| \simeq  \Omega^\infty \left( (\bc \bp^\infty)^{-L}\right)   \times \prod_{D \in \scrb} Q(\bc \bp^\infty_+)    $$  where, as usual, $X_+$ denotes $X$    with a disjoint basepoint and $Q(Y) =    \varinjlim_{k}\Omega^k\Sigma^k(Y)$.   Here $\Omega^\infty \left( (\bc \bp^\infty)^{-L}\right) $ is the zero space of the   Thom spectrum
 of the virtual bundle $-L$, where $L \to \bc \bp^\infty$ is the canonical line bundle. 
 \end{thm}
 
 \med

 In a theory with $D$-branes, associated to   an open-closed cobordism $F$ there is an operator,
 $$\Phi_F : V_{\p_{in}(F)} \la V_{\p_{out}(F)}.$$  Of course such a theory must respect gluing of open-closed cobordisms.  
 
 Such a theory with $D$-branes  has been put into the  categorical language of PROPs by Ramirez \cite{ramirez} extending notions of Segal and Moore \cite{segal}.  He called such a field theory a $\scrb$-topological quantum field theory.
 
 In the setting of string topology, operators $\Phi_F$ were defined by Sullivan \cite{sullivan:oc} using transversal intersections of chains.  They were defined via Thom-collapse maps by Ramirez in \cite{ramirez}. We will illustrate his definitions by the following examples.  
    
\begin{figure}[ht]  \centering  \includegraphics[height=4cm]{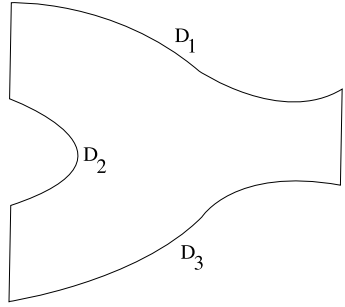}
 \caption{The cobordism $C_1$ }   \label{cob1}\end{figure}

    Consider the  genus zero open-closed cobordism, $C_1$, with free boundary components labeled by $D$-branes,
    $D_1$, $D_2$, and $D_3$ as indicated in Figure \ref{cob1}.
    
 This cobordism defines an operation
 $$
 \Phi_{C_1} : H_*( \cp_M(D_1, D_2)) \otimes H_*( \cp_M(D_2, D_3))  \la H_*( \cp_M(D_1, D_3))
 $$
 in the following way.  
 Let $Map_{\scrb} ( C_1 , M)$ denote the  space of smooth maps  $\gamma : C_1 \to M$,  where the restriction of $\gamma$ to a boundary interval  labeled by the $D$-brane  $D_i \in \scrb$,    takes values in $D_i \subset M$.  Notice we have a diagram of restriction maps,
 \begin{equation}
 \cp_M(D_1, D_3) \xleftarrow{\rho_{out}} Map_{\scrb} ( C_1 , M) \xr{\rho_{in}}  \cp_M(D_1, D_2) \times  \cp_M(D_2, D_3)
 \end{equation}
 
 As in the case of closed string operations, the main idea is to construct an umkehr map
 $$
( \rho_{in})_! :  H_*( \cp_M(D_1, D_2)) \otimes H_*( \cp_M(D_2, D_3)) \to H_*(Map_{\scrb} ( C_1 , M) )
$$
and then define the operation $\Phi_{C_1}$ to be the composition,
\begin{align}
\Phi_{C_1} = (\rho_{out})_* \circ ( \rho_{in})_! :  H_*( \cp_M(D_1, D_2)) &\otimes H_*( \cp_M(D_2, D_3)) \to H_*(Map_{\scrb} ( C_1 , M) ) \notag \\
& \to H_*( \cp_M(D_1, D_3)). \notag
\end{align}

The umkehr map $(\rho_{in})_!$ was defined by replacing the mapping space $ Map_{\scrb} ( C_1 , M)  $
by the path space
$$
\cp_M(D_1, D_2, D_3) = \{\alpha : [0,1] \to M\, : \, \alpha(0) \in D_1, \, \alpha (\frac{1}{2}) \in D_2, \, \alpha (1) \in D_3 \}.
$$
We notice that there is a restriction map $r :  Map_{\scrb} ( C_1 , M)  \to \cp_M(D_1, D_2, D_3)$
which is a homotopy equivalence.  Furthermore we observe that there is a pullback diagram of fibrations,
$$
\begin{CD}
\cp_M(D_1, D_2, D_3)   @>>>    \cp_M(D_1, D_2) \times  \cp_M(D_2, D_3) \\
@Vev_{\frac{1}{2}} VV      @VVev_1 \times ev_0 V \\
D_2 @>>\Delta >   D_2 \times D_2,
\end{CD}
$$ where the vertical fibrations are evaluation maps at the times given by the subscripts.  As argued previously, the existence of this pullback square allows for the definition of a Thom collapse map,
$$
\tau :  \cp_M(D_1, D_2) \times  \cp_M(D_2, D_3) \to (\cp_M(D_1, D_2, D_3))^{\nu_\Delta} = (\cp_M(D_1, D_2, D_3))^{TD_2},
$$
which, in homology defines the umkehr map we are looking for.    

\begin{figure}[ht]  \centering  \includegraphics[height=2cm]{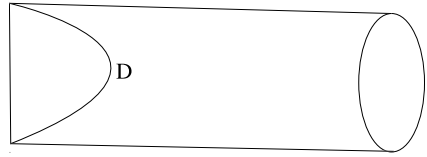}
 \caption{The cobordism $C_2$ }   \label{cob2}\end{figure}

We consider one more example.  Consider the open-closed  cobordism $C_2$ between an interval, whose boundary is labeled by a $D$-brane $D$, and a circle.  This   cobordism is  pictured in Figure \ref{cob2}.
    
 As in the previous example, we consider the mapping space, $Map_{\scrb} ( C_2 , M)$ consisting of maps
 from the surface to $M$, which map  the free part of the boundary to $D \subset M$.  Then there is a diagram of restriction maps,
  \begin{equation}
LM \xleftarrow{\rho_{out}} Map_{\scrb} ( C_2 , M) \xr{\rho_{in}}  \cp_M(D , D ).
 \end{equation}
 As before, the operation $\Phi_{C_2}$ will be defined to be the composition,
$$
\Phi_{C_2} = (\rho_{out})_* \circ ( \rho_{in})_! :  H_*( \cp_M(D , D ))    \to H_*(Map_{\scrb} ( C_2 , M) )  
  \to H_*( LM). 
$$
The umkehr map $(\rho_{in})_!$ is defined as above, except that we now replace $Map_{\scrb} ( C_2 , M) $ by the path space
$$
L_D(M) = \{\alpha \in LM \, : \, \alpha (0) \in D \}.
$$
Again, there is a restriction map $r : Map_{\scrb} ( C_2 , M) \to L_D(M)$ that is a homotopy equivalence.  Furthermore, there is a pullback square,
$$
\begin{CD}
L_D(M) @>>>  \cp_M(D,D) \\
@V ev_0 VV   @VV ev_0 \times ev_1 V \\
D @>>\Delta > D \times D
\end{CD}
$$
from which we construct a Thom collapse map, $\tau :  \cp_M(D,D) \to (L_D(M))^{TD}$ and  the induced 
umkehr map, $(\rho_{in})_! : H_*( \cp_M(D,D))   \to H_{*-dim (D)} (L_D(M))$.  

\med
In order to construct operations $\Phi_\Sigma$ for an arbitrary open-closed cobordism $\Sigma$ as in Figure \ref{openclosed},  Ramirez replaced mapping spaces $Map_{\scrb}(\Sigma, M)$ by homotopy equivalent
spaces, $Map_{\scrb}(\G, M)$, where $\G$ is an appropriate fat graph with boundary labels.  Ramirez
defined the appropriate concept of these graphs, which he called ``$\scrb$-fat graphs", studied the moduli space of such graphs, used them to construct     operations $\Phi_\Sigma$, and proved the following theorem \cite{ramirez}.

 \med
 \begin{thm} Given a set of $D$-branes $\scrb$ in a manifold $M$ and a generalized homology theory $h_*$ that supports orientations of $M$ and all the submanifolds of $\scrb$, then the open-closed string topology operations define a positive boundary $\scrb$- topological quantum field theory.
 \end{thm}

%%%%%% Version of 3/21/05

  \chapter{A Morse theoretic viewpoint}

 In this chapter we describe a Morse theoretic approach to string topology
 and its relationship to the Floer theory of the cotangent bundle.  We will survey
 work contained in \cite{cohennorbury}, \cite{cohenmontreal}, \cite{cohen2}, \cite{salamonweber}, \cite{abschwarz}, \cite{schwarzmontreal}.

 \section[Cylindrical gradient graph flows]{Cylindrical gradient graph flows and string topology}
 
 As already pointed out, string topology has  many of the same formal
 properties and structure as more geometric theories, such as  Floer theory, and Gromov-Witten theory.  In these theories, invariants are obtained by counting, in an appropriate
 sense, maps of surfaces to the manifold that satisfy  the Cauchy- Riemann equations, or certain pertubations of these equations.   The point of this section is to outline work
 in \cite{cohen2} that shows that the string topology invariants also can be computed
 by counting maps of certain topological surfaces to a manifold $M$, that satisfy
 certain differential equations.  The differential equations in this case are gradient flow equations of Morse functions on the loop space.  In the next section we will use this point of view to describe recent work exploring relationships between string topology and Gromov-Witten theory of the cotangent bundle.
 
 The surfaces in this study come as thickenings of fat graphs, or more specifically, marked Sullivan chord diagrams as described in chapter 3.  Recall from section 2 of that chapter that a reduced Sullivan chord diagram $\G$ of type $(g; p,q)$ is a fat graph whose $p+q$ boundary cycles are subdivided into $p$ incoming, and $q$ outgoing cycles, and because of their markings,  come equipped with parameterizations, which
 we designate by $\alpha^-$ for the incoming, and $\alpha^+$ for the outgoing boundaries:
 
  $$
  \alpha^- : \coprod_p S^1 \la \G, \quad \alpha^+ : \coprod_q S^1 \la \G.
  $$
  
  If we endow $\G$ with a metric, we can take the circles to have circumference equal to the sum of the lengths of the edges making up the boundary cycle it parameterizes, each component of $\alpha^+$ and $\alpha^-$ is a local isometry.   
  
  Define the surface $\Sigma_\G$ to be the mapping cylinder of these parameterizations,
  \begin{equation}\label{sigmap}
  \Sigma_\G = \left( \coprod_p S^1 \times (-\infty, 0]\right)   \sqcup \left(\coprod_q S^1 \times [0, +\infty)\right)  \bigcup  \G / \sim
  \end{equation}
  where $(t,0) \in S^1 \times (-\infty, 0] \, \sim \alpha^- (t) \in \G$, and $(t,0) \in S^1 \times [0, +\infty)  \, \sim \alpha^+ (t) \in \G$
  
  \med
  
  Notice that the figure 8  is a  fat graph representing a surface of genus $g=0$ and $3$ boundary components.  This graph has two edges, say $A$ and $B$, and has boundary cycles $(A), (B), (\bar A, \bar B)$.  If we let $(A)$ and $(B)$ be the incoming cycles and $(\bar A, \bar B)$ the outgoing cycle, then the figure 8 graph becomes a chord diagram.
Figure \ref{sigmap}  is a picture of the surface $\Sigma_\G$, for $\G$ equal to the figure 8.

\begin{figure}[ht]
  \centering
  \includegraphics[height=5cm]{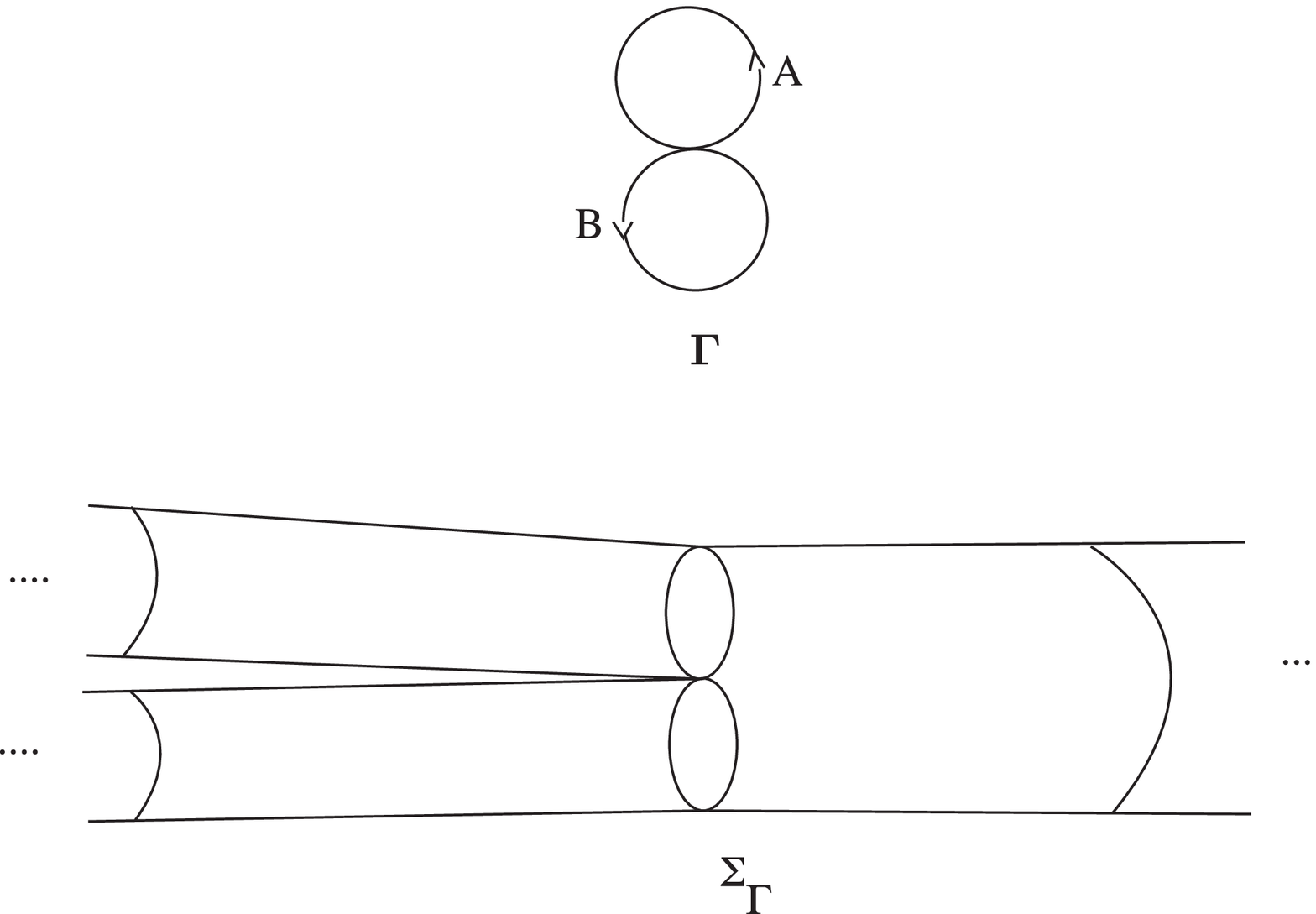}
 \caption{$\Sigma_\G$ }
   \label{sigmag}
\end{figure}

We want to study maps $\gamma : \Sigma_\G \to M$.  Notice that since $\Sigma_G$ is made up of $p+q$ half cylinders, this is equivalent to considering $p$-maps,
$\gamma_i : (-\infty, 0] \times S^1 \to M$, and $q$-maps $\gamma_j : [0, +\infty) \times S^1 \to M$ that satisfy an intersection condition at time $t=0$ determined by the combinatorics of the fat graph $\G$.  Equivalently, these are maps from half lines to the loop space, 
$$
\gamma_i : (-\infty, 0] \to LM,  \, i = 1, \cdots, p, \quad \text{and} \quad  \gamma_j : [0, +\infty)\to LM,  \, j = 1, \cdots , q
$$
that satisfy the appropriate intersection properties at $t=0$.  We will study those maps
$\gamma : \Sigma_\G \to M$ so that the curves $\gamma_i$ and $\gamma_j$ satisfy the gradient flow lines of certain Morse functions on $LM$.  

The constructions we are about to describe are motivated by the use of moduli spaces of ``gradient graph flows"  that have been used to define classical cohomology operations on compact manifolds in \cite{betzcohen}, \cite{fukaya}, and \cite{cohennorbury}.  We also refer the reader to \cite{cohenmontreal} for a general description of this theory, well as its adaptations to string topology.

\med
Assume $M$ is endowed with a Riemannian metric, and $LM$ has the induced $L^2$-metric. 
Let $f : LM \to \br$ be a Morse function on the loop space, which is bounded below, and that satisfies the Palais - Smale condition.  Recall that this condition states that if
$\{a_i\}$ is a sequence of points on which the gradient tends to zero, $\lim_{i\to \infty} \nabla f (a_i)  = 0$, then there is a subsequence that converges to a critical point.

Recall that given such a function $f$, a gradient flow line of $f$ is a curve
$\gamma : \br \to LM$ that satisfies the gradient flow equation:
\begin{equation}
\frac{d \gamma}{dt} + \nabla f (\gamma (t)) = 0.
\end{equation}
If $a \in LM$ is a critical point of $f$, let $W^u(a)$ be the unstable manifold, and $W^s(a)$ be the stable manifold.  Recall that $x \in W^u(a)$ if and only if
there is a gradient flow $\gamma : \br \to LM$ satisfying the initial condition,
$\gamma (0) = x$, and the asymptotic condition, $\lim_{t\to -\infty} \gamma (t) = a$.
The stable manifold $W^s(a)$ is defined similarly except the boundary condition is that
$\lim_{t \to +\infty}\gamma (t) = a$.  

In the case we are considering, ($f : LM \to \br$ Morse, bounded below, and satisfies
the Palais-Smale condition), then $W^u(a)$ is diffeomorphic to a disk of dimension $\lambda (a)$ called the index of $a$, and $W^s(a)$ is infinite dimensional, but has
codimension equal to $\lambda (a)$.  When, in addition, $f$ satisfies the Morse-Smale transversality conditions, i.e the unstable and stable manifolds intersect transversally, one can consider the moduli spaces of gradient flow lines,
$$
\cm (a,b) = W^u(a) \cap W^s(b)/ \br
$$
where $\br$ acts by reparameterizing the flow lines.  This is a manifold of dimension
$\lambda (a) - \lambda (b) -1$.  

In this case there 
  is a $CW$ complex homotopy equivalent to $LM$, which is built out of one cell of dimension $p$ for every critical point $a$ of index $\lambda (a) = p$.   The cellular chain complex is referred to as the Morse complex, $C_*^f(LM)$:
\begin{equation}
\to \cdots \xr{\p_{p+1}} C_p^f(LM) \xr{\p_p} C_{p-1}^f(LM) \xr{\p_{p-1}} \cdots
\end{equation}
Here $C_p^f(LM)$ is the free abelian group generated by the critical points of index $p$.
As usual, if $M$ is oriented, the  boundary homomorphism can be computed by the formula
$$
\p_p (a) = \sum_{\lambda (b) =p-1}n_{a,b}\cdot b
$$
where $n_{a,b} = \#\cm(a,b)$.  This number is counted with sign, which makes sense because  in this case $\cm (a,b)$ is a compact, oriented, zero dimensional manifold.  We will let $C^*_f(LM)$ denote
the dual cochain complex for computing $H^*(LM)$.  

 \med
 We now discuss a plentiful supply of Morse functions on $LM$.   Consider a potential function on $M$, defined to be a smooth map
 $$
 V : \br/\bz \times M \la \br.
 $$ 
 We can then define the classical energy functional
 \begin{align}\label{energy}
 \cs_V :  LM &\la \br \notag \\
 \gamma &\la \int_0^1\left(\frac{1}{2} |\frac{d\gamma}{dt}|^2 - V(t, \gamma (t)) \right) dt.
 \end{align}
 
 For a generic choice of $V$,  $\cs_V$ is a Morse function \cite{weber} satisfying the Palais-Smale condition. 
 Its critical points are those $\gamma  \in LM$ satisfying the ODE
 \begin{equation}\label{critical}
 \nabla_t\frac{d\gamma}{dt} = -\nabla V_t(x)
 \end{equation}
 where $\nabla V_t(x)$ is the gradient of the function $V_t (x) = V(t,x)$, and $ \nabla_t\frac{d\gamma}{dt} $ is the Levi-Civita covariant derivative.

 \med
 We will now describe a metric and Morse theoretic structure on a marked chord
 diagram $\G$.  This will allow us to define the differential equations that we would like
 the maps $\Sigma_\G \to M$ to satisfy.
 
 \med
 \begin{definition}\label{lmmorse}  Given a marked chord diagram $\G$ with $p$-incoming and $q$-outgoing boundary cycles, we define
an $LM$-Morse structure $\sigma$  on $\G$ to be a metric  on $\G$ together with a labeling of each boundary cycle   of $\G$  by a distinct  Morse function $f : LM \to \br$ which is bounded below, and that satisfies the Palais-Smale condition.    
\end{definition}

\med
Notice that we can think of such a labeling of boundary cycles, as a labeling
of the boundary \sl cylinders \rm of the surface $\Sigma_\G$.  Notice also that
we can choose our Morse functions to be energy functions of the sort mentioned above,
in which case the labeling can be taken to be by potential functions, $V : \br/\bz \times M \to \br$. 

This leads to the following definition of the moduli space of cylindrical flows.

\begin{definition}\label{cylgraph}  Let $\G$ be a marked chord diagram as above.  Let $\sigma $ be a $LM$-Morse structure on $\G$.  Suppose 
$\phi : \Sigma_\G \to M$ is a continuous map, smooth on the cylinders.   Let $\phi_i : S^1 \times (-\infty, 0] \to M  \, $ be the  restriction of $\phi$ to the $i^{th}$ incoming cylinder, $ i = 1, \cdots, p$,  and $\phi_j : S^1 \times [0, +\infty) \to M  \,$ be the restriction to the $j^{th}$ outgoing cylinder,  $ j  =  1, \cdots, q$. We consider the $\phi_i$'s and $\phi_j$'s as curves in the loop space, $LM$.  Then the moduli space of cylindrical flows  is defined to be
\begin{align}
 \cm^\sigma_\G(LM) = \{\phi : \Sigma_\G \to M \, : \,   \frac{d\phi_i}{dt} +  &\nabla f_i(\phi_i(t)) = 0 \quad  \text{and} \quad   \frac{d\phi_j}{dt} + \nabla f_j (\phi_j(t)) =0 \notag \\
  &\text{for} \quad  i = 1, \cdots , p \quad \text{and} \quad  j = 1, \cdots, q . \} \notag
\end{align}
\end{definition}

\med
  Let $\phi \in \cm^\sigma_\G(LM) $.  For $i = 1, \cdots , p$, let $\phi_{i,-1} : S^1 \to M$ be the restriction of $\phi_i : S^1 \times (-\infty, 0] \to M$ to $S^1 \times \{-1\}$.  Similarly, for $j = 1, \cdots, q$,
  let $\phi_{j,1} : S^1 \to M$ be the restriction of $\phi_j$ to $S^1 \times \{1\}$.  These restrictions define  the following maps.   
    
  \begin{equation}\label{restrict2}
 (LM)^q \xleftarrow{\rho_{out}}  \cm^\sigma_\G(LM) \xr{\rho_{in}} (LM)^p .
  \end{equation}
  
  \med
  In \cite{cohen2} it is shown that one can define a Thom collapse map,
  $$
 \tau_\G:  (LM)^p \to (\cm^\sigma_\G(LM)^\nu
 $$
 where $\nu$ is a certain vector bundle of dimension $-  \chi (\G)\cdot d$. This can be thought of as a normal bundle in an appropriate sense.    This allows the definition of an 
 umkehr map 
\begin{equation} \label{umkehrmorse}
 (\rho_{in})_! : h_*((LM)^p) \to h_{*+ \chi (\G) \cdot d}(\cm^\sigma_\G(LM))
\end{equation}
 for any homology theory $h_*$ supporting an orientation of $M$.  One can
 then define an operation
\begin{equation}\label{morseop}
 q_\G^{morse} : h_*((LM)^p) \xr{ (\rho_{in})_!}  h_{*+ \chi (\G) \cdot d}(\cm^\sigma_\G(LM)) \xr{(\rho_{out})_*} h_{*+\chi (\G)\cdot d}((LM)^q).
\end{equation}

 Consider the inclusion map, $j : \cm^\sigma_\G(LM)  \hk Map (\Sigma_G, M)$.  In
 With respect to this map, diagram \ref{restrict2} is the restriction of diagram \ref{restrict1} described in the last chapter.  Furthermore, the construction of the Thom collapse map  $\tau_G$ and resulting umkehr map $(\rho_{in})_!$ in \cite{cohen2} is compatible with the corresponding Thom collapse map and umkehr map from \cite{cohen-godin}
 described in the last chapter.    This then yields the following theorem, 
  proved in \cite{cohen2}.
   
   \begin{thm}\label{morsestring} For any marked chord diagram $\G$, the Morse theoretic operation 
   $$q_\G^{morse} : h_*((LM)^p) \la  h_{*+\chi (\G)\cdot d}((LM)^q)$$
   given in (\ref{morseop}) is equal to the string topology operation
   $$
q_\G  : h_*((LM)^p) \la  h_{*+\chi (\G)\cdot d}((LM)^q)
$$
defined  in theorem \ref{stringop}.
\end{thm}  

This Morse theoretic viewpoint of the string topology operations has another, more
geometric due to Ramirez \cite{ramirez}.  
It is a  direct analogue of the perspective on the graph
operations in \cite{betzcohen}.   

As above, let $\G$ be a marked chord diagram, and   
  $\sigma$ be an $LM$-Morse structure on $\G$.   Let $ (f_1, \cdots, f_{p+q})$  be the Morse functions on $LM$ labeling the $p+q$   cylinders of $\Sigma_\G$.  As above, the first $p$ of these cylinders are incoming, and the remaining $q$ are outgoing.

Let $\vec{a} = (a_1, \cdots , a_{p+q})$ be a sequence of loops such that $a_i  \in LM$
is a critical point of $f_i : LM \to \br$.   Then define
 
 \med
$\cm^\sigma_\G (LM, \vec{a}) = \{  
  \phi : \Sigma_\G \to M $  that satisfy the following two conditions:   
  \begin{enumerate} 
  \item $\frac{d\phi_i}{dt} + \nabla f_i (\phi_i(t)) = 0 $     for   $  i = 1, \cdots , p+q$  \\
  \item $ \phi_i \in W^u(a_i)$ for $ i = 1, \cdots , p$, and $\phi_j \in W^s(a_j)$ for $j = p+1, \cdots p+q$.\}
  \end{enumerate}
  
  Ramirez then proved that under sufficient transversality conditions described in \cite{ramirez} then $\cm^\sigma_\G (LM, \vec{a})$ is a smooth manifold of dimension
\begin{equation}\label{dimension}
  dim (\cm^\sigma_\G (LM, \vec{a})) = \sum_{i=1}^p Ind (a_i) - \sum_{j=p+1}^{p+q} Ind (a_j)   + \chi (\G)\cdot d.
 \end{equation}
 
 Moreover, an orientation on $M$ induces an orientation on $\cm^\sigma_\G (LM, \vec{a})$.  Furthermore  compactness issues are addressed, and it is shown
 that if $dim (\cm^\sigma_\G (LM, \vec{a})= 0$ then it is compact.  This leads to the following definition.
 For $f_i$ one of the labeling Morse functions, let $C_*^{f_i}(LM)$ be the 
Morse chain complex for computing $H_*(LM)$, and let $ C^*_{f_i}(LM)$ be the corresponding cochain complex.
Consider the chain

\begin{align}\label{chain}
q_\G^{morse}(LM) =  &\sum_{dim (\cm^\sigma_\G (LM, \vec{a}))= 0}  \#\cm^\sigma_\G (LM, \vec{a}) \cdot [\vec{a}]   \\
 & \in \quad  \bigotimes_{i=1}^p C^*_{f_i}(LM)  \otimes \bigotimes_{j=p+1}^{p+q}C_*^{f_j}(LM) \notag  \end{align}
  
  We remark that the (co)chain complexes $C^*_{f_i}(LM)$ are generated by critical points, so this large tensor product of chain complexes is generated by vectors of critical points $[\vec{a}]$.     It is shown in \cite{ramirez} that this chain is a cycle and if one uses (arbitrary) field coefficients this defines   a class 
\begin{align}\label{geodef}
 q_\G^{morse}(LM)& \in   (H^*(LM))^{\otimes p} \otimes (H_*(LM))^{\otimes q}  \\
 &= Hom ((H_*(LM))^{\otimes p}, (H_*(LM))^{\otimes q}). \notag
 \end{align}

 Ramirez then proved that 
  these operations are the same as those defined by (\ref{morseop}), and hence by theorem \ref{morsestring}  is equal to the string topology operation.  In the case when $\G$ is the figure 8,   this operation is the same as that defined by Abbondandolo and Schwarz \cite{abschwarz} in the Morse homology of the loop space. 
 
 \section{Cylindrical holomorphic curves in $T^*M$.}

This section is somewhat speculative.  It is based on conversations with Y. Eliashberg
and is motivated by the work of Salamon and Weber \cite{salamonweber}.  The goal
of the work described in this section can be divided into two parts.

\begin{enumerate}
\item In the previous section, string topology operations were defined in terms of the  moduli
space  of  cylindrical graph flows, $\cm^\sigma_\Gamma (LM)$. Here $\G$ is a marked  chord diagram, and $\sigma$ is an $LM$-Morse structure on $\Gamma$ (see definitions \ref{lmmorse} and \ref{cylgraph}.)   These consisted of maps $\gamma : \Sigma_\Gamma \to M$, that satisfy appropriate gradient flow equations on the 
cylinders, dictated by the structure $\sigma$.  We would like to replace this moduli
space by a space of maps to the cotangent bundle,  $\phi : \Sigma_\G \to T^*M$ that satisfy appropriate Cauchy-Riemann equations when restricted to the cylinders.  These equations are determined by the structure $\sigma$, and an almost complex structure on $T^*M$.
\item We would like to understand how string topology type invariants defined using these moduli spaces of  ``cylindrical holomorphic curves",  are related to invariants such as the Gromov-Witten invariants, which are defined using moduli spaces of holomorphic  curves from a Riemann surface.
\end{enumerate}

We will just give outlines of the ideas of this program in this section. This program is described in more detail in \cite{cohenmontreal}.

\med
The fact that the cotangent bundle $T^*M$ has an almost complex structure comes
from the existence of its canonical symplectic structure.  This structure is defined as follows:

  Let $p : T^*M \to M$ be the projection map.
Let $x\in M$ and $u \in T^*_xM$.  Consider the composition
$$
\alpha (x,u) : T_{(x,u)}(T^*(M)) \xr{Dp} T_x M \xr{u} \br
$$
where $T_{(x,u)}(T^*(M))$ is the tangent space of $T^*(M)$ at $(x,u)$, and $Dp$ is the derivative of $p$. Notice that $\alpha$ is a one form,
$ 
\alpha \in \Omega^1(T^*(M)),
$ 
and we define 
$$
\omega = d\alpha \in \Omega^2(T^*(M)).
$$
 The form  $\omega$ is a nondegenerate symplectic form on $T^*(M)$.  Now given a Riemannian metric on $M$,
$g : TM \xr{\cong} T^*M$, one gets a corresponding almost complex structure $J_g$ on $T^*(M)$ defined as follows.

\med
The Levi-Civita connection defines a splitting of the tangent bundle of the $T^*(M)$,
$$
T(T^*(M)) \cong p^*(TM) \oplus p^*(T^*(M)).
$$
With respect to this splitting,  $J_g : T(T^*(M)) \to T(T^*(M))$ is defined by the matrix,
$$
J_g = \left( \begin{matrix}
0 & -g^{-1} \\
g & 0
\end{matrix}\right).
$$
The induced metric on $T^*(M)$   is defined by
$$
G_g = \left( \begin{matrix}
g & 0 \\
0 &  g^{-1}
\end{matrix}\right).
$$

Now let $\G$ be a marked chord diagram, and $\sigma$ an $LM$-Morse structure on $\G$.  Recall that this consists of a metric on $\G$ and a labeling of the boundary cycles
by Morse functions $f_i : LM \to \br$.  We now assume that these Morse functions
are of the form
$$
f_i = S_{V_i} : LM \to \br
$$
where $S_{V_i}$ is the energy functional given in definition \ref{energy} using a potential function $V_i : \br/\bz \times M \to \br$.

In \cite{salamonweber} Salamon and Weber showed  how,  given such a potential function $V$, one  can define a Hamiltonian function on the cotangent bundle,
  $H_V : \br/\bz \times T^*(M) \to \br$ by the formula
\begin{equation}\label{hamilton}
H_V(t, (x,u)) = \frac{1}{2}|u|^2 + V(t, x).
\end{equation}

Using this Hamiltonian, Salamon and Weber studied the perturbed symplectic action functional on the loop space of the cotangent bundle,

 \begin{align}\label{AV}
\ca_V : L(T^*M) &\to \br \\
(\gamma, \eta) &\to \ca (\gamma, \eta) - \int_0^1H(t, (\gamma (t), \eta (t)))dt.
\end{align}
  
Here $(\gamma, \eta)$ represents a loop in $T^*M$, where  $\gamma \in LM$ and $\eta (t) \in T^*_{\gamma (t)}M$.
The classical symplectic action $\ca : L(T^*M) \to \br$ is defined by
$$
\ca ( \gamma, \eta) = \int_0^1 \langle \eta (t), \frac{d\gamma}{dt}(t) \rangle dt.
$$
 
 Following Floer's original construction, one can define a ``Floer complex", $CF^V_*(T^*M)$,   generated by the critical points of $\ca_V$, and whose boundary operator is defined by  by counting gradient flow lines of $\ca_V$.   As shown in \cite{salamonweber}, these are   
are curves $(u,v): \br \to  \ L(T^*M)$, or equivalently,
$$
(u,v) : \br \times S^1 \to T^*M
$$ that satisfy the perturbed Cauchy Riemann equations,
\begin{equation}\label{cauchyriemann}
\p_su - \nabla_t v - \nabla V_t (u) = 0   \quad \text{and} \quad \nabla_s v +\p_t u - v = 0.
\end{equation}
We refer to these maps as   holomorphic cylinders in $T^*M$ with respect to the almost complex structure $J_g$ and the Hamiltonian
$H_V$. 

Salamon and Weber also observed that  the critical points of $\ca_V$ are loops
$(\gamma, \eta)$, where $\gamma \in LM$ is a critical point of the energy functional $\cs_V : LM \to \br$, and $\eta$ is determined
by the derivative $\frac{d\gamma}{dt}$ via the metric, $\eta (v) = \langle v,\frac{d\gamma}{dt}\rangle$.
Thus the critical points of $\ca_V$ and those of $\cs_V$ are in bijective correspondence.  The following result is stated in a form proved by by Salamon and Weber in \cite{salamonweber}, but the conclusion of the theorem was first proved by Viterbo \cite{viterbo}.

\begin{thm}\label{floermorse} The Floer chain complex  $CF^{ V}_*(T^*M)$ and the Morse complex $C^V_*(LM)$
are chain homotopy equivalent.  There is a resulting isomorphism of the Floer homology of the cotangent bundle with the 
homology of the loop space,
$$
HF^{V}_*(T^*M) \cong H_*(LM).
$$
\end{thm}

This result was also proved using somewhat different methods by Abbondandolo and Schwarz \cite{abschwarz}.

The Salamon-Weber argument involved scaling the metric on $M$,  $g \to \frac{1}{\eps}g$, which  scales the almost complex
structure $J \to J_\eps$,  and the metric on $T^*M$, $G \to G_\eps =   \left( \begin{matrix}
\frac{1}{\eps}g & 0 \\
0 &  \eps g^{-1}
\end{matrix}\right).$ Notice that in this metric, the ``vertical" distance in the cotangent space is scaled by $\eps$.  

    Salamon and Weber proved that there is an $\eps_0 >0$ so that for $\eps < \eps_0$,  the set of these  holomorphic cylinders     defined with respect to the metric $G_\eps$, that  connect  critical points $(a_1, b_1)$ and  $(a_2, b_2)$ where $a_1$ and $a_2$ have relative Morse index one with respect to the action functional $\cs_V$, (or equivalently  $(a_1, b_1)$ and  $(a_2, b_2)$ have relative Conley-Zehnder index one)  is in bijective correspondence with the set of gradient 
trajectories of the energy functional $\cs_V : LM \to \br$ defined with respect to the metric $\frac{1}{\eps}g$ that connect $a_1$ to $a_2$.
Theorem \ref{floermorse} is then a consequence. 

\med
 Now again consider a marked chord diagram $\G$ with an $LM$ Morse structure $\sigma$ whose
labeling Morse functions are of the form $f_i = \cs_{V_i}$ for some potential $V_i : \br/\bz \times M \to \br$. Using the Salamon-Weber idea,  we replace the moduli space of  cylindrical graph flows, $\cm^\sigma_\Gamma (LM)$, which consists of functions $\gamma : \Sigma_G \to M$  so that the restriction to the $i^{th}$ boundary cylinder  is  a  gradient trajectories  the classical energy functional $\cs_{V_i}$, by the moduli space of ``cylindrical holomorphic curves" in the cotangent bundle $T^*(M)$, $\cm^{hol}_{(\G, \sigma, \eps)}(T^*M)$,  which consists of  continuous maps,
$$
   \phi : \Sigma_\G \to T^*(M)
   $$
 such that the restrictions to the cylinders, 
$$
\phi_i : (-\infty, 0] \times S^1_{c_i} \to T^*M \quad \text{and} \quad \phi_j : [0, +\infty) \times S^1_{c_j} \to T^*M
$$
are holomorphic with respect to the almost complex structure $J_{\eps}$ and the Hamiltonians $H_{V_i}$ and $H_{V_j}$ respectively.
Here the circles $S^1_{c}$ are round with circumference $c_j$ determined by the metric given by the structure $\sigma$.

\med
Like in the last section we have restriction  maps (compare (\ref{restrict2}))
 \begin{equation}\label{restrict3}
 (L(T^*M))^q \xleftarrow{\rho_{out}} \cm^{hol}_{(\G, \sigma, \eps)}(T^*M) \xr{\rho_{in}} (L(T^*M))^p .
  \end{equation}
 $\rho_{in}$ is  defined
by  sending a cylindrical flow $\phi$ to  $\prod_{i=1}^p \phi_{i,-1} : \{-1\} \times S^1  \to T^*M$ and $\rho_{out}$ sends $\phi$ to $\prod_{j=p+1}^{p+q} \phi_{j,1}:  \{1\} \times  S^1\to T^*M.$

Motivated by the umkehr maps $ (\rho_{in})_! : h_*((LM)^p) \to h_{*+ \chi (\G) \cdot d}(\cm^\sigma_\G(LM))$ given in the last section, as well as the Salamon-Weber results, 
we conjecture the following analogue of the existence of the string topology operations, and their field theoretic properties.  

\begin{conj}\label{basic} For every marked chord diagram $\G$, there is an umkehr map
$$
(\rho_{in})_! : (HF_*(T^*M))^{\otimes p} \to H_{*+\chi (\G)\cdot d}( \cm^{hol}_{(\G, \sigma, \eps)}(T^*M))
$$
and a homomorphism
$$
(\rho_{out})_* :  H_{* }( \cm^{hol}_{(\G, \sigma, \eps)}(T^*M)) \to (HF_*(T^*M))^{\otimes q}
$$
so that the operations
$$
\theta_\G  = (\rho_{out})_* \circ (\rho_{in})_! :  (HF_*(T^*M))^{\otimes p}  \to (HF_*(T^*M))^{\otimes q}
$$
satisfy the following properties:
\begin{enumerate}
\item The maps $\theta$ fit together to define a positive boundary, topological field theory.
\item With respect to the Salamon-Weber isomorphism $HF_*(T^*M) \cong H_*(LM)$ (theorem \ref{floermorse})
the Floer theory operations $\theta_\G$ equal the string topology operations $q_\G$ studied in the last two sections.
\end{enumerate}
\end{conj}

\bfl
\bf Remark. \rm
   The existence of a field theory structure on the Floer homology of a closed symplectic manifold was established by Lalonde \cite{lalonde}.  The above conjecture should be directly related to Lalonde's constructions.
   
   \efl
    
 \med
 Now one might also take the more geometric approach to the construction of these Floer theoretic operations, analogous to Ramirez's  geometrically defined Morse  theoretic constructions of string topology operations.  This would involve the study of the space of cylindrical holomorphic curves in $T^*M$, with boundary conditions   in stable
 and unstable manifolds of critical points, $ \cm^{hol}_{(\G, \sigma, \eps)}(T^*M, \vec{a})$. Smoothness and compactness properties need to be established for these moduli spaces.  In particular, in a generic situation  their dimensions should
 be given by the formula
 $$
 dim \,( \cm^{hol}_{(\G, \sigma, \eps)}(T^*M, \vec{a})) = \sum_{i=1}^p Ind (a_i)  - \sum_{j=p+1}^{p+q}Ind (a_j)  + \chi (\G)\cdot d
 $$
 where $Ind (a_i)$ denotes the Conley-Zehnder index. 
 
 \med
 We remark that in the case of the figure 8,  this analysis has all been worked out by Abbondandolo and Schwarz \cite{schwarzmontreal}.  In this case $\Sigma_\G$ is a Riemann surface structure on the pair of pants.  They proved the existence of a ``pair of pants" algebra structure on $HF_*^V(LM)$ and with respect to their isomorphism,
 $HF_*^V(LM) \cong H_*(LM)$ it is isomorphic to the pair of pants product on the Morse homology of $LM$.  In view
 of the comment following definition \ref{geodef} we have the following consequence.
 
 \med
 \begin{thm} With respect to the isomorphism $HF^V_*(T^*M) \cong H_*(LM)$, the pair of pants product in the Floer homology of the cotangent bundle corresponds to the Chas-Sullivan string topology product.
 \end{thm}

 \med
 Another aspect of the relationship between the symplectic structure of the cotangent bundle and the string topology
 of the manifold,  has to do with the relationship between the moduli space of $J$-holomorphic curves with cylindrical
 boundaries, $\cm_{g,n}(T^*M)$, and moduli space of cylindrical holomorphic curves, $ \cm^{hol}_{(\G, \sigma, \eps)}(T^*M)$, where we now let $\G$ and $\sigma$ vary over the appropriate space of metric graphs.  We conjecture that these moduli spaces are related as a parameterized version of the relationship between the moduli space of Riemann surfaces and the space of metric fat graphs (theorem \ref{penner:CMP}). 
 
 Once established, this conjecture would give a direct relationship between Gromov-Witten invariants of the cotangent bundle, and the string topology of  the underlying manifold. In this setting the Gromov-Witten invariants would be  defined using moduli spaces of curves with cylindrical ends rather than marked points, so that the invariants would be defined in terms of the homology of the loop space (or, equivalently, the Floer homology of the cotangent bundle), rather than the homology of the manifold. 
 
  We believe that there is a very deep  relationship between the symplectic topology of the cotangent bundle and the string topology of the underlying manifold.  There is much work yet to be done in understanding the extent of this relationship.

% Version of 03/17/05

\chapter{Brane topology}
%\label{branetopology}

\section{The higher-dimensional cacti operad}

String topology may be generalized to higher-dimensional sphere spaces
$SM: = M^{S^n} = \Map (S^n, M)$ for $n \ge 1$, see
\cite{sullivan-voronov}. See also the paper \cite{chataur} by Chataur,
in which a string (i.e., ``dot'') product on the homology of sphere
spaces $M^{S^n}$ was defined, as well as the other string topology
operations for $n=3$. In the paper \cite{kallel-salvatore} by Kallel
and Salvatore, the string product for $n=2$ and $M = \nc P^N$ was
shown to be coming from the holomorphic mapping space
$\operatorname{Hol} (\nc P^1, \nc P^N)$ via Segal's homological
approximation theorem of the continuous maps by holomorphic ones.

The corresponding cacti operad does not admit a nice conbinatorial
description available for $n=1$, so that we will consider an awfully
big, but neat, infinite dimensional operad, which will do the job.

\begin{figure}
\centerline{\includegraphics[width=1.2in]{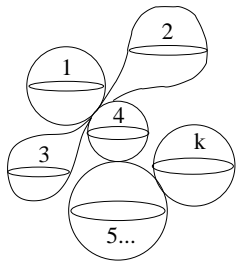}}
\caption{A 2d-cactus.}
\label{dcactus}
\end{figure}

For $n \ge 1$, the \emph{$n$-dimensional cacti operad} is an operad of
topological spaces. It may be described as a collection of topological
spaces $\CC^n(k)$, $k \ge 1$, defined as follows. The space $\CC^n(k)$
is the space of all continuous maps from the unit $n$-sphere $S^n$ to
the union of $k$ labeled (by the numbers $1, \dots, k$) $n$-spheres,
called the \emph{lobes}, joined at a few points, such that every two
lobes intersect at most at one point and the \emph{dual graph} (whose
vertices are the lobes and the intersection points and whose edges
connect the lobes with the adjacent intersection points) of this union
is a connected tree. One can think of a point of the space $\CC^n(k)$
as a pair $(c,\phi)$, where $c$ is a join of $k$ spheres, as above,
called a \emph{cactus}, see Figure~\ref{dcactus}, and $\phi: S^n \to
c$ a continuous map, called a \emph{structure map} from the unit
$n$-sphere to the cactus.

The topology on $\CC^n(k)$ is given as follows. The topology on the
set of cacti $c$ is induced from the following inclusion into
$(S^n)^{k(k-1)}$: for each ordered pair $(i,j)$, $1 \le i \ne j \le
k$, one takes the point on the $i$th lobe of a given cactus at which
this lobe is attached to a lobe eventually leading to the $j$th
lobe. The topology on the ``\emph{universal cactus}'' $\UU$ defined as
the set of pairs $(c,x)$, where $c$ is a cactus and $x$ is a point on
it, is induced from the inclusion of $\UU$ into $(S^n)^{k^2} =
(S^n)^{k(k-1)} \times (S^n)^k$, the second factor governing the
position of $x$ on $c$ from the point of view of each lobe. The
universal cactus projects naturally onto the \emph{space of cacti},
the space of tree-like joins of $k$ labeled spheres, with the
``fiber'' being the corresponding cactus. Finally, a basis of the
topology on the set $\CC^n(k)$, which may be identified with the set
of continuous maps from the unit sphere $S^n$ to the fibers of the
universal cactus, is formed by finite intersections of the sets of
such maps $\Phi: S^n \to \UU$ satisfying $\Phi(K) \subset U$ for given
compact $K \subset S^n$ and open $U \subset \UU$.

An operad structure on the collection $\CC^n = \{ \CC^n(k) \; | \; k
\ge 1 \}$ is defined as follows. An action of the symmetric group is
defined via changing the labels of the lobes. A unit element $\id \in
\CC^n(1)$ is the identity map $\id : S^n \to S^n$. Operad compositions
$\circ_i : \CC^n(k) \times \CC^n(l) \to \CC^n(k+l-1)$, $i = 1, \dots,
k$, are defined for two cacti $C_1$ with $k$ lobes and $C_2$ with $l$
lobes as follows. First, construct a new cactus by attaching $C_1$ to
$C_2$ via a certain attaching map from the $i$th lobe of $C_1$ to the
cactus $C_2$. This attaching map is nothing but the structure map
$f_2$ from the unit sphere to the cactus $C_2$. At the level of dual
graphs, we glue in the tree corresponding to $C_2$ in place of the
vertex corresponding to the $i$th lobe of $C_1$ by connecting the
incoming edges into this vertex in $C_1$ with certain vertices of the
tree of $C_2$, as prescribed by the map $f_2$. This procedure, similar
to ones appearing in \cite{connes-kreimer,kon-soib}, results in a tree
again at the level of dual graphs and gives a cactus $C_1 \circ_i C_2
:= C_2 \cup_{f_2} C_1$ with $k+l-1$ lobes. Label the lobes by
inserting the labels $1, \dots, l$ of $C_2$ at the $i$th space in the
set of labels $1, \dots, k$ of $C_1$ and relabeling the resulting
linearly ordered set as $1, \dots, k+l-1$. Then, define a structure
map for this cactus as the composition of the structure map $f_1: S^n
\to C_1$ and the natural map $C_1 \to C_2 \cup_{f_2} C_1$.

This operad is a close relative of the \emph{framed little
$(n+1)$-disks operad} $f\DD^{n+1}$, see Section~\ref{GBV}. Little is
known about a direct relation, except the following theorem and
Remark~\ref{Salvatore} of Salvatore, which one might try to generalize
to higher-dimensional cacti.

\begin{thm}[Sullivan-AV \cite{sullivan-voronov}]
\label{hard}
There exists an operad morphism
\[
H_* (f\DD^{n+1}; \nq) \to H_* (\CC^n ; \nq)
\]
from the rational homology framed little $(n+1)$-disks operad to the
rational homology $n$-cacti operad.
\end{thm}

\begin{proof}[Idea of proof]
One constructs maps $f\DD^n (k) \to \CC^n(k)$, $k = 1, 2, 3$, from the
first three components of the little disks operad to the cacti
operad. To do that, first construct $\SO(n+1)$-equivariant maps
$\DD^{n+1}(k) \to \CC^n(k)$, $k = 1, 2, 3$, for the (nonframed) little
disks operad $\DD^{n+1}$ and then extend them to $f\DD^{n+1}(k)$ by
equivariance. The point is that for the framed little disks operad,
the first three components contain all the information about its
homology operad: the generators and relations of $H_*
(f\DD^{n+1}; \nq)$ lie in $H_* (f\DD^{n+1}(k); \nq)$ for $k \le
3$. This follows from the fact, noticed by Salvatore and Wahl
\cite{salvatore-wahl}, that the operad $H_* (f\DD^{n+1}; \nq)$
is a semidirect product of the homology little disks operad $H_*
(\DD^{n+1}; \nq)$ and the rational homology of the group $\SO
(n+1)$. They make an explicit computation of the rational homology
$H_* (f\DD^{n+1}; \nq)$, see Theorem~\ref{sw:BV}, which we will use
later. The proof is completed by showing that the constructed operad
maps respect the operad structures up to homotopy. This implies that
there is a homology operad morphism $H_* (f\DD^{n+1}; \nq) \to
H_* (\CC^n; \nq)$.
\end{proof}

A complete proof along these lines is quite long and ad hoc, as it
involves a number of explicit constructions with three disks in the
unit disk and explicit homotopies.

In principle, one may give a more direct and, perhaps, shorter proof
of this theorem over the rationals, by mapping the generators,
explicit in Salvatore-Wahl's theorem, of the operad $H_*
(f\DD^{n+1}; \nq)$ to specific cycles in $\CC^n$ and checking that the
relations are satisfied by showing that the corresponding cycles are
homologous in $\CC^n$ explicitly. However, the proof outlined earlier
seems to be more conceptual, avoids computing the homology of the
cacti operad, which might make a good research problem, and hints on
the possibility of an operad morphism $f\DD^{n+1} \to \CC^n$.

It would also be interesting to study the relationship between the
$k$th connected component $\CC^n(k)_0$ of the cacti operad and the
classifying space of the group of homeomorphisms (or diffeomorphisms)
of $S^{n+1}$ with $k+1$ disks removed. In the $n=1$ case, all these
spaces are known to be homotopy equivalent to the framed little disks
space $f\DD^2 (k+1)$ and the moduli space of $(k+1)$-punctured $\nc
P^1$'s with a real tangent direction at each puncture.

\section{The cacti action on the sphere space}

Let $M$ be an oriented manifold of dimension $d$. We would like to
study continuous $k$-ary operations on the (\emph{free}) \emph{sphere
space} $SM: = M^{S^n} := \Map (S^n, M)$, functorial with respect to
$M$. By passing to singular chains or homology, these operations will
induce functorial operations on the chain and homology level,
respectively. We will generalize the motivic approach described in
Section~\ref{cacti:action} in the $n=1$ case, following
\cite{sullivan-voronov}. We will use the notation $\CC := \CC^n$
throughout this section.

Like for $n=1$, consider the following diagram
\begin{equation}
\label{rho:n}
\CC(k) \times (SM)^k \xleftarrow{\rho_{\varin}} \CC (k)M
\xrightarrow{\rho_{\out}} SM
\end{equation}
for each $k \ge 1$, where $\CC (k)M$ is the space of triples
$(c,\phi,f)$ with $(c,\phi) \in \CC(k)$ and $f: c \to M$ a continuous
map from the corresponding cactus to $M$. The map $\rho_{\varin}: \CC
(k)M \to \CC(k) \times (SM)^k$ takes a triple $(c,\phi,f)$ to
$(c,\phi) \in \CC(k)$ and the restrictions of $f$ to the $k$ lobes of
$c$. It is an embedding of codimension $d(k-1)$. The map $\rho_{\out}:
\CC (k)M \to SM$ takes $(c,\phi,f)$ to the composition of the
structure map $\phi: S^n \to c$ with $f: c \to M$. The above diagram
defines an operad action in the category $\Corr$ of correspondences,
see Section~\ref{cacti:action}, on the sphere space $SM$. Then the
proof of the following theorem is no different from the proof of
Theorem~\ref{homotopyaction} for $n=1$.

\begin{thm}[\cite{sullivan-voronov}]
\label{homotopyaction:n}
\begin{enumerate}
\item
Diagram \eqref{rho:n}, considered as a morphism $\CC(k) \times (SM)^k
\to SM$ in $\Corr$, defines the structure of a $\CC$-algebra on the
sphere space $SM$ in $\Corr$.
\item
This $\CC$-algebra structure on the sphere space $SM$ in $\Corr$
induces an $h_* (\CC)$-algebra structure on the shifted homology
$h_{*+d} (SM)$ for any multiplicative generalized homology theory
$h_*$ which supports an orientation of $M$.
\end{enumerate}
\end{thm}

\section{The algebraic structure on homology}

Combining Theorems \ref{hard} and \ref{homotopyaction:n} we observe the
following result.
\begin{crl}
\label{disks}
The shifted rational homology $H_{*+d} (SM; \nq):= H_* (SM; \nq)[d]$
of the sphere space $SM = M^{S^n}$ is an algebra over the operad $H_*
(f\DD^{n+1}; \nq)$.
\end{crl}

Taking into account the explicit generators and relations for the
operad $H_* (f\DD^{n+1}; \nq)$ given by Theorem~\ref{sw:BV}, we obtain
the following algebraic structure on the rational homology of the
sphere space $SM$.
\begin{crl}
\label{algebraic}
  The \textup{(}shifted\textup{)} rational homology $H_{*+d} (SM;
  \nq)$ admits the algebraic structure of Theorem~\ref{sw:BV}. This
  includes operations $\cdot$, $[,]$, $\Beta_1, \dots,$ and, for $n$
  odd, $\Delta$, satisfying the corresponding identities.
\end{crl}

\begin{rem}
  The dot product, like in the $n=1$ case, see \cite{cohen-jones},
  comes from a \emph{ring} spectrum structure
\[
(SM)^{-TM} \wedge (SM)^{-TM}  \to (SM)^{-TM},
\]
which is constructed using the a standard pinch map $S^n \to S^n
\vee S^n$ as a point in $\CC^n(2)$ and the diagram
\[
SM \times SM \hookleftarrow SM \times_M SM \to SM,
\]
which is a specification of \eqref{rho:n} at that point of $\CC^n (2)$
for $k=2$.
\end{rem}

In reality, though, all the unary operations but $\Delta$ in the odd
$n$ case and $B_{n/2}$ in the even $n$ case vanish on $H_{*+d} (SM;
\nq)$, as we will see soon. So, perhaps, an analogue of the previous
corollary with integral or finite coefficients would bring a lot more
information. What we have over $\nq$ can be described using the
following definition.

\begin{definition}
  A BV$_{n+1}$-\emph{algebra} over a field of zero characteristic is a
  graded vector space $V$, along with a dot product $ab$, a bracket
  $[a,b]$, and a BV operator $\Delta$ of degree $n$, if $n$ is odd,
  and an operator $\Beta = \Beta_{n/2}$ of degree $2n-1$, if $n$ is
  even, satisfying the properties (1)--(3) and (6)--(7) of
  Theorem~\ref{sw:BV}.
\end{definition}

\begin{thm}[\cite{sullivan-voronov}]
  The shifted rational homology $H_{*+d}(SM; \nq)$ of the sphere space
  $SM = M^{S^n}$ of an oriented manifold $M$ of dimension $d$ has a
  natural structure of a \textup{BV}$_{n+1}$-algebra.
\end{thm}

This theorem follows trivially from Corollary~\ref{algebraic}; in
other words, in view of Theorem~\ref{sw:BV}, any $H_* (f\DD;
\nq)$-algebra is a BV$_{n+1}$-algebra. Moreover, one has the following
vanishing result.

\begin{prop}[\cite{sullivan-voronov}]
  The operators $B_i$, $i \ne n/2$, vanish on $H_{*+d} (SM; \nq)$.
\end{prop}

\begin{proof}
  These operators come from an action of the group $\SO(n+1)$, which
  is obviously homotopy equivalent to the monoid $f\DD^{n+1}(1)$. This
  monoid maps naturally to $\CC^n(1) = \Map (S^n, S^n)$ and acts
  through it on the sphere space $SM$. This action is the $k=1$ part
  of the action from Theorem~\ref{homotopyaction:n}. Thus, it suffices
  to perform the following computation of the rational homology of
  $\Map (S^n, S^n)$.

\begin{lm}
\begin{enumerate}
\item We have an isomorphism
\[
H_* (\Map(S^n, S^n); \nq) =
\begin{cases}
\nq[\Delta, q, q^{-1}] & \text{for $n$ odd},\\
\nq[\Beta, q, q^{-1}] & \text{for $n$ even}.\
\end{cases}
\]
bIn either case the right-hand side is a graded commutative algebra on
one generator $q$ of degree zero and one odd generator $\Delta$ of
degree $n$ for $n$ odd or one odd generator $\Beta$ of degree $2n-1$
for $n$ even, localized at $q=0$.
\item
Under the map $\SO(n+1) \to \Map (S^n, S^n)$, the elements $\Delta$
and $\Beta_{n/2} \in H_* (\SO(n+1); \nq)$ map to the above
$\Delta$ and $\Beta$, while the other $B_i$'s map to zero.
\end{enumerate}
\end{lm}
\begin{proof}[Idea of proof of Lemma, after Sullivan]
The proof uses a method of rational homotopy theory, which finds the
(commutative, minimal) DGA representing the structure group of a fiber
bundle, given a DGA representing a generic base and a DGA representing
the fiber. We will compute the homology of $\Map (S^4, S^4)$ as an
example. Here we are talking about an $S^4$-bundle over a base
represented by a DGA $A$. The DGA of the fiber is $\nq[u,v]$ with
degrees $\abs{u} = 4$ and $\abs{v}=7$ and differential $du=0$ and $dv
= u^2$. Then the DGA model of the total space may be obtained as
$A[u,v]$ with a differential determined by $du = a$ for some $a \in A$
and $dv = u^2 + pu + q$ for some $p,q \in A$. Of course, we should
have $d^2 = 0$, so $d^2 v =0$ yields $2ua + pa + u dp + dq = 0$ or $2a
+ dp =0$ and $pa + dq = 0$. If we set $u' := u - p/2$, we will have
$du' = 0$ and $dv = (u')^2 + q - p^2/4$. Thus, the total space DGA is
now universally presented as $A[u',v]$ with the above
differential. Note that for $q' := q - p^2/4$ we have $dq' =
0$. Therefore, an $S^4$-bundle is determined in rational homotopy
theory by a characteristic class $q' \in A$ of degree eight and the
classifying space must rationally be $K(\nq,8)$. By transgression, the
structure group must rationally be $K(\nq,7) \sim_\nq S^7$. One needs
some extra work to see that a homology generator of this space comes
from the generator $\Beta$ of $H_7 (\SO(5); \nq)$.
\end{proof}
\end{proof}

\section{Sphere spaces and Hochschild homology}

Here we will announce a result of the second author, which provides an
approximation of the sphere space $SX = X^{S^n}$ via configuration
spaces with labels and implies a computation of the homology of $SX$
as the Hochschild homology of the $n$-algebra $C_* (\Omega^n X)$,
where $\Omega^n X = \Map_* (S^n, X)$ is the based $n$-fold loop space
of a space $X$. These results directly generalize the results of
Burghelea and Fiedorowicz \cite{bf} for loop spaces ($n=1$). The first
computation of the homology of $SX$ as an $n$-algebra Hochschild
homology was done by Po Hu \cite{pohu}, who generalized the
Chen-Getzler-Jones-Petrack-Segal method. Unfortunately, we can only
speculate that our version of the $n$-algebra Hochschild homology is
isomorphic to the one used by Po Hu.

\begin{thm}
The topological Hochschild complex (see below) of the based $n$-fold
loop space $\Omega^n X$, considered as an $f\DD^n$-algebra is
homeomorphic to the space $X^{S^n}_0$ of maps $S^n \to X$ passing
through the basepoint of a space $X$.
\end{thm}

\begin{proof}
Let us first define the topological Hochschild complex; when this is
done, the proof will be fairly straightforward. 

From this point on, we need to consider a homotopy equivalent version
of the framed little $n$-disks operad, which is the operad of labeled,
framed little disks inside the standard (unit) upper hemisphere of
$S^n$ with a chosen orientation. A little disk on the sphere is
understood as a ball in the standard metric on the sphere induced from
the ambient Euclidean $\nr^{n+1}$. A frame in a little disk is a
positively oriented orthonormal frame in the tangent space to the
sphere at the center of the disk. The operad composition is given by
shrinking the hemisphere along the great arcs towards the north pole
to fit the prescribed little disk on another unit hemisphere,
transporting the shrunk disk to the place of that prescribed little
disk and rotating the shrunk disk to fit the frames. Here we assume
that we choose a standard orthonormal frame at the north pole of the
standard sphere. We will still denote the resulting operad $f\DD^n$.

Consider a similar space $f\DD S^n = \{ f\DD S^n (k) \; | \; k \ge
1\}$ of labeled, framed little disks on the whole unit sphere
$S^n$. It is a (right) module over the operad $f\DD^n$ via similar
combination of dilations, translations, and rotations, as in the
definition of the operad structure on $f\DD^n$ above.

Suppose a space $A$ is an algebra over the operad $\OO$ in the
category of topological spaces and $M$ an $\OO$-module. Then we can
form the \emph{tensor product of $M$ and $A$ over $\OO$} as follows:
\[
M \otimes_\OO A := \coprod_k M(k) \times_{\Si_k} A^k /\sim,
\]
where $\sim$ denotes the equivalence relation which may be described
roughly as $(mo, a) \sim (m,oa)$ for $m \in M$, $o \in \OO$, and $a
\in A$. In other words, $M \otimes_\OO A$ is the coequalizer of the
diagram $\xymatrix{M \times_\Si \OO \times_\Si \times A \ar@<.7ex>[r]
\ar@<-.3ex>[r] & M \times_\Si A}$.

Now let $A$ be an algebra over the operad $f\DD^n$, for example, the
based $n$-fold loop space $\Omega^n X$ of a pointed space $X$, and
$f\DD S^n$ the $f\DD^n$-module of framed little disks on the
sphere. Then, by definition, the \emph{topological Hochschild homology
of the $f\DD^n$-algebra $A$} is the tensor product $f\DD S^n
\otimes_{f\DD^n} A$.

To prove the theorem, we thereby need to present a homeomorphism
\[
f\DD S^n \otimes_{f\DD^n} \Omega^n X \cong X^{S^n}_0.
\]
We define it by assigning to a configuration of $k$ framed little
disks on $S^n$ and a given collection of $k$ maps of the upper
hemisphere in $S^n$ to $X$ (representing $k$ points in $\Omega^n X$)
to a maps from the sphere $S^n$ to $X$ by using the $k$ maps on the
little disks (prepended by the appropriate translations, dilations,
and rotations on the sphere, as above) and extending these maps by a
constant to a map $S^n \to X$. This gives a point in $X^{S^n}_0$.

The inverse of this homeomorphism is given by taking, for a given map
$S^n \to X$ passing through the basepoint, a point on $s \in S^n$
mapping to the basepoint and thinking of its complement $S^n \setminus
\{s\}$ as one little disk on $S^n$, giving it any framing, which
determines a map from the standard upper hemisphere to $S^n \setminus
\{s\}$ mapping the equator to $s$, and thus getting an element of
$\Omega^n X$. This combines into a map $X^{S^n}_0 \to f\DD S^n (1)
\times \Omega^n X$, which we map into the tensor product $f\DD S^n
\otimes_{f\DD^n} \Omega^n X$ naturally.
\end{proof}

The relation between the free sphere space $X^{S^n}$ and the
``semifree'' one $X^{S^n}_0$ is amazingly simple.
\begin{prop}
For a path connected space $X$ having numerable category, in
particular, a connected CW complex, the inclusion $X^{S^n}_0 \subset
X^{S^n}$ is a homotopy equivalence.
\end{prop}
\begin{proof}
Consider the map $p: X^{S^n} \to X$ evaluating a map $S^n \to X$ at
the north pole. One checks that it is a fibration by verifying the
homotopy lifting property directly. There is a similar fibration
$X^{S^n}_0 \to X$ and a commutative square
\[
\begin{CD}
X^{S^n}_0 @>>> X^{S^n}\\
@VVV @VVV\\
X @>\id>> X.
\end{CD}
\]
It suffices to show that the map induced on the fiber is a homotopy
equivalence. Indeed the map induced between the fibers over the
basepoint is just the identity map.
\end{proof}

% If we use a cofibrant version $Q f \DD^n$ of the operad $f\DD^n$ and a
% cofibrant version $Q f \DD S^n$ of the module $f\DD S^n$ (e.g., given
% by the real version of the Fulton-MacPherson compactification of the
% configuration spaces), 
We can pass to singular cochains $C_*$ and homology $H_*$ and obtain
the following corollary.
\begin{crl}
\begin{enumerate}
\item
The Hochschild chain complex $C_*^{(n)} (C_* \Omega^n X, C_* \Omega^n
X)$ of the $C_* f\DD^n$-algebra $C_* \Omega^n X$ is homotopy
equivalent to the singular cochain complex $C_* X^{S^n}_0$.
\item The Hochschild homology $H_*^{(n)} (C_* \Omega^n X, C_* \Omega^n
X)$ of the $C_* f\DD^n$-algebra $C_* \Omega^n X$ is isomorphic to the
homology $H_* (X^{S^n}_0; \nz)$.
\end{enumerate}
If $X$ is path connected and has numerable category, one can replace
$X^{S^n}_0$ by the sphere space $X^{S^n}$ in the above.
\end{crl}
Here the \emph{Hochschild chain complex} $C_*^{(n)} (A,A)$ of a $C_* f
\DD^n$-algebra $A$ is understood as $C_* f\DD S^n \otimes_{C_* f
\DD^n} A$ and the \emph{Hochschild homology} $H_*^{(n)} (A,A)$ as the
homology of this complex. Motivation for this definition may be the
$n=1$ case, in which $C_* f\DD S^1$ and $C_* f \DD^1$ must be replaced
by the cofibrant models $C_* \overline{C}(S^1)$ and $C_*
\overline{C}(D^1)$, where $\overline{C}$ stands for the
Fulton-MacPherson compactification and $C_*$ for cellular
chains. After these changes, the tensor product becomes a
configuration space with summable labels, see \cite{bf,salvatore}, and
we get the usual Hochschild complex.

As we mentioned earlier, the first result of this type was obtained by
Po Hu \cite{pohu}, who used a different notion of the Hochschild
homology. It will be very interesting to find a relationship between
them. Hopefully, our result is a geometric incarnation of that of Po
Hu.

\subsection{Brane topology and Kontsevich's Hochschild cohomology conjecture}

We would like to end the discussion with the following table of
analogies between algebra and topology. Since the discussion is highly
speculative, we will freely confuse the notion of an $n$-algebra,
which is an algebra over the operad $H_*(\DD^n; \nz)$, see
Section~\ref{n-algebras}, with that of an algebra over the chain
operad $C_* (\DD^n; \nz)$ or the topological operad $\DD^n$ and an
$n_\infty$-algebra. Recall that a 1-algebra is the same as an
associative algebra and 2-algebra is the same as a G-algebra.

\begin{center}
\begin{tabular}{|p{2.2in}|p{2.3in}|}
\hline
\multicolumn{1}{|c|}{Algebra} & \multicolumn{1}{c|}{Topology}\\
\hline
a 1-algebra $A$ & a loop space $\Omega M$ is a 1-algebra\\
\hline
$H^*(A,A)$ is a 2-algebra by Gerstenhaber &
$H_{*+d} (LM)$ is a 2-algebra
by Chas-Sullivan\\
\cline{2-2}
& $H_{*+d} (C_* \Omega M, C_* \Omega M)$, a graded abelian group
isomorphic to $H_{*+d} (LM)$ by Burghelea-Fiedorowicz, is a 2-algebra,
via Poincar\'e duality and the fact that $C_* \Omega M$
is a 1-algebra\\
\hline
an $n$-algebra $A$ & a loop space $\Omega^n M$ is an $n$-algebra\\
\hline
$H^*_{(n)} (A,A)$ is an $(n+1)$-algebra by Kontsevich &
$H_{*+d} (M^{S^n})$ is an $(n+1)$-algebra
by Sullivan-AV\\
\cline{2-2}
& $H^{(n)}_{*+d} (C_* \Omega^n M, C_* \Omega^n M)$, a graded abelian
group isomorphic to $H_{*+d} (M^{S^n})$ by Po Hu, is an
$(n+1)$-algebra, via Poincar\'e duality and the fact that $C_*
\Omega^n M$ is an $n$-algebra\\
\hline
\end{tabular}
\end{center}
It would be very interesting to find a relation of our notion of
Hochschild $n$-algebra homology to an algebraic one, see, for example,
\cite{kon:after,tamarkin:koncon,hkv}, and check if
the isomorphisms in the right column of the table respected the
$n$-algebra structures. One may think of such a statement as a
topological incarnation of Kontsevich's conjecture on Hochschild
cohomology, proven in \cite{tamarkin:koncon,hkv}: the Hochschild
cochain complex of an $n$-algebra is naturally an
$(n+1)_\infty$-algebra.

\backmatter
%\bibliographystyle{amsalpha}
%\bibliography{ralphbiblio,../dennis/op,../grant/nsf}

\providecommand{\bysame}{\leavevmode\hbox to3em{\hrulefill}\thinspace}
\providecommand{\MR}{\relax\ifhmode\unskip\space\fi MR }
% \MRhref is called by the amsart/book/proc definition of \MR.
\providecommand{\MRhref}[2]{%
  \href{http://www.ams.org/mathscinet-getitem?mr=#1}{#2}
}
\providecommand{\href}[2]{#2}

\end{document}